\documentclass[reqno,12pt]{amsart}
\usepackage{amsmath,amssymb,amsfonts,amscd}

\usepackage{graphicx,color}
\input xy
\xyoption{all}

{
\newtheorem{thm}{Theorem}[section]

\newtheorem{cor}[thm]{Corollary}
\newtheorem{lem}[thm]{Lemma}
\newtheorem{prop}[thm]{Proposition}

\theoremstyle{definition}

\newtheorem{defn}[thm]{Definition}

\newtheorem{rem}[thm]{Remark}

\newtheorem*{pf}{Proof}
}

\numberwithin{equation}{section}
\newcommand{\bp}{\begin{pmatrix}}
\newcommand{\ep}{\end{pmatrix}}
\newcommand{\bps}{\begin{smallmatrix}}
\newcommand{\eps}{\end{smallmatrix}}
\def\C{{\mathbb C}}

\def\P{{\mathbb P}}

\def\R{{\mathbb R}}
\def\Z{{\mathbb Z}}
\def\A{{\mathcal A}}

\def\E{{\mathcal E}}

\def\S{{\mathcal S}}
\def\T{{\mathcal T}}

\def \fm{{\mathfrak m}}
\def \sp{{\frak{sp}}}

\def\p{{\partial }}

\def\bb{{\bar b}}

\def\sb{{\bar s}}

\def\qb{{\bar q}}

\def \0{{\bf 0}}
\def \1{{\bf 1}}

\def \la{\langle}
\def \ra{\rangle}

\def \Der{\mathrm{Der}}
\def\Tr{\mathrm{Tr}}
\def\Hom{\mathrm{Hom}}
\def\hom{\mathrm{hom}}

\def\Ext{\mathrm{Ext}}
\def \Mat{\mathrm{Mat}}
\def \Ker{\mathrm{Ker}}

\def \Ob{\mathrm{Ob}}
\def \diag{\mathrm{diag}}
\def \rank{\mathrm{rank}}
\def \Id{\mathrm{Id}}

\def \Homub{\underline{\Hom}}

\def \CMg{\mathrm{CM}^{gr}}
\def \CMgb{\underline{{\rm CM}^{gr}}}
\def \syz{\mathrm{syz}}

\def \quiver{\vec{\Delta}_W}
\def \quiverT{\vec{\Delta}_W^T}
\def \quiver'{\vec{\Delta}_W'}

\def \PAT{\C\quiverT/I_{W,\Lambda}}
\def \PA'{\C\quiver'/I'_{W,\Lambda}}

\def \mf#1#2#3#4{
\xymatrix{
{#1}\  \ar@<0.4ex>[r]^{{#2}} & \ {#4}
\ar@<0.4ex>[l]^{{#3}}
}
}

\def \mfs#1#2#3#4{\!
\xymatrix@C=1,5em{{#1} \! \ar@<0.2ex>[r]^{{#2}} & \! {#4}
\ar@<0.2ex>[l]^{{#3}}
}
\!}

\def \mfl#1#2#3#4{
\xymatrix@C=2.6em{{#1}\  \ar@<0.4ex>[r]^{{#2}} &\  {#4}
\ar@<0.2ex>[l]^{{#3}}
}
}

\def \mfss#1#2#3#4{\!
\xymatrix@C=1.5em{{#1} \ar@<0.3ex>[r]^{{#2}} & {#4}
\ar@<0.3ex>[l]^{{#3}}
}
\!}

\def \olF{\overline{F}}
\def \olFF{{\overline{F'}}}

\def \tMF#1#2{\mathit{MF^{gr}_{#1}(#2)}}

\def \tHMF#1#2{{HMF}^{gr}_{#1}({#2})}

\def \ov#1{\frac{1}{#1}}

\def \Vt{V_1}
\def \Vb{V_{\bar 1}}
\def \vt{v_1}
\def \vb{v_{\bar 1}}

\newcommand{\gr}{\ensuremath{\mathrm{gr\text{--}}}}
\newcommand{\Dbsing}{\ensuremath{D^{\mathrm{gr}}_{\mathrm{Sg}}}}

\begin{document}
\begin{flushleft}
\hfill RIMS-1600\, \\
\end{flushleft}

\vspace*{0.4cm}

\title[Triangulated categories of matrix factorizations for 
RSW with $\varepsilon=-1$]
{Triangulated categories of matrix factorizations for 
regular systems of weights with $\varepsilon=-1$}

\date{August 1, 2007.}
\author{Hiroshige Kajiura}
\address{Research Institute for Mathematical Sciences, Kyoto University, 
606-8502, Japan}
\email{
kajiura@kurims.kyoto-u.ac.jp}
\author{Kyoji Saito}
\address{Research Institute for Mathematical Sciences, Kyoto University,
606-8502, Japan}
\email{saito@kurims.kyoto-u.ac.jp}

\author{Atsushi Takahashi}
\address{Department of Mathematics, Graduate School of Science, 
Osaka University, 
Toyonaka Osaka, 
560-0043, Japan}
\email{takahashi@math.sci.osaka-u.ac.jp}
\begin{abstract}
We construct a full strongly exceptional collection 
in the triangulated category of graded matrix factorizations of 
a polynomial associated to a non-degenerate 
regular system of weights whose smallest exponents are equal to $-1$. 
In the associated Grothendieck group, 
the strongly exceptional collection defines 
a root basis of a generalized root system of sign $(l,0,2)$ and a 
Coxeter element of finite order, whose primitive eigenvector is
a regular element in the expanded symmetric domain of type IV 
with respect to the Weyl group.
\end{abstract}

\maketitle
%
%

\tableofcontents

 \section{Introduction}

A quadruple of positive integers $W:=(a,b,c;h)$ is called 
a {\it regular system of weights} if the rational function 
$\chi_W\!=T^{-h}\frac{(T^h-T^a)(T^h-T^a)(T^h-T^c)}{(T^a-1)(T^b-1)(T^c-1)}$ 
develops in a Laurent polynomial 
and satisfies a suitable reducedness condition 
(\cite{sa:weight}, see subsection \ref{ssec:w}).
Then, $\chi_W$ is a sum of Laurent monomials and 
the exponents of the monomials are called 
the {\it exponents} of $W$. 
The {\it smallest exponent}, given by $a+b+c-h$, 
is denoted by $\varepsilon_W$. 
The regularity condition on $W$ is equivalent to that 
a degree $h$ weighted homogeneous polynomial $f_W\in A:=\C[x,y,z]$ 
in three variables $x$, $y$ and $z$ 
of weights $a$, $b$ and $c$, respectively 
with a generic choice of coefficients 
defines a hypersurface in ${\mathbb A}^3$ 
having an isolated singular point at the origin.

Motivated by the theory of primitive forms associated to the polynomial $f_W$, 
we asked to construct a generalization of a root system 
and a Lie algebra for any regular system of weights $W$ 
(\cite{sa:around,sa:primitive}).
In fact, by taking the set of vanishing cycles in the Milnor fiber of $f_W$, 
a finite root system of type ADE 
or an elliptic root system of type 
$E_6^{(1,1)}$, $E_7^{(1,1)}$ or $E_8^{(1,1)}$   
\cite{sa:elliptic} is 
associated to a regular system of weights $W$ 
with $\varepsilon_W\!=\!1$ or $0$, respectively, 
where the set of exponents of the weight system coincides 
with the set of Coxeter exponents of the root system. 
However, the vanishing cycles are transcendental 
object and are hard to study further for the cases of $\varepsilon_W<0$. 
Then, based on the duality theory of the weight systems 
\cite{sa:duality,t:1} and the (homological) mirror symmetry \cite{ko}, 
the third author \cite{t:2} proposed 
to use the triangulated category of graded matrix 
factorizations of $f_W$ 
introduced by \cite{t:2} and Orlov \cite{o:2} independently, 
where the root system appears as the set of 
the isomorphism classes of the exceptional objects via 
the Grothendieck group of the category. 

In our previous paper \cite{KSTade}, we showed that, 
for any regular system of weights $W$ of type ADE (i.e., $\varepsilon_W=1$), 
the triangulated category $\tHMF{A}{f_W}$ of graded matrix factorizations 
of $f_W$ is equivalent to the bounded derived category 
of finitely generated modules over the path algebra 
of a Dynkin quiver of the corresponding ADE-type.  
Due to a theorem of Gabriel \cite{g:1}, 
this implies that one gets the root system of type ADE 
in the Grothendieck group $K_0(\tHMF{A}{f_W})$ as expected 
(see \cite{t:2} for $A_l$ case). 

The present paper studies the 
category $\tHMF{A}{f_W}$ associated to regular systems of weights $W$ 
with $\varepsilon_W=-1$  and $a_0=0$ (the second condition 
means, by definition, there are no exponents equal to 0, 
and we call such $W$ 
{\it nondegenerate}). 
The set of indecomposable objects of $\tHMF{A}{f_W}$ 
is no more simple to describe 
as opposed to the case of type ADE. 
However, we can still find
a strongly exceptional collection which generates 
the category and gives a good basis of the generalized 
root system in $K_0(\tHMF{A}{f_W})$. 
More precisely, the main theorem (Theorem \ref{thm:main}) states that, 
for any regular system of weights $W$ with $\varepsilon_W=-1$ 
and $a_0=0$, 
there exists a strongly exceptional collection 
in the category $\tHMF{A}{f_W}$, 
the associated quiver (in a generalized sense, see 
subsection \ref{ssec:GDD}) to which is given by the following diagram 

\includegraphics{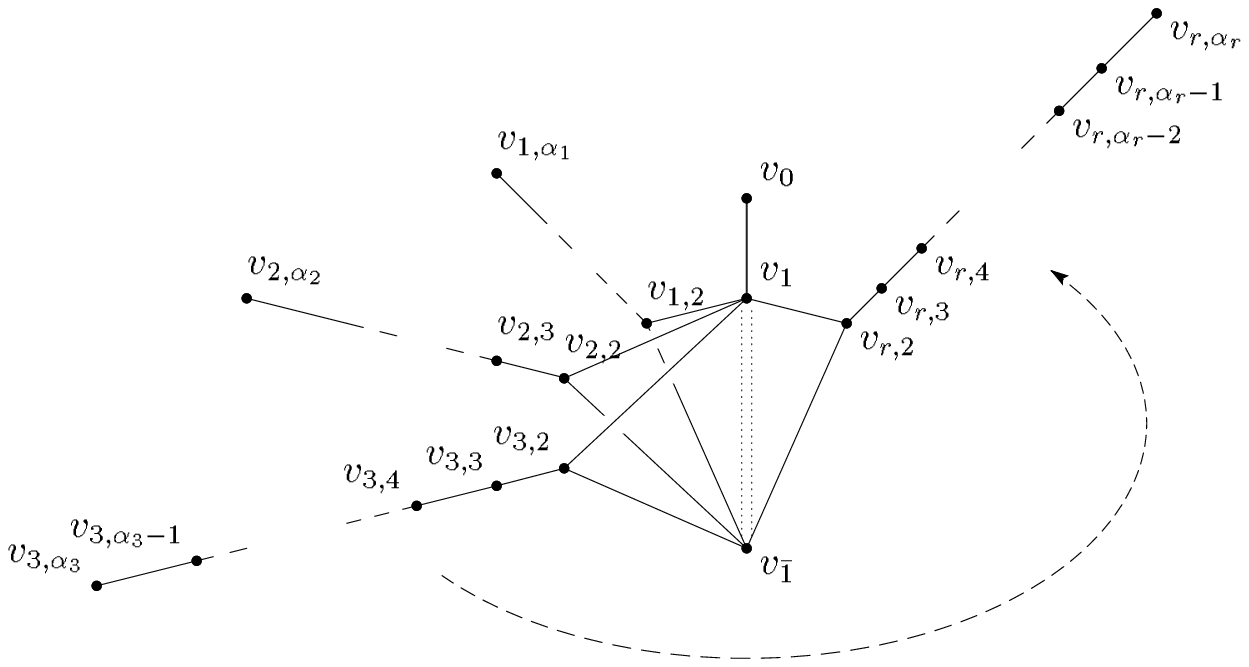}

\noindent
with appropriate orientations of arrows for the edges 
(see Figures \ref{fig:quiver1}, \ref{fig:quiver1'} and \ref{fig:quiver2} 
in subsection \ref{ssec:DeltaW}), 
where the multi-set $A_W=\{\alpha_1,\dots,\alpha_r\}$ of positive integers 
is the {\it signature of} $W$ (see Definition \ref{defn:Aw}).

\medskip
Let us discuss some background and consequences of Theorem \ref{thm:main}. 

1. There are 14+8 regular systems of weights 
with $\varepsilon_W\!=\!-1$ and $a_0\!=\!0$. The first 14 cases 
define exceptional unimodular 
singularities of Arnold \cite{arnold}, who
found an involutive bijection, called the strange duality, among their 
numerical invariants. 
The strange duality was reconstructed by the $*$-duality 
among regular systems of weights \cite{sa:duality} in terms of 
the characteristic polynomial, 
which is understood as a mirror symmetry 
in \cite{kawai-yang,t:1}. 
Then, Theorem \ref{thm:main} implies that 
{\it the lattice of the vanishing cycles of $f_W$ is obtained 
by the Grothendieck group of 
$\tHMF{A}{f_{W^*}}$ for the dual weight system $W^*$} 
as explained in subsection \ref{ssec:strange}. 
In particular, 
{\it the set of exponents of $W$ coincides with 
the set of Coxeter exponents of the root system in $K_0(\tHMF{A}{f_{W^*}})$} 
(see Remark \ref{rem:key}).

2. For those fourteen regular systems of weights $W$ 
corresponding to exceptional unimodular singularities, 
the root lattice $(K_0(\tHMF{A}{f_{W^*}}),\chi+{}^t\chi)$ 
is an indefinite lattice of sign $(2,\mu_{W^*}-2)$, where 
$\mu_{W^*}$ is the Milnor number of the Milnor fiber of $f_{W^*}$ and 
$\chi$ is the Euler pairing. 
The Coxeter transformation defined as 
a product of reflections associated to 
the underlying graph of the above quiver is identified with 
the Auslander-Reiten translation $\tau_{AR}$ 
and, hence, is of {\em finite} order $h$ (see Remark \ref{rem:key}). 
The Weyl group ${\mathcal W}$ generated by those reflections 
acts on the expanded symmetric domain 
${\mathcal B}_\C\!
:=\!\{\varphi\!\in\! \Hom_\R(K_0(\tHMF{A}{f_{W^*}})\otimes_\Z \R,\C) 
\mid \Ker(\varphi)\!>\!0\}$ of type IV 
\cite{sa:primitive, sa:period}. It is shown 
that the eigenvectors of the Coxeter transformation whose eigenvalues 
belong to the primitive $h$-th roots of unity are regular 
with respect to ${\mathcal W}$ 
(see \cite{sa:around}). 
This fact will conjecturally play an important role to construct 
a flat structure on the quotient variety 
${\mathcal B}_\C/{\mathcal W}$ which should be 
the space of stability conditions 
on $\tHMF{A}{f_{W^*}}$ in the sense of Bridgeland \cite{bd:1}.

3. Due to a theorem of Bondal and Kapranov \cite{bo:1,bk:1}, we see
that {\it the triangulated category $\tHMF{A}{f_W}$ is equivalent to 
the bounded derived category of finitely generated 
modules over the path algebra with relations 
corresponding to the quiver above}
(Corollary \ref{cor:main}). 
Recently, a parallel statement is proven independently by 
Lenzing and de la Pena  \cite{L-delaP} 
in the framework of the weighted projective lines 
by Geigle and Lenzing \cite{GeLe1} by combining it 
with Orlov's arguments in \cite{o:2}.

\medskip
The construction of the present paper is as follows. 

Section \ref{sec:cat} is devoted to the preparation of the 
categories of our study in three equivalent formulations.
In subsections \ref{ssec:Dbsing} and \ref{ssec:MCM}, 
following Orlov \cite{o:2} 
(see also Buchweitz \cite{buch:1}), 
we recall the triangulated category 
$\Dbsing(R_W)$ of singularity and 
the triangulated category $\CMgb(R_W)$ 
of graded maximal Cohen-Macaulay modules over $R_W=A/(f_W)$, 
respectively.
Then, in subsection \ref{ssec:MF}, we recall 
the triangulated category $\tHMF{A}{f_W}$ of 
graded matrix factorizations from \cite{KSTade}.

In section \ref{sec:Serre}, we show the existence of 
the Serre duality and 
the Auslander-Reiten triangles 
in the triangulated category $\CMgb(R)$ 
(Proposition \ref{prop:SerreDual}). 
This fact may be well-known among experts. We are grateful to
Prof. Iyama who explained the result to us. 

In section \ref{sec:mainI}, we show the other basic result which we 
use in the proof of the main theorem:
{\it a right admissible full triangulated subcategory $\T'$ of 
$D^{gr}_{Sg}(R)$ satisfying the conditions$:$
\begin{enumerate}
\item the degree shift functor $\tau$ is an autoequivalence of $\T'$,
\item $\T'$ has an object $E$ which is isomorphic to $R/{\mathfrak m}$ 
in $D^{gr}_{Sg}(R)$, 
\end{enumerate}
is equivalent to the category $D^{gr}_{Sg}(R)$ itself 
as a triangulated category 
}
(Theorem \ref{thm:mainI}).

Section \ref{sec:main} is devoted to stating our main results. 
In subsection \ref{ssec:w}, we recall 
regular systems of weights and related notion. 
In subsection \ref{ssec:GDD}, 
we prepare a generalized notion of quivers. 
In this formulation, 
we define the quivers associated to 
regular systems of weights $W$ of $\varepsilon_W=-1$ with genus $a_0=0$ 
in subsection \ref{ssec:DeltaW}. 
In subsection \ref{ssec:main}, 
we state the main Theorem (Theorem \ref{thm:main}).  
It is obtained as a consequence of the structure theorem
(Theorem \ref{thm:key}) on the category $\tHMF{A}{f_W}$, 
where we describe the Auslander-Reiten triangles for all objects 
forming the exceptional collection. 
The proof of Theorem \ref{thm:key}, stated in subsection \ref{ssec:proof-key}, 
is based on explicit data 
of graded matrix factorizations of $f_W$ for 
a regular system of weights $W$ with $\varepsilon_W=-1$ and $a_0=0$, 
which are given in section \ref{sec:list}.
We hope these data can give a well-defined stability condition \cite{bd:1}, 
which is one of future directions. 

Throughout the paper, we denote by $k$ 
an algebraic closed field of characteristic zero.

{\bf Acknowledgment} :\  We would like to 
thank O.~Iyama for valuable comments 
on the relation between the Auslander-Reiten transformation 
and the Serre duality. 
This work was partly supported by Grant-in Aid for Scientific Research 
grant numbers 16340016, 17654015, 17740036 and 19740038 
from the Ministry of Education, Culture, Sports, Science and Technology, 
Japan.

\section{Some equivalent categories}
\label{sec:cat}

Let $k$ be an algebraic closed field of characteristic zero. 
For a positive integer $h$, 
let $R:=\oplus_{s\in \frac{2}{h}\Z_{\ge 0}}R_s$ be a 
commutative Noetherian $(2\Z/h)$-graded 
ring of dimension $d (\ge 0)$ with $R_0=k$. 
This ring $R$ defines a graded isolated singularity, i.e., 
the graded localization $R_{({\mathfrak p})}$ is regular 
for any graded prime ${\mathfrak p}\ne {\mathfrak m}$, 
where ${\mathfrak m}:=\oplus_{s\in \frac{2}{h}\Z_{> 0}}R_s$. 
By a graded $R$-module, 
we always mean a graded $R$-module with degrees 
only in $2\Z/h$. 
Namely, a graded $R$-module $M$ decomposes 
into the direct sum $M=\oplus_{s\in\frac{2}{h}\Z}M_s$.
For two graded $R$-modules $M$ and $N$, 
a graded $R$-homomorphism $g$ of degree $t\in 2\Z/h$ 
is an $R$-homomorphism $g:M\to N$ such that 
$g(M_s)\subset N_{s+t}$ for any $s\in 2\Z/h$. 
\subsection{Category of graded singularities}
\label{ssec:Dbsing}

\begin{defn}\label{defn:tau}
Denote by $\gr{R}$ the abelian category of finitely generated graded 
$R$-modules, in which morphisms are $R$-homomorphism of degree zero. 
The degree shift of $M\in\gr{R}$, denoted by $\tau M$, is defined by 
$(\tau M)_s:=M_{s+\frac{2}{h}}$. 
This $\tau$ naturally induces an auto-equivalence functor on $\gr{R}$,
which we denote by the same symbol $\tau$.
\end{defn}
We have $\Ext^i_R(M,N)\simeq \oplus_{n\in\Z} 
\Ext^i_{\gr{R}}(\tau^{-n}M,N)\simeq\oplus_{n\in\Z} 
\Ext^i_{\gr{R}}(M,\tau^n N)$ since $R$ is Noetherian.
In particular, for $i=0$, 
\begin{equation*}
\Hom_R(M,N)\simeq \oplus_{n\in\Z} \Hom_{\gr{R}}(\tau^{-n}M,N)
\simeq\oplus_{n\in\Z}\Hom_{\gr{R}}(M,\tau^n N)
\end{equation*}
forms a graded $R$-module, where the grading of each homogeneous piece is 
defined as $2n/h$. 
Note also that any graded projective module is free since 
$R$ is finitely generated over $R_0=k$.
Denote by ${\rm grproj-}R$ the exact category 
(= the extension-closed full 
additive subcategory in $\gr{R}$) of graded projective modules.
\begin{defn}[Orlov \cite{o:2}]
The triangulated category $D^{gr}_{Sg}(R)$, called 
the category of the singularity $R$, is defined as the 
quotient $D^b(\gr{R})/D^b({\rm grproj-}R)$.
We denote by $T$ the translation functor
\footnote{It is often called the shift functor and denoted by $[1]$.
In this paper, in order to avoid confusions, 
we always mean by $\tau$ the degree shift functor and 
by $T$ the translation functor. 
Also, the Auslander-Reiten translation will be denoted 
by $\tau_{AR}$. } on the triangulated category 
$D^{gr}_{Sg}(R)$.
\end{defn}
\subsection{Category of graded maximal Cohen-Macaulay modules}\hfill
\label{ssec:MCM}

The above definition of $D^{gr}_{Sg}(R)$ is simple, however, 
it is not easy to understand morphisms between objects 
since they are defined in the localized category.
Therefore, we recall 
some categories equivalent to $D^{gr}_{Sg}(R)$.
In this subsection, we recall the triangulated category $\CMgb(R)$ 
of graded maximal Cohen-Macaulay modules; 
we refer to \cite{y} for terminologies and the statements presented here. 
\begin{defn}
An element $M\in\gr{R}$ 
is called a {\em graded maximal Cohen-Macaulay module} 
if $\Ext^i_R(R/{\mathfrak m},M)=0$ for $i<d$, where 
$d$ is the dimension of the ring $R$.
We denote 
the full subcategory of $\gr{R}$ consisting of 
all graded maximal Cohen-Macaulay modules over $R$ 
by $\CMg(R)$, which forms an exact category. 
\end{defn}
Recall that an element $K_R\in\CMg(R)$ is called 
a {\em canonical module} of $R$ if 
$\Ext^i_R(R/\fm,K_R)\simeq 0$ for $i\ne d$ 
and $\Ext^d_R(R/\fm,K_R)\simeq k$. 
\begin{lem}\label{lem:CM}
The following conditions are equivalent$:$
\begin{enumerate}
\item $M$ is a graded maximal Cohen-Macaulay module,
\item $H_{\fm}^i(M)=0$ for $i\ne d$, 
where $H_{\fm}^\bullet$ is the local cohomology functor 
with support on $\{\fm\}$ defined by 
$H_{\fm}^i(M):=\displaystyle{\lim_{\longrightarrow}}\  
\Ext^i_R(R/R_{\ge n},M)$, $R_{\ge n}:=\oplus_{i\in \frac{2}{h}\Z_{\ge n}}R_i$,
\item $\Ext^i_R(M,K_R)=0$ for $i>0$.
\end{enumerate}
\qed\end{lem}
\begin{defn}\label{defn:G}
The ring $R$ is called {\em Gorenstein} 
if the injective dimension of $R$ is finite and the canonical module $K_R$
is isomorphic to $\tau^{-\varepsilon(R)}R$ for some $\varepsilon(R)\in\Z$. 
The integer $\varepsilon(R)$ is called the {\em Gorenstein parameter} of $R$. 
\end{defn}
\begin{lem}
For a Gorenstein ring $R$, 
$\CMg(R)$ is a Frobenius category, i.e., it has enough projectives 
and enough injectives and the projectives coincide with the injectives.
\qed\end{lem}
\begin{defn}
For a Gorenstein ring $R$, 
we define an additive category $\CMgb(R)$ as follows: 
objects of it are graded maximal Cohen-Macaulay modules over $R$ 
and, for any $M,N\in\CMgb(R)$, 
the space of morphisms $\Homub_{\gr{R}}(M,N)$ 
is given by $\Hom_{\gr{R}}(M,N)/{\mathcal P}(M,N)$, where 
${\mathcal P}(M,N)$ is the subspace consisting of elements 
factoring through projectives, i.e., 
$g\in {\mathcal P}(M,N)$ if and only if $g=g''\circ g'$ for $g':M\to P$ and 
$g'':P\to N$ with a projective object $P$. 
\end{defn}
The stable category of a Frobenius category forms a 
triangulated category (Happel \cite{h})). 
Since $\CMgb(R)$ is the stable category of the 
Frobenius category $\CMg(R)$, one obtains that: 
\begin{prop}
The stable category $\CMgb(R)$ forms 
a triangulated category. 
\qed\end{prop}
The following important fact 
is implicit in Orlov \cite{o:2}:
\begin{thm}[Section 1.3. in \cite{o:2} (see also 
Buchweitz \cite{buch:1})]
For a Gorenstein ring $R$, 
there is an equivalence 
$\CMgb(R)\simeq D^{gr}_{Sg}(R)$
as triangulated categories.
\qed\end{thm}
\subsection{Category of graded matrix factorizations}\hfill
\label{ssec:MF}

Consider the case when $R$ is a quotient algebra $A/(f)$ 
of a graded Noetherian 
regular algebra $A=\oplus_{i\in \frac{2}{h}\Z_{\ge 0}}A_i$ with $A_0=k$ and 
an element $f\in A_2$ which is a non-zero divisor. 
Since $A/(f)$ defines a hypersurface, $R$ is Gorenstein. 
Recall that $\tau$ is the degree shifting operator 
defined in Definition \ref{defn:tau}. 
\begin{defn}
For a non-zero element $f\in A_2$, 
we define an additive category $\tMF{A}{f}$ as follows. 
Objects of it are graded matrix factorizations $M$ of $f$ defined by 
\begin{equation*}
\olF:=\Big(\mf{F_0}{f_0}{f_1}{F_1}\Big)\ ,
\end{equation*}
where $F_0$ and $F_1$ are graded free $A$-modules of finite rank, 
$f_0:F_0\to F_1$ is a graded ${A}$-homomorphism of degree zero, 
$f_1:F_1\to F_0$ is a graded ${A}$-homomorphism of degree two 
such that $f_1f_0=f\cdot\Id_{F_0}$ and $f_0f_1=f\cdot\Id_{F_1}$.
A morphism $g:\olF\to\olFF$ in the category $\tMF{A}{f}$ 
is a pair $g=(g_0,g_1)$ 
of graded $A$-homomorphisms 
$g_0:F_0\to F'_0$ and $g_1:F_1\to F'_1$ 
of degree zero satisfying 
$g_1f_0=f'_0g_0$ and $g_0f_1=f'_1g_1$.
\end{defn}
For a graded matrix factorization $\overline{F}$, 
by definition, 
the rank of $F_0$ coincides with that of $F_1$, 
which we call the {\em rank} 
of the matrix factorization $\overline{F}$.

Componentwise monomorphisms 
and epimorphisms equip the additive category $\tMF{A}{f}$ 
with an exact structure. 
Moreover, we have the following:
\begin{lem}[\cite{o:2}]
$\tMF{A}{f}$ is a Frobenius category.
\qed\end{lem}
A morphism $g=(g_0,g_1):\overline{F}\to \overline{F'}$ 
is called null-homotopic if there are graded $A$-homomorphisms 
$\psi_0:F_0\to F'_1$ of degree minus two and 
$\psi_1:F_1\to F'_0$ of degree zero such that 
${g}_0=f'_1\psi_0+\psi_1f_0$ and 
${g}_1=\psi_0f_1+f'_0\psi_1$.
Morphisms factoring through projectives in $\tMF{A}{f}$ are 
null-homotopic morphisms.
\begin{defn}
We denote by $\tHMF{A}{f}$ the stable (homotopy) 
category of the Frobenius category $\tMF{A}{f}$. 
\end{defn}
\begin{prop}
[See Eisenbud \cite{e:1}, Orlov \cite{o:2} and Yoshino \cite{y}, for example.]
$\tHMF{A}{f}$ is a triangulated category 
which is equivalent to $\CMgb(A/(f))$. 
The equivalence is given by the correspondence
$$
\overline{F}=\Big(\mf{F_0}{f_0}{f_1}{F_1}\Big) \mapsto 
M:={\rm Coker}(f_1).
$$
\qed\end{prop}
\begin{rem}
An object $M\in \tHMF{A}{f}$ is zero if and only if it is a direct sum of 
the graded matrix factorizations of the forms 
$(\mfs{\tau^n(A)}{1}{f}{\tau^n(A)})\in\tMF{A}{f}$ and 
$(\mfs{\tau^{n'}(A)}{f\ }{1\ }{\tau^{n'+h}(A)})\in\tMF{A}{f}$ 
for some $n,n'\in\Z$.
\end{rem}
\subsection{Some basic properties 
of the category of graded matrix factorizations}\hfill
\label{ssec:NF-detail}

For later necessity, 
we discuss the structure of $\tHMF{A}{f}$ in more detail.

The auto-equivalence functor $\tau$ on $\gr{A}$ induces 
an auto-equivalence on $\tHMF{A}{f}$, 
which we denote by the same notation $\tau$.
Explicitly, the action of $\tau$ takes an object $\overline{F}$ to the 
object
\begin{equation*}
\tau \overline{F}:=
\Big(\mf{\tau F_0}{\tau(f_0)}{\tau(f_1)}{\tau F_1}\Big),
\end{equation*}
and takes a morphism ${g}=({g}_0,{g}_1)$ to the morphism 
$\tau({g}):=(\tau({g}_0),\tau({g}_1))$.
The translation functor $T$ on $\tHMF{A}{f}$ takes an object 
$\overline{F}$ to the object
\begin{equation*}
T\overline{F}
:=\Big(\mf{F_1}{-f_1\quad}{-\tau^h(f_0)\quad}{\tau^hF_0}\Big),
\end{equation*}
and takes a morphism ${g}=({g}_0,{g}_1)$ to the morphism 
$T({g}):=({g}_1,\tau^h({g}_0))$.
The following fact is straightforward by definition, 
but plays an important role 
in the study of $\tHMF{A}{f}$:
\begin{prop}
$T^2=\tau^h$ on $\tHMF{A}{f}$.
\qed\end{prop}
Next, we explain 
the triangulated structure in $\tHMF{A}{f}$. 
First, we recall the mapping cone. 
\begin{defn}
For a morphism ${g}=({g}_0,{g}_1)\in 
\Hom_\tMF{A}{f}(\overline{F},\overline{F'})$,  
we define a mapping cone $C({g})\in\tMF{A}{f}$ as 
\begin{equation*}
C({g}):=\Big(\mf{F_1\oplus F_0'}{c_0}{c_1}{\tau^h F_0\oplus F'_1}\Big), 
\qquad 
 c_0:=\bp -f_1 & 0 \\ {g}_1 & f'_0\ep ,\quad 
 c_1:=\bp -\tau^h(f_0) & 0 \\ \tau^h({g}_0) & f'_1\ep .
\end{equation*}
\label{defn:cone_Z}
\end{defn}
We sometimes denote this cone by $C(\overline{F}\to\overline{F'})$ 
when omit writing the morphism explicitly. 

Note that there exist morphisms 
$\overline{F'}\stackrel{{}^t({\rm id},0)}{\to} C({g})$ and 
$C({g})\stackrel{(0,-{\rm id})}{\to} T\overline{F}$.
By definition of the triangulated structure on $\tHMF{A}{f}$, 
one can easily see that 
\begin{prop}
Each exact triangle in $\tHMF{A}{f}$ 
is isomorphic to a triangle of the form 
\begin{equation*}
\overline{F}\stackrel{{g}}{\longrightarrow}
\overline{F'}\stackrel{{}^t({\rm id},0)}{\longrightarrow}
C({g})\stackrel{(0,-{\rm id})}{\longrightarrow}
T\overline{F}
\end{equation*}
for some $\overline{F},\overline{F'}\in\tMF{A}{f}$ 
and ${g}\in\Hom_\tMF{A}{f}(\overline{F},\overline{F'})$. 
\qed\end{prop}
Let $\overline{F}=(\mfs{F_0}{f_0}{f_1}{F_1})\in\tHMF{A}{f}$
be a graded matrix factorization of rank $r$. 
Choose homogeneous free basis 
$(b_1,\dots,b_r;\bb_1,\dots,\bb_r)$ such that 
$F_0=b_1{A}\oplus\cdots\oplus b_r{A}$ and 
$F_1=\bb_1{A}\oplus\cdots\oplus\bb_r{A}$. 
Then, the graded matrix factorization $\overline{F}$ is expressed 
as a pair $(Q,S)$ of $2r$ by $2r$ matrices, 
where $S$ is the diagonal matrix of the form 
$S:=\diag(s_1,\dots,s_r; \sb_1,\dots,\sb_r)$ such that 
$s_i=\deg(b_i)$ and $\sb_i=\deg(\bb_i)-1$, 
$i=1,2,\dots,r$, and $Q$ is given by 
\begin{equation}\label{Q}
Q= \bp \0 & {q_0} \\ {q_1} & \0 \ep,\qquad {q_0},{q_1}\in\Mat_r({A}), 
\end{equation}
with $q_0$ and $q_1$ the matrix expressions of the 
graded $A$-homomorphisms $f_0: F_0\to F_1$ and $f_1: f_1\to f_0$, 
respectively. 
Namely, they are defined as 
$f_0 (b_1,\dots,b_r)=(\bb_1,\dots,\bb_r)q_0$ and 
$f_1 (\bb_1,\dots,\bb_r)=(b_1,\dots, b_r)q_1$. 
By definition, $(Q,S)$ satisfies  
\begin{equation}\label{S}
Q^2=f\cdot\1_{2r},\qquad -SQ+QS+2EQ= Q, 
\end{equation}
where $E\in\Der_k(A)$ is the derivation corresponding to 
the infinitesimal generator of $k^\times$-action 
(see eq.(\ref{Euler-vf})). 
We call this $S$ a {\em grading matrix} of $Q$. 
This procedure $\overline{F}\mapsto (Q,S)$ 
gives a triangulated equivalence between 
the triangulated category $\tHMF{A}{f}$ and 
the triangulated category $D^b_\Z(\A_f)$ introduced in \cite{t:2}. 
In particular, the latter category $D^b_\Z(\A_f)$ is defined 
as the cohomology of a DG-category of twisted complexes. 
This implies that $D^b_\Z(\A_f)$ and then $\tHMF{A}{f}$ are 
enhanced triangulated categories 
in the sense of Bondal-Kapranov \cite{bk:1}. 
In this paper, we often represent a graded matrix factorization 
$\overline{F}=(\mfs{F_0}{f_0}{f_1}{F_1})$ by its matrix 
representation $(Q,S)$. 
\begin{defn}\label{defn:t}
Let $t:\Ob(\tHMF{A}{f})\to \Ob(\tHMF{A}{f})$ be 
the bijection induced by the correspondence 
$t:(Q,S)\mapsto (^tQ, -S)$, 
where $^tQ$ is the transpose of the matrix $Q$. 
This $t$ is lifted to be a contravariant equivalence functor 
on $\tHMF{A}{f}$, which we denote by the same notation $t$. 
\end{defn}
\begin{prop}\label{prop:t}
On $\tHMF{A}{f}$, one has the following identities$:$ 
\begin{equation*}
 \tau t \tau=t,\qquad T t T =t. 
\end{equation*}
\qed\end{prop}

\subsection{Further remark}\hfill
\label{ssec:future-remark}

Let $\T$ be one of the equivalent triangulated categories 
$D^{gr}_{Sg}(A/(f))$, $\CMgb(A/(f))$ and $\tHMF{A}{f}$.
Since we assume that the ring $R=A/(f)$ defines an isolated singularity, 
$\T$ is Krull-Schmidt (see \cite{KSTade}), that is, \\
(a)\ for any two objects $M,M'\in\T$, 
$\Hom_\T(M,M')$ is of finite rank over $k$; \\
(b)\ for any object $M\in\T$ and 
any idempotent $e\in{\rm Hom}_\T(M,M)$, 
there exists an object $M'\in \T$ 
and a pair of morphisms 
${g}\in {\rm Hom}_\T(M, M')$, ${g}'\in {\rm Hom}_\T(M',M)$ 
such that ${g}'{g}=e$ and ${g}{g}'=\Id_{M'}$. 

\section{Serre duality}
\label{sec:Serre}

In this section, we assume $R$ is a Gorenstein ring and 
show the existence of the Serre duality and 
the Auslander-Reiten triangles 
in the triangulated category $\CMgb(R)$. 
For terminologies, we again refer to \cite{y}. 
\begin{defn}
Consider a finite presentation of $M\in\gr{R}$ by graded free modules,
$F_1\stackrel{f}{\to} F_0\to M\to 0$.
Define ${\rm tr}(M)$ by the following exact sequence
$$
0\to\Hom_{{R}}(M,R)\to \Hom_{{R}}(F_0,R)
\stackrel{\Hom_{{R}}(f,R)}{\longrightarrow}
\Hom_{{R}}(F_1,R)\to {\rm tr}(M)\to 0,
$$
i.e., ${\rm tr}(M)={\rm Coker}(\Hom_{{R}}(f,R))$. 
The graded module 
${\rm tr}(M)$ is called the {\em Auslander transpose} of $M$.
\end{defn}
The Auslander transpose ${\rm tr}(M)$ is unique up to free summands. 
Since we shall only deal with properties that are 
independent of free summands of ${\rm tr}(M)$, the above definition 
will be sufficient.
\begin{defn}
For a graded $R$-module $M\in \gr R$, 
consider a long exact sequence 
\begin{equation*}
0\to N\to F_{n-1}\to F_{n-2}\to \dots \to F_1\to F_0\to M\to 0
\end{equation*}
in $\gr{R}$, where each $F_i$ is graded free.
The {\em reduced $n$-th syzygy} ${\rm syz}^n(M)$ of $M$ 
is the graded $R$-module obtained from $N$ by deleting 
all graded free summands.
\end{defn}
The reduced $n$-th syzygy ${\rm syz}^n(M)$ 
is uniquely determined by $M$ and $n$ up to isomorphism.
\begin{defn}
For a graded $R$-module $M\in \gr R$, 
the {\em Auslander-Reiten (AR-) translation} 
$\tau_{AR}(M)\in\gr R$ is defined by 
\begin{equation*}
\tau_{AR}(M):={\rm Hom}_{R}({\rm syz}^d({\rm tr}(M)),K_{R}).
\end{equation*}
\end{defn}
\begin{rem}
For $M\in \CMg(R)$ which is reduced, i.e., 
has no free direct summands, 
we have ${\rm syz}^2({\rm tr}(M))\simeq {\rm Hom}_{R}(M,R)$.
\end{rem}
\begin{lem}
For $M\in\CMgb(R)$, we have
\begin{equation*}
\tau_{AR}(M)\simeq T^{d-2}\tau^{-\varepsilon(R)}M, 
\end{equation*}
where $\varepsilon(R)$ is the Gorenstein parameter of $R$ 
defined in Definition \ref{defn:G}. 
\end{lem}
\begin{pf}
This easily follows from that $R$ is Gorenstein and the definition of 
the translation functor $T$ on $\CMgb(R)$.
\qed\end{pf}
\begin{cor}\label{cor:tau-AR}
Suppose that $R$ defines a weighted homogeneous hypersurface 
as in subsection \ref{ssec:MF}. 
Let $[\tau_{AR}]$ denotes the induced map of 
$\tau_{AR}:\CMgb(R)\to\CMgb(R)$ on the Grothendieck group. 
Then, $[(\tau_{AR})^h]=(-\Id)^{hd}$ holds 
and hence $[\tau_{AR}]$ is of finite order. 
\qed\end{cor}
\begin{prop}[Auslander-Reiten duality \cite{ar:lmn1273}]\label{prop:SerreDual}
Let $R$ be a graded Cohen-Macaulay ring of dimension $d$ which defines 
an isolated singularity and has the canonical module $K_R$.
Then, there exists the following bi-functorial isomorphism of degree zero
\begin{equation}
 \Ext^d_{{R}}(\Homub_{{R}}(M,N),K_R)\simeq 
\Ext^1_{{R}}(N,\tau_{AR}(M)).
\end{equation}
\qed\end{prop}
By this Proposition, we see that the triangulated category 
$\CMgb(R)\simeq D^{gr}_{Sg}(R)$ has a 
Serre functor$:$
\begin{thm}\label{thm:SerreDual}
\footnote{
We thank O.~Iyama for explaining to us that Auslander-Reiten duality 
implies the Serre duality.
}
The functor $\S:=T\tau_{AR}=T^{d-1}\tau^{-\varepsilon(R)}$ 
is the Serre functor on $\CMgb(R)$. 
More precisely, $\S$ is an auto-equivalence functor which induces 
bi-functorial isomorphisms
$$
\Hom_k(\Homub_{\gr{R}}(M,N),k)\simeq
\Homub_{\gr{R}}(N,\S M),
\quad M,N\in \CMgb(R).
$$
\end{thm}
\begin{pf}
Note that $\Homub_{{R}}(M,N)$ is a graded $R$-module of 
finite length since $R$ is an isolated singularity.
Hence, we have the following isomorphism of degree zero by the 
local duality theorem
\begin{align*}
\Hom_k(\Homub_{\gr{R}}(M,N),k)\simeq &
\Hom_{\gr{R}}(\Homub_{{R}}(M,N),E_{R}(R/{\mathfrak m}))\\
\simeq & \Hom_{\gr{R}}
(H_{{\mathfrak m}}^0(\Homub_{{R}}(M,N)),
E_{R}(R/{\mathfrak m}))\\
\simeq & \Ext^d_{\gr{R}}
(\Homub_{{R}}(M,N),K_{R}),
\end{align*}
where $E_{R}(R/{\mathfrak m})$ is the injective envelope of 
the graded $R$-module $R/{\mathfrak m}$.
Since $R$ is Gorenstein, by Lemma \ref{lem:CM} (iii), 
one has $\Ext^1_R(N,F)=0$ for 
$N\in\CMg(R)$ and a free module $F$,  
and hence one sees that $\Ext^1_{\gr{R}}(N,\tau_{AR}(M))\simeq 
\Homub_{\gr{R}}(N, T\tau_{AR}(M))$. 
(Recall that there exists an exact sequence $0\to \tau_{AR}(M)\to F
\to T\tau_{AR}(M)\to 0$ in $\CMg(R)$.)
Therefore, we have the canonical isomorphism 
$\Hom_k(\Homub_{\gr{R}}(M,N),k)\simeq
\Homub_{\gr{R}}(N, T\tau_{AR}(M))$ of degree zero.
\qed\end{pf}
\begin{rem}
This theorem holds true even if 
we replace the $\Z$-grading by $L(p)$-grading in the sense of 
Geigle-Lenzing \cite{GeLe1,GeLe2}, 
since the generalization of the Auslander-Reiten duality (Proposition \ref{prop:SerreDual}) 
is straightforward. 
\end{rem}
Recall the notion of Auslander-Reiten (AR-)triangles 
(see \cite{h},\cite{y}; for an Auslander-Reiten sequence or 
equivalently an almost split sequence, see\cite{ar:ass}.).
A morphism $g$ is called {\em irreducible} 
if $g$ is neither a split monomorphism nor 
a split epimorphism but for any factorization $h=g_1g_2$ 
either $g_1$ is a split epimorphism or $g_2$ is a split monomorphism.
\begin{defn} 
An exact triangle in a Krull-Schmidt triangulated category $\T$ 
\begin{equation}\label{ARtri-general}
X\overset{u}{\to} Y\overset{v}{\to} Z\overset{w}{\to} T(X)
\end{equation}
is called an {\em Auslander-Reiten (AR-)triangle} if 
the following conditions are satisfied: 
\begin{itemize}
 \item[(AR1)] $X, Z$ are indecomposable objects in $\T$. 
 \item[(AR2)] $w\ne 0$
 \item[(AR3)] If ${g}:W\to Z$ is not a split epimorphism, then 
there exists ${g}':W\to Y$ such that $v{g}'={g}$. 
\end{itemize}
We call such a triangle (\ref{ARtri-general}) an AR-triangle of $Z$. 
\end{defn}
\begin{prop}[{Happel \cite[Proposition 4.3]{h}}]
Suppose given an AR-triangle (\ref{ARtri-general}) 
in a Krull-Schmidt triangulated category $\T$. 
\begin{enumerate}
\item Any AR-triangle of $Z$ is isomorphic 
to the AR-triangle (\ref{ARtri-general}) 
as exact triangles. 
\item The morphisms $u$ and $v$ 
in the AR-triangle (\ref{ARtri-general}) are irreducible morphisms.
\end{enumerate}
\qed\end{prop}
We say that a Krull-Schmidt triangulated category {\em $\T$ has AR-triangles} 
if there exists an Auslander-Reiten (AR-)triangle (\ref{ARtri-general}) of $Z$ 
for any indecomposable object $Z\in\T$. 

Now, for $\T=\CMgb(R)$, 
Theorem \ref{thm:SerreDual} implies the followings. 
\begin{cor}\label{cor:SerreDual-AR}
The triangulated category $\CMgb(R)$ has AR-triangles. 
\end{cor}
\begin{pf}
For any indecomposable object $Z\in\T$, 
the AR-triangle is given by 
\begin{equation*}
\tau_{AR}(Z)\to AR(Z)\to Z\to T\tau_{AR}(Z)
\end{equation*}
for some $AR(Z)\in \CMgb(R)$, 
where the morphism $Z\to T\tau_{AR}(Z)$ is given by the Serre dual, 
in the sense in Theorem \ref{thm:SerreDual}, of the identity morphism 
on $Z$. This fact can be shown just by following the same argument 
as that in \cite[Chapter I. 4.6]{h}. 
\qed\end{pf}
Thus, by definition of AR-triangles, one obtains 
the following which will be employed in the proof of Theorem \ref{thm:key}. 
For $X,Y\in\T$, we denote $\hom_\T(X,Y):=\dim_k(\Hom_\T(X,Y))$. 
\begin{cor}\label{cor:AR}
For any indecomposable object $Z\in\T:=\CMgb(R)$, consider the AR-triangle 
\begin{equation*}
\tau_{AR}(Z)\to AR(Z)\to Z\to T\tau_{AR}(Z). 
\end{equation*}
Then, for any indecomposable object $W\in\T$, one has  
\begin{equation*}
 \begin{split}
 & \hom_\T(W,AR(Z)) 
 = (\hom_\T(W,Z)-\sigma) + (\hom_\T(W,\tau_{AR}(Z))-\sigma'), \\
 & \hom_\T(AR(Z),W) 
 = (\hom_\T(Z,W)-\sigma) + (\hom_\T(\tau_{AR}(Z),W)-\sigma'), 
 \end{split}
\end{equation*}
where $\sigma:=1$ if $W\simeq Z$ and zero otherwise, 
and $\sigma':=1$ if $W\simeq T^{-1}(Z)$ and zero otherwise. 
\qed\end{cor}


\section{Category generating theorem}
\label{sec:mainI}

In this section, 
we discuss about the generation of the category $\Dbsing(R)$. 
We show Theorem \ref{thm:mainI} 
and then Corollary \ref{cor:mainI} 
which is necessary to prove the structure theorem (Theorem \ref{thm:key}) 
of $\tHMF{A}{f}$. 

We first recall some definitions and facts concerning admissible categories 
and exceptional collections from \cite{bo:1,o:2}. 
\begin{defn}
Let $\T$ be a triangulated category and $\T'\subset \T$ 
a full triangulated subcategory.
The {\em right orthogonal} to $\T'$ 
is a full subcategory $(\T')^\perp\subset \T$
consisting of all objects $M$ such that ${\rm Hom}_\T(N,M)=0$ 
for any $N\in \T'$.
\end{defn}
\begin{defn}
Let $\T$ be a triangulated category and $\T'\subset \T$ 
a full triangulated subcategory.
We say that $\T'$ is {\em right admissible} 
if, for any $X\in \T$, there 
is an exact triangle $N\to X\to M\to TN$ 
with $N\in \T'$ and $M\in (\T')^\perp$.
\end{defn}
\begin{defn}\label{defn:EC}
An object $E$ of a $k$-linear triangulated 
category $\T$ is called {\em exceptional} 
if $\Hom_\T(E,T^nE)=0$ when $n\ne 0$ and $\Hom_\T(E,E)\simeq k$.
An {\em exceptional collection} is a sequence of exceptional objects 
$\E:=(E_1, \dots, E_l)$ satisfying the condition $\Hom_\T(E_i,T^nE_j)=0$
for all $n$ and $i>j$. 
Furthermore, an exceptional collection $\E=(E_1, \dots, E_l)$ 
is called a {\em strongly exceptional collection} 
if $\Hom_\T(E_i,T^nE_j)=0$ for all $i,j$ and all $n$ except 
for $n=0$. 
\end{defn}
We say that a triangulated category $\T$ is of {\em finite type} if, 
for any $E, E'\in \T$, $\Hom_\T(E,T^n E')$ is of finite rank over $k$ 
and in particular zero for almost all $n\in\Z$.
\begin{prop}[Bondal {\cite[Theorem 3.2]{bo:1}} (see also \cite{bp:1})]
Let $\T$ be a triangulated category of finite type and 
$\T':=\left<E_1,\dots, E_l\right>\subset \T$ a full 
triangulated subcategory generated by an exceptional collection 
$(E_1,\dots, E_l)$. 
Then, $\T'$ is right admissible.
\qed\end{prop}
The following result is the key lemma of this paper.
\begin{thm}\label{thm:mainI}
Let $\T'$ be a right admissible full triangulated subcategory of 
$\Dbsing(R)$, 
with $R$ a Gorenstein ring, 
satisfying the following conditions$:$
\begin{enumerate}
\item The shift functor $\tau$ on $\Dbsing(R)$ induces an autoequivalence of $\T'$,
\item $\T'$ has an object $E$ which is isomorphic to $R/{\mathfrak m}$ 
in $\Dbsing(R)$.
\end{enumerate}
Then $\T'$ is equivalent to $\Dbsing(R)$ as a triangulated category.
\end{thm}
\begin{pf}
By the equivalence $\CMgb(R)\simeq \Dbsing(R)$, 
we shall often represent objects in $\Dbsing(R)$ by the corresponding 
graded maximal Cohen-Macaulay modules.
First, recall the following characterization of free modules:
\begin{lem}[(see \cite{y})]\label{lem:free}
An object $M\in \CMg(R)$ is graded free if and only if
$\Ext^i_R(R/\fm,M)=0$ for $i\ne d$.
\qed\end{lem}
Take a minimal graded free resolution of $R/\fm\in\Dbsing(R)$ 
\begin{equation*}
\dots \to F_2\to F_1 \to F_0 \to R/{\mathfrak m}\to 0.
\end{equation*}
By definition of syzygy, one has 
\begin{equation}\label{syz-free}
 0\to \syz^{i+1}(R/\fm)\to F_i \to\syz^i(R/\fm)\to 0
\end{equation} 
for any $i\ge 0$, which implies 
$\syz^i(R/\fm)\simeq T^{-i}(R/\fm)$ in $\Dbsing(R)$. 
For $N\in\CMg(R)$ and $i\ge 0$, 
the long exact sequence obtained from eq.(\ref{syz-free}) 
yields 
\begin{equation}\label{derived-LES}
 \begin{array}{cccccccc}
0 &\to & {\rm Hom}_R({\rm syz}^i(R/{\mathfrak m}),N)&\to&
 \Hom_R(F_i,N)&\to&
\Hom_R(\syz^{i+1}(R/\fm),N)& \to\\
 &\to & \Ext^1_R(\syz^i(R/\fm),N) &\to& 0 & & & 
 \end{array}
\end{equation}
and 
\begin{equation}\label{Ext-recursive}
 \Ext^n_R(\syz^{i+1}(R/\fm),N)\simeq
\Ext^{n+1}_R(\syz^i(R/\fm),N), \quad n\ge 1
\end{equation}
since $\Ext^k(F_i,N)=0$ for $k\ge 1$. 
Note that $\syz^i(R/\fm)\in\CMg(R)$ for 
$i\ge d=\dim R$ and then the exact sequence (\ref{syz-free}) 
becomes the one in $\CMg(R)$. 

Now, in the exact sequence (\ref{derived-LES}) 
with $i\ge d$ and $N\in (\T')^\perp$, 
any morphism in $\Hom_R(\syz^{i+1}(R/\fm),N)$ 
factors through a projective-injective object $I$ in 
the Frobenius category $\CMg(R)$.\footnote{In fact, $I\in\CMg(R)$ is 
a projective-injective object if and only if $I$ is graded 
free. }
Moreover, any morphism in $\Hom_R(\syz^{i+1}(R/\fm),I)$ 
factors through $F_i$ with the injection 
$\syz^{i+1}(R/\fm)\to F_i$ 
since $I$ is injective. 
Therefore, any morphism in $\Hom_R(\syz^{i+1}(R/\fm),N)$ 
factors through $F_i$, which implies that 
the map $\Hom_R(F_i,N)\to\Hom_R(\syz^{i+1}(R/\fm),N)$ 
in eq.(\ref{derived-LES}) 
is surjective and hence 
$\Ext^1_R(\syz^i(R/\fm),N)=0$. 

Here, by the isomorphisms (\ref{Ext-recursive}), 
$\Ext^1_R(\syz^i(R/\fm),N)\simeq\Ext^{i+1}_R(R/\fm,N)$ holds 
and hence $\Ext^{i+1}_R(R/\fm,N)=0$ 
for $N\in (\T')^\perp$ and $i\ge d$. 
By Lemma \ref{lem:free}, 
this means that $N$ is graded free, which is isomorphic to zero in 
$\CMgb(R)$. The theorem follows.
\qed\end{pf}
If $R$ defines a hypersurface singularity $A/(f)$, 
then, by the isomorphism of functors $T^2\simeq \tau^h$, 
we see that $D^{gr}_{Sg}(R)$ is of finite type.
\begin{cor}\label{cor:mainI}
Let $\left<E_1,\dots,E_l\right>$ be a full triangulated subcategory 
of $\T=D^{gr}_{Sg}(R)$ generated by an exceptional collection 
$(E_1, \dots, E_l)$ which is closed under the action of $\tau$ and contains 
an object isomorphic to $R/{\mathfrak m}$.
Then $\left<E_1,\dots,E_l\right>\simeq\T$ as 
a triangulated category.
\qed\end{cor}
\begin{rem}
In this paper, we shall apply 
this category generating lemma (Theorem \ref{thm:mainI} 
or Corollary \ref{cor:mainI}) 
together with the Serre functor in Theorem \ref{thm:SerreDual} 
to the corresponding triangulated categories 
associated to the regular systems of weight 
with $\varepsilon_W=-1$ and $a_0=0$ (see subsection \ref{ssec:w}). 
However, these two theorems theirselves can be applied to 
the cases of any $\varepsilon_W$ with $a_0=0$. 
By these theorems, the proof of the main theorem in \cite{KSTade} 
(ADE case: $\varepsilon_W=1$ and $a_0=0$) can be simplified. 
Moreover, the category generating lemma (Theorem \ref{thm:mainI}) 
holds true even if we place the $\Z$-grading by $L(p)$-grading 
as Theorem \ref{thm:SerreDual} does. 
Thus, we can apply these two theorems to $a_0>0$ cases including 
the elliptic cases ($\varepsilon_W=0$ and $a_0=1$), 
which simplifies the proof of the main theorem of \cite{u}. 
\end{rem}

 \section{Strongly exceptional collections in $\T_W$ and 
the associated quivers}
\label{sec:main}

In this section, we formulate our main result 
on the structure of the triangulated category 
($\Dbsing(R_W)\simeq\CMgb(R_W)\simeq\tHMF{A}{f_W}$) 
associated to a regular system of weights $W$ 
of $\varepsilon_W=-1$ and $a_0=0$ 
with a fixed weighted homogeneous polynomial $f_W$. 
In subsection \ref{ssec:w}, we recall the definition 
of the regular system of weights. 
In subsection \ref{ssec:GDD}, we prepare a 
generalized notion of quivers; the quivers associated to 
regular systems of weights $W$ of $\varepsilon_W=-1$ with genus $a_0=0$ 
are defined in this notion in subsection \ref{ssec:DeltaW}. 
Then, in subsection \ref{ssec:main}, 
we state the main theorem of the present paper.

 \subsection{Regular system of weights}\hfill
\label{ssec:w}

In this subsection,  we recall the definition and some of basic facts 
on the regular systems of weights. 
A quadruple $W:=(a,b,c;h)$ of positive integers 
with $a,b,c<h$ and $\mathrm{g.c.d.}(a,b,c)=1$ is called 
a {\em weight system}. 
For a weight system $W$, we define the 
{\em Euler vector field} $E=E_W$ by 
\begin{equation}\label{Euler-vf}
  E:= \frac{a}{h}x\frac{\p}{\p x}+\frac{b}{h}y\frac{\p}{\p y}+
 \frac{c}{h}z\frac{\p}{\p z}\ .
\end{equation}
For a given weight system $W$, 
the regular $\C$-algebra $A=\C[x,y,z]$ becomes a graded ring by putting 
$\deg(x)=2a/h$, $\deg(y)=2b/h$ and 
$\deg(z)=2c/h$. 
Let ${A}=\oplus_{s\in\frac{2}{h}\Z_{\ge 0}} {A}_s$ 
be the graded piece decomposition, 
where ${A}_s:=\{f\in {A}\ |\ 2Ef=s\cdot f\}$\ .
A weight system $W$ is called {\em regular} (\cite{sa:weight}) 
if the following equivalent conditions are satisfied: 
\begin{itemize}
 \item[(a)] $\chi_W(T):=
T^{-h}\frac{(T^h-T^a)(T^h-T^b)(T^h-T^c)}{(T^a-1)(T^b-1)(T^c-1)}$ 
has no poles except at $T=0$. 

 \item[(b)] 
A generic element $f_W$ of the space $A_2=\{f\in A\ |\ Ef= f\}$ 
has an isolated critical point at the origin, 
{\it i.e.}, 
the Jacobi ring 
$A/\left(\frac{\p f_W}{\p x},\frac{\p f_W}{\p y},\frac{\p f_W}{\p z}\right)$ 
is finite rank over $\C$. 

 \item[(c)]
There exists a finite sequence of integers 
$m_1\le m_2\le\cdots\le m_{\mu_W}$ 
for some $\mu_W\in\Z_{>0}$ such that 
the function $\chi_W(T)$ has a Laurent polynomial expansion: 
\begin{equation*}
 \chi_W(T)=T^{m_1}+T^{m_2}+\cdots +T^{m_{\mu_W}}. 
\end{equation*}
\end{itemize}
Here, the number $\mu_W$, called the rank of $W$, 
is given by $(h-a)(h-b)(h-c)/abc$, 
and $\{m_1,\dots,m_{\mu_W}\}$ 
is called the set of {\em exponents} of $W$. 
The smallest exponent $m_1$ is given 
by $\varepsilon_W:=a+b+c-h$. 
An element $f_W$ in $A_2$ as in (b) 
is called a {\em polynomial of type $W$}. 
The quotient ring $R_W:=A/(f_W)$ is Gorenstein, 
whose Gorenstein parameter $\varepsilon(R_W)$ is given 
by the smallest exponent $\varepsilon_W$ because 
the canonical module $K_{R_W}$ is given by the residue 
$Res [A \frac{dx dy dz}{f_W}]=\tau^{-\varepsilon_W}R_W$.

The regular systems of weights are classified as follows. 
Regular systems of weights $W$ with $\varepsilon_W>0$ 
automatically have $\varepsilon_W=1$. 
They are called of type ADE 
via the identification of their exponents 
with those of the root systems of type ADE. 
In our previous paper \cite{KSTade}, 
we studied the triangulated category 
$\tHMF{A}{f_W}$ of this type 
and obtained the root systems of type ADE, as expected. 
Next, regular systems of weights with $\varepsilon_W=0$ are called 
of {\em elliptic} 
since they are associated with simply elliptic singularities. 
There are three such regular systems of weights, which are often 
denoted by $E_6^{(1,1)}$, $E_7^{(1,1)}$ and $E_8^{(1,1)}$ 
according to the classification of elliptic root systems \cite{sa:elliptic}. 

In this paper, we discuss the triangulated category $\tHMF{A}{f_W}$ 
attached to regular systems of weights with $\varepsilon_W=-1$. 
We concentrate on such regular systems of weights with genus $a_0=0$, 
where the genus $a_0$ of $W$ is defined as 
the number of zero exponents of $W$. 
There are twenty two regular systems of weights 
with $\varepsilon_W=-1$ and genus zero, which include fourteen ones, 
so-called, of {\em exceptional unimodular} type 
(see subsection \ref{ssec:strange}).

In order to describe our strongly exceptional collections 
in $\tHMF{A}{f_W}$ for 
a regular system of weights $W$ of $\varepsilon_W=-1$, 
we recall the notion of the signature. 
\begin{defn}[Signature of a regular system of weights 
(\cite{sa:weight} eq.(5.3.2))]
\label{defn:Aw}
For a given regular system of weights $W$, consider 
the following multi-sets of positive integers
\begin{equation*}
 A'_W:= 
 \{a_i\, |\, h/a_i \notin\Z\, ,\, i=1,2,3\} 
\coprod
 \{\gcd(a_i,a_j)^{\left(m(a_i,a_j:h)-1\right)}\, |\, 1\le i<j\le 3\}
\end{equation*}
where $a_1=a$, $a_2=b$, $a_3=c$, and 
$\gcd(a_i,a_j)^{\left(m(a_i,a_j:h)-1\right)}$ indicates that 
we include $\left(m(a_i,a_j:h)-1\right)$ copies of $\gcd(a_i,a_j)$ in 
$A_W'$ with 
$m(a_i,a_j:h):=\sharp\{(u,v)\in(\Z_{\ge 0})^2\, |\, a_iu+a_jv=h \}$. 
We exclude elements equal to one in $A'_W$, 
and denote the result by 
$A_W=(\alpha_1,\dots,\alpha_r)$ for some $r\in\Z_{\ge 0}$, 
where $\alpha_1,\dots,\alpha_r$ are the remaining elements in $A_W$ so that 
$\alpha_1\le\alpha_2\le\cdots\le\alpha_r$. 
The pair $(A_W; a_0)$ is called the {\em signature} of $W$. 
\end{defn}
\begin{defn}[Dual rank $\nu_W$]\label{defn:nuW}
For a given regular system of weights $W$, 
we call the integer 
\begin{equation*}
 \nu_W:=\sum_{i=1}^r(\alpha_i-1)+2(1-a_0)-\varepsilon_W 
\end{equation*}
the {\em dual rank} of $W$. 
\end{defn}
\begin{rem}\label{rem:AW}
Historically, the pair $(A_W, a_0)$ was called 
the signature (Fricke-Klein, Magnus). 
We omit $a_0$ and call $A_W$ the signature of $W$ in this paper, 
since we discuss regular systems of weights $W$ with $\varepsilon_W=-1$ 
and $a_0=0$ only. 

The dual rank $\nu_W$ was originally introduced in \cite{sa:duality} 
as a ``virtual'' rank of the $*$ dual weight 
system $W^*$ of a regular system of weights $W$. 
Since the original formula needs a slight preparation, 
we employ the above formula. The equivalence between those 
two formulas shall be discussed elsewhere.
\end{rem}
For a given $W$, consider a polynomial $f_W$ of type $W$. 
The quotient of the hypersurface $\{(x,y,z)\in\C^3\, |\, f_W=0 \}$ 
by the $\C^\times$-action defined by the weight $(a,b,c)$ 
turns out to be a curve of genus $a_0$ 
having $r$ distinct orbifold points 
$(\lambda_1,\dots,\lambda_r)$ 
with order $A_W=(\alpha_1,\dots,\alpha_r)$. 
See Remark 10.3. of \cite{sa:duality}. 
The relation of these curves with the weighted projective line \cite{GeLe1, GeLe2} 
will be discussed in \cite{KSTwpl}.

 \subsection{Quivers and path algebras of relations}\hfill 
\label{ssec:GDD}

For a triangulated category $\T$, 
assume there exists a strongly exceptional collection 
$\E:=(E_1,\dots,E_l)$. 
We denote 
\begin{equation*}
 \Hom_\T(\E,\E):=\oplus_{i,j=1}^l\Hom_\T(E_i,E_j) 
\end{equation*}
and call it the {\em homomorphism algebra} of $\E$. 
Then, it is known that 
there is a unique quiver such that 
its path algebra with some relations 
is isomorphic to the homomorphism algebra 
$\Hom_\T(\E,\E)$ (see Gabriel \cite{g:2}). 
Although those terminologies are enough for our purpose, 
in this subsection, 
we introduce a modified notion of quivers and give a more explicit rule 
to attach such a quiver to an exceptional collection $\E$. 
In particular, when $\E$ is a strongly exceptional collection, 
we define path algebras of relations of the corresponding quivers.

Let $\T$ be a triangulated category of finite type. 
For any two objects $X, Y\in\T$, the Euler pairing is defined by 
\begin{equation*}
 \chi(X,Y):=\sum_{n\in\Z}(-1)^n\hom_{\T}(X,T^n(Y)). 
\end{equation*}
Suppose there exists an exceptional collection $\E=(E_1,\dots,E_l)$ 
in $\T$. 
The $l$ by $l$ matrix 
\begin{equation*}
 \chi:=\{\chi_{ij}\},\qquad \chi_{ij}:=\chi(E_i,E_j), 
\end{equation*}
is an upper half triangular matrix with $\chi_{ii}=1$ for any 
$i=1,\dots, l$. 
Then, the inverse matrix $C:=\chi^{-1}$ is also 
an upper half triangular matrix with $C_{ii}=1$, 
$i=1,\dots, l$. 
Here we define a quiver associated to such a matrix $C$. 
\begin{defn}\label{defn:quiver}
Let $C=\{C_{ij}\}_{i,j=1,\dots,l}$ be 
an upper triangular $l$ by $l$ matrix 
of integer valued such that $C_{ii}=1$ for any $i=1,\dots,l$. 

The {\em quiver $\vec{\Delta}_{C}=(\Delta_0,\Delta_1;s,e,d)$ 
associated to $C$} is 
the set $\Delta_0=\{1,\dots,l\}$ of vertices and 
the set $\Delta_1$ of arrows with maps 
$s:\Delta_1\to\Delta_0$, $e:\Delta_1\to\Delta_0$, 
$d:\Delta_1\to\{\pm 1\}$ such that 
$s(\rho)\ne e(\rho)$ for any $\rho\in\Delta_1$ and 
\begin{equation*}
 (\sharp\{\rho\in\Delta_1\, |\, (s,e,d)(\rho)=(i,j,+1)\}, 
 \sharp\{\rho\in\Delta_1\, |\, (s,e,d)(\rho)=(i,j,-1)\} ) 
 = \begin{cases}
    (-C_{ij},0) & C_{ij}<0 \\
    (0,C_{ij}) & C_{ij}> 0 \\
    (0,0) & C_{ij}=0 \ 
   \end{cases}
\end{equation*}
for any $i<j$. 
\end{defn}
Now, suppose that we start from a quiver associated to $C$ and 
denote the inverse matrix of $C$ by $\chi:=\{\chi_{ij}\}$. 
By definition one has $\chi_{ii}=1$ 
for any $i$ and $\chi_{ij}=0$ for $i>j$. 

When $\chi_{ij}\ge 0$ for any $i$ and $j$, 
a {\em path algebra with relation} 
of the quiver $\vec{\Delta}_C$ is defined as follows. 
Let $\C\vec{\Delta}_C$ be the path algebra 
defined by arrows $\rho$ of $d(\rho)=+1$. 
For each arrow $\rho$ such that $(s,e,d)(\rho)=(i,j,-1)$, 
a relation $I_{\rho}$ is given 
as a $\C$-linear combination of 
all paths from $i$ to $j$. 
Then, the path algebra with relations 
is the quotient algebra $\C\vec{\Delta}_C/I$, 
where $I$ is the ideal generated by $\{I_\rho\}_{d(\rho)=-1}$. 
\begin{rem}
For a given exceptional collection $\E:=(E_1,\dots,E_l)$ 
in a triangulated category $\T$ of finite type, 
the above procedure actually gives an explicit way to define 
a quiver $\vec{\Delta}_C$, and, in particular, 
if $\E$ is a strongly exceptional collection, 
one can define a path algebra with relations of 
the quiver $\vec{\Delta}_C$. 
However, in general, the homomorphism algebra 
$\Hom_\T(\E,\E)$ of $\E$ is not 
isomorphic to a path algebra with relations of $\vec{\Delta}_C$. 
This is related to that 
Definition \ref{defn:quiver} forbids 
the existence of both arrows $\rho,\rho'\in\Delta_1$ 
such that $(s,e,d)(\rho)=(i,j,+1)$ 
and $(s,e,d)(\rho)=(i,j,-1)$ for the same $i<j$. 
As we shall see in subsection \ref{ssec:main}, 
in our situation, 
this procedure gives an explicit way to give the correct quiver 
associated to an exceptional collection $\E$. 
In particular, the matrix element $\chi_{ij}$ of the inverse matrix 
$\chi$ of $C$ is calculated by 
counting paths from $i$ to $j$, 
consisting of arrows $\rho\in\Delta_1$ of both $d(\rho)=\pm 1$, 
with sign associated to $d=\pm 1$. 
\end{rem}
Hereafter we call a quiver associated to $C$ just a quiver. 
We sometimes drop $C$ 
when we do not give the explicit form of 
the corresponding upper triangular matrix $C$.

 \subsection{Path algebras with relations of 
quivers $\vec{\Delta}_W$, $\vec{\Delta}_W^T$ and 
$\vec{\Delta}_W'$}\hfill
\label{ssec:DeltaW}

\begin{figure}[h]
 \hspace*{-2.0cm} \includegraphics{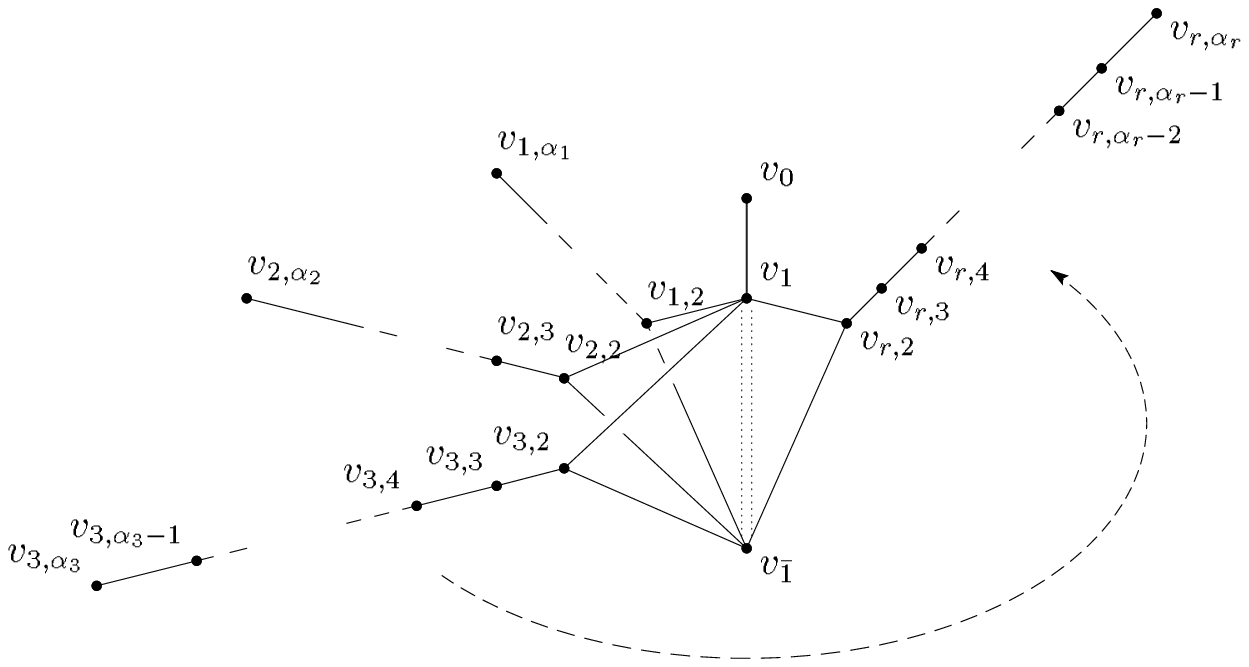}
 \caption{The diagram $\Delta_W=\Delta_{\alpha_1,\dots,\alpha_r}$, 
$A_W=(\alpha_1,\dots,\alpha_r)$, 
obtained from $\vec{\Delta}_W$ or $\vec{\Delta}_W'$ 
by removing the orientation of the arrows. }\label{fig:GD1}
\end{figure}

Now, 
for a regular system of weights $W$ with $\varepsilon_W=-1$ and $a_0=0$, 
we define quivers $\vec{\Delta}_W$, $\vec{\Delta}^T_W$, 
$\vec{\Delta}'_W$ and associated path algebras with relations 
which are necessary to state the main theorem. 
\begin{figure}[h]
 \hspace*{-2.0cm}
 \scalebox{1.0}[1.0]{\includegraphics{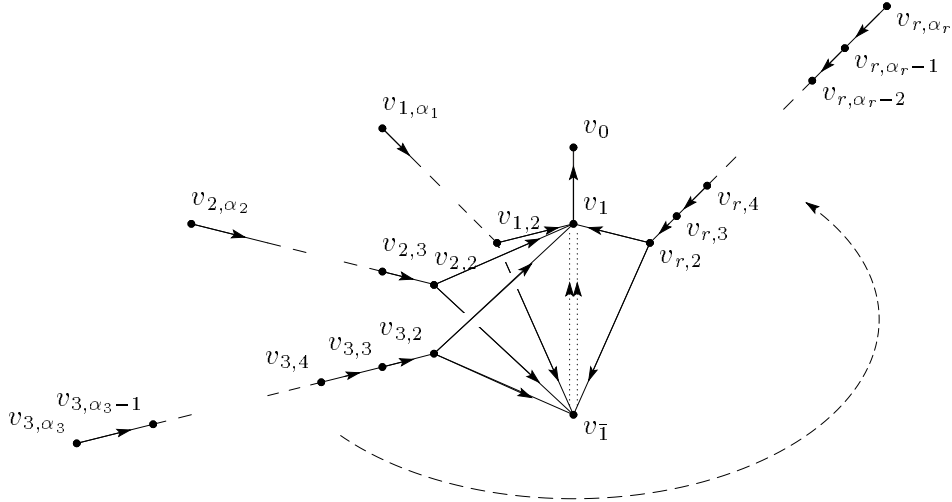}}
 \caption{Figure of $\vec{\Delta}_W$}\label{fig:quiver1}
\end{figure}
\begin{figure}[h]
 \hspace*{-2.0cm}
 \scalebox{1.0}[1.0]{\includegraphics{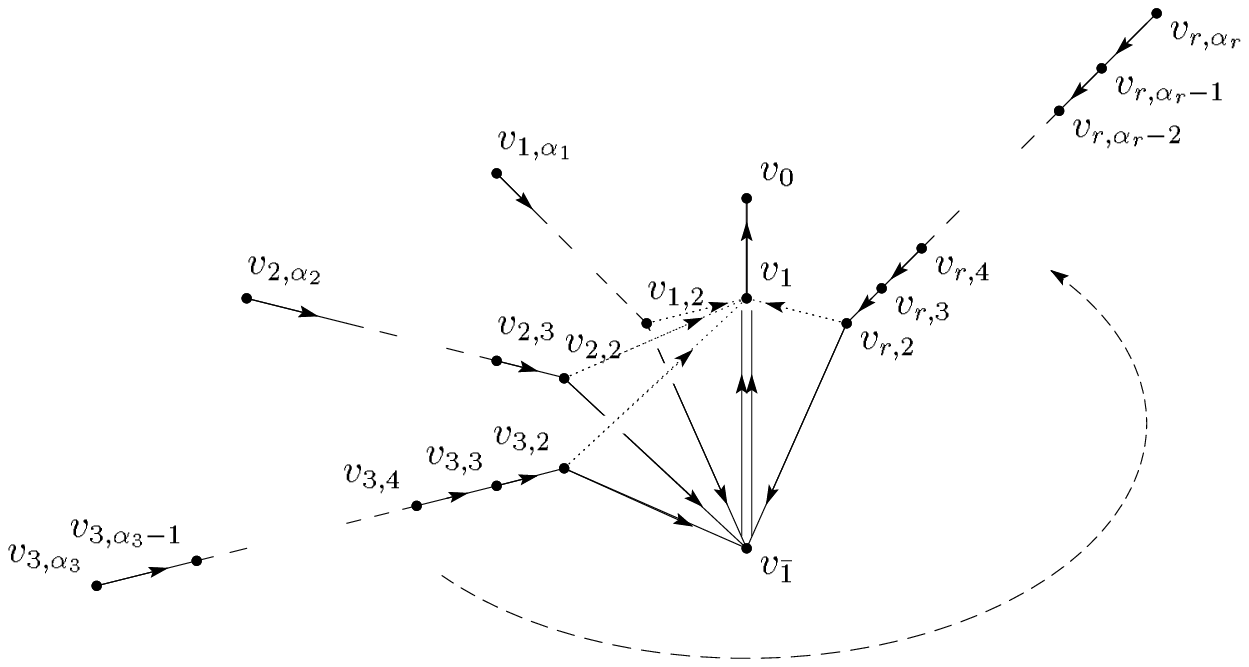}}
 \caption{Figure of $\vec{\Delta}^T_W$}\label{fig:quiver1'}
\end{figure}
\begin{figure}[h]
 \hspace*{-2.0cm}
 \scalebox{1.0}[1.0]{\includegraphics{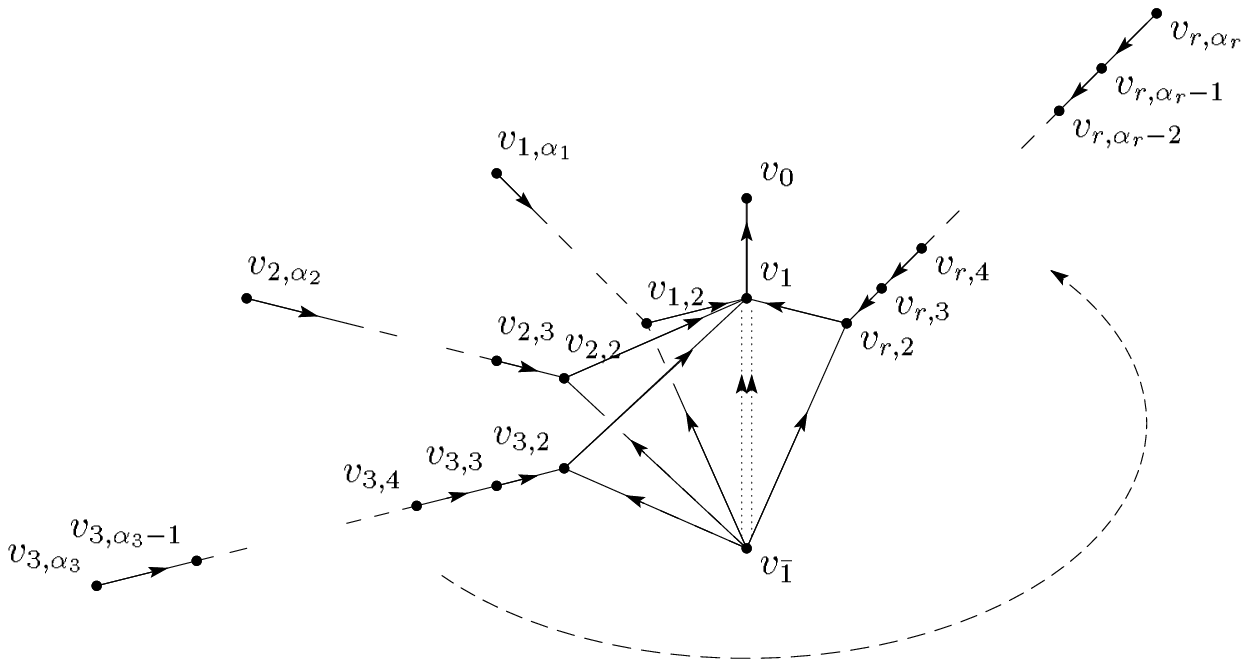}}
 \caption{Figure of $\vec{\Delta}'_W$}\label{fig:quiver2}
\end{figure}
\begin{defn}[Quivers $\vec{\Delta}_W$, $\vec{\Delta}^T_W$, 
$\vec{\Delta}'_W$]\label{defn:DeltaW}
For a regular system of weights $W$ of $\varepsilon_W=-1$ and genus $a_0=0$ 
with signature $A_W=(\alpha_1,\dots,\alpha_r)$, 
we define quivers 
$\vec{\Delta}_W$, $\vec{\Delta}^T_W$ and $\vec{\Delta}'_W$ 
by those in Figure \ref{fig:quiver1}, \ref{fig:quiver1'} and 
\ref{fig:quiver2}, respectively, where 
$\Delta_0=\Pi_W
:=\coprod_{i=1}^r\{v_{i,2},\dots,v_{i,\alpha_i}\}\coprod\{v_0,\vt,\vb\}$ 
is the vertex set which is isomorphic to $\{1,\dots,\nu_W\}$ 
as sets, 
the arrows denote elements $\rho\in\Delta_1$ of $d(\rho)=+1$ 
and the dotted arrows denote elements $\rho\in\Delta_1$ of $d(\rho)=-1$. 
Under the identification $\Pi_W\simeq\{1,\dots,\nu_W\}$, 
for any $v,v'\in\Pi_W$, we sometimes denote by $C(v,v')$ or $\chi(v,v')$ 
the corresponding matrix elements of $C$ or $\chi$. 
\end{defn}
If we forget the orientation of the quiver 
$\vec{\Delta}_W$ or $\vec{\Delta}'_W$, 
we can recover the diagram $\Delta_W$ in Figure \ref{fig:GD1} 
which 
was used to appear as the intersection matrix of the vanishing cycles 
of an exceptional unimodular singularity 
(see subsection \ref{ssec:strange}). 
We sometimes denote this diagram more explicitly by 
$\Delta_{\alpha_1,\dots,\alpha_r}:=\Delta_W$ for 
$A_W=(\alpha_1,\dots,\alpha_r)$. 

Note that, 
for the quiver $\vec{\Delta}_W$, 
$C(v_{i,2},\vb)=C(v_{i,2},v_1)=-1$ for any $i=1,\dots,r$ 
and $C(\vb,v_1)=2$. 
Thus, $\chi(\vb,v_1)=-2$ and also $\chi(v_{i,2},v_1)=-1$. 
On the other hand, for the quiver $\vec{\Delta}_W^T$, 
one has $C(v_{i,2},v_1)=1$ and $C(\vb,v_1)=-2$. 
Then, the inverse matrix $\chi$ of $C$ has non-negative elements 
only. In particular, $\chi(\vb,v_1)=2$ and $\chi(v_{i,2},v_1)=1$. 
For the quiver $\vec{\Delta}_W'$, 
$C(v_{i,2},\vb)=0$ but $C(\vb,v_{i,2})=-1$. 
Then, again the inverse matrix $\chi$ has non-negative elements only. 
In particular, $\chi(\vb,v_{i,2})=\chi(v_{i,2},\vt)=1$ and 
$\chi(\vb,v_1)=r-2$. 

For each of the quivers $\vec{\Delta}^T_W$ and $\vec{\Delta}'_W$, 
we define a path algebra with relations as follows. 
\begin{defn}[
$\PAT$, $\PA'$]
For the quiver $\vec{\Delta}^T_W$, 
we define $r$ relations corresponding to 
the dotted arrows $\rho(v_{i,2},\vt)$ 
from $v_{i,2}$ to $\vt$, $i=1,\dots,r$, by 
\begin{equation}\label{relation1}
 u_{1,i}\cdot(\rho_1(\vb,\vt)\circ \rho(v_{i,2},\vb)) 
 + u_{2,i}\cdot(\rho_2(\vb,\vt)\circ \rho(v_{i,2},\vb)), 
\end{equation}
where $\rho_1(\vb,\vt)$ and $\rho_2(\vb,\vt)$ are two arrows 
from $\vb$ to $\vt$, and 
$[u_{1,i}:u_{2,i}]=\lambda_i$ is a point in $\P^1$. 
We denote by $I_{W,\Lambda}$, $\Lambda:=(\lambda_1,\dots,\lambda_r)$, 
the ideal generated by these relations, 
and by $\PAT$ the corresponding path algebra with relations.  

For the quiver $\quiver'$, 
we define two relations corresponding to 
the dotted arrows $\rho_1(\vb,\vt)$, $\rho_2(\vb,\vt)$ 
from $\vb$ to $\vt$ by 
\begin{equation}\label{relation2}
 \sum_{i=1}^r u_{1,i}\cdot(\rho(v_{i,2},\vt)\circ\rho(\vb,v_{i,2})),
 \qquad 
 \sum_{i=1}^r u_{2,i}\cdot(\rho(v_{i,2},\vt)\circ\rho(\vb,v_{i,2})),  
\end{equation}
where $\rho(v_{i,2},\vt)$ is the arrow from $v_{i,2}$ to $\vt$, 
and $\rho(\vb,v_{i,2})$ is the arrow from $\vb$ to $v_{i,2}$. 
We denote by $I_{W,\Lambda}'$ the ideal generated by these relations and 
by $\PA'$ the corresponding path algebra with relations.  
\end{defn}

 \subsection{The main theorem (Theorem \ref{thm:main})}\hfill
\label{ssec:main}

Recall that $\T_W$ is one of the equivalent 
triangulated categories $\tHMF{A}{f_W}\simeq\Dbsing(R_W)\simeq\CMgb(R_W)$ 
with $R_W=A/(f_W)$ for a fixed polynomial $f_W$ of $W$. 

The following is the main theorem of the present paper. 
\begin{thm}\label{thm:main}
Let $f_W$ be a polynomial of type $W$ with $\varepsilon_W=-1$ 
and $a_0=0$. 
Then, there exist distinct $r$-points 
$\Lambda=(\lambda_1,\dots,\lambda_r)$ in $\P^1$ and 
a full strongly exceptional collection 
$\E^T$ in $\T_W$ such that 
the homomorphism algebra $\Hom_{\T_W}(\E^T,\E^T)$ is isomorphic to 
the path algebra with relations $\PAT$.

The same statement holds true even if $\PAT$ is replaced by 
$\PA'$. 
\end{thm}
This, together with Bondal-Kapranov's theorem \cite{bo:1, bk:1}, 
implies the following: 
\begin{cor}\label{cor:main}
The triangulated category $\T_W$ is triangulated equivalent 
to the derived category 
of finitely generated modules over the algebra $\PAT$ or $\PA'$. 
\end{cor}
Recently, this corollary is proved independently by 
Lenzing and de la Pena \cite{L-delaP} 
in the framework of the weighted projective lines by 
Geigle-Lenzing together with Orlov's arguments in \cite{o:2}.

The main theorem is obtained as a consequence of the 
following structure theorem 
of the triangulated category $\T_W$. 
\begin{thm}\label{thm:key}
Let $f_W$ be a polynomial of type $W$ with $\varepsilon_W=-1$, 
$a_0=0$ and the virtual rank $\nu_W$. 
Then, there exist distinct $r$-points $\Lambda=(\lambda_1,\dots,\lambda_r)$ 
in $\P^1$ and a full exceptional collection 
$\E:=(E_1,\dots, E_{\nu_W})$ in $\T_W$ 
which satisfies the following properties$:$ 
\begin{itemize}
 \item[(i)] For the quiver 
$\vec{\Delta}_W=(\Delta_0=\Pi_W,\Delta_1;s,e,d)$, 
there exists a map $V:\Pi_{W}\to\Ob(\T_W)$ and one has 
\begin{equation*}
\{\E\}:=\{E_1,\dots,E_{\nu_W}\}= 
\coprod_{i=1}^r\{V_{i,2},\dots,V_{i,\alpha_i}\}\coprod\{V_0,\Vt,\Vb\}
\end{equation*}
as sets, 
where 
$V_{i,j}:=V(v_{i,j})$, $V_0:=V(v_0)$, 
$\Vt:=V(\vt)$, $\Vb:=V(\vb)$, 
and $\nu_W$ is the virtual rank of $W$ 
defined in Definition \ref{defn:nuW}. 

 \item[(ii)] For any $v,v'\in\Pi_W$, 
one has $\Hom_\T(V(v),T^n(V(v')))\ne 0$ only if 
$n=0$ or $n=1$. Thus, 
$\chi(V(v),V(v'))=\hom_\T(V(v),V(v'))-\hom_\T(V(v),T(V(v')))$ and 
then 
\begin{equation*}
 \chi(V(v),V(v'))=\chi(v,v') 
\end{equation*}
holds, where $\chi(v,v')$ is the matrix element 
in Definition \ref{defn:DeltaW}. 
We denote by $f_i$, $i=1,\dots,r$, a basis of $\Hom_{\T_W}(V_{i,2},\Vb)$ 
and by $g_i$, $i=1,\dots,r$, a basis of $\Hom_{\T_W}(V_{i,2},T\Vt)$, 
where note that $\chi(v_{i,2},\vt)=\chi(V_{i,2},\Vt)=-1$.

 \item[(iii)]
Under the transpose in Definition \ref{defn:t}, 
the objects in $\{\E\}$ satisfy 
\begin{equation*}
 \begin{split}
 & t(V_{i,j})= \tau^{-(j-1)}(V_{i,j}) ,\\
 & t(V_0)= T\tau^{-1}(V_0) ,\\
 & t(\Vt)= T(\Vt) ,\\
 & t(\Vb)= \tau^{-1}(\Vb') ,
 \end{split}
\end{equation*}
where $T(\Vb')\in\Ob(\T)$ is the cone of the morphisms 
$\oplus_{i=1}^r f_i: \oplus_{i=1}^r V_{i,2}\to\Vb$.

 \item[(iv)] For each element ${V}\in\{\E\}$, 
the AR-triangle is given as follows. 
Recall that $\tau_{AR}=\tau$. 

(iv-a)\ For $V_{i,j}\in\{\E\}$, $i=1,2,\dots, r$, $j=2,3,\dots, \alpha_i$, 
\begin{equation}\label{AR1}
 \tau_{AR}(V_{i,j})\to\tau V_{i,j-1}\oplus V_{i,j+1} \to V_{i,j} \to 
 T\tau_{AR}(V_{i,j}) ,
\end{equation}
where we put $V_{i,\alpha_i+1}=0$ 
for $j=\alpha_i$ and $V_{i,1}$ is defined by 
the above triangle with $j=2$.  

(iv-a')\ For $V_0\in\{\E\}$, 
\begin{equation}\label{AR2}
  \tau_{AR}(V_0)\to\Vt\to V_0 \to T\tau_{AR}(V_0), 
\end{equation}

(iv-b)\ For $\Vt,\Vb\in\{\E\}$,  
\begin{equation}\label{AR3}
\tau_{AR}(\Vt)\to 
 V_2\oplus \tau V_0\to\Vt\to T\tau_{AR}(\Vt),
\end{equation}
\begin{equation}\label{AR3'}
\tau_{AR}\Vb\to 
AR(\Vb)\to\Vb\to T\tau\Vb. 
\end{equation}
Here, in eq.(\ref{AR3}),  
$T(V_2)$ is defined by the cone of 
$(V_{1,2}\oplus V_{2,2}\oplus\cdots\oplus V_{r,2})\to 
(\Vb)^{\oplus 2}$ with the morphisms given by 
\begin{equation*}
u_{1,i}\cdot f_i\oplus u_{2,i}\cdot f_i: V_{i,2}\to\Vb\oplus\Vb
\end{equation*}
for $\lambda_i=[u_{1,i}:u_{2,i}]\in\P^1$ 
and $f_i$, $i=1,\dots,r$, a basis of $\Hom_{\T_W}(V_{i,2},\Vb)$. 
In eq.(\ref{AR3'}), $AR(\Vb)$ is defined by the cone of 
$(\tau\Vt)^{\oplus 2}\to
(V_{1,2}\oplus V_{2,2}\oplus \cdots\oplus V_{r,2})$
with the morphism given by 
\begin{equation*}
u_{1,i}\cdot \tau t(g_i)\oplus u_{2,i}\cdot \tau t(g_i): 
(\tau\Vt)^{\oplus 2}\to V_{i,2}
\end{equation*}
for $g_i$, $i=1,\dots,r$, a basis of $\Hom_{\T_W}(V_{i,2},T\Vt)$, 
where 
recall that $t$ is the transpose in Definition \ref{defn:t}. 

\item[(iv')]
There exists a triangle 
\begin{equation}\label{coneXX*}
   \Vt\to V_{i,1}\to\Vb\to T\Vt 
\end{equation}
for each $i=1,2,\dots r$, 
where the morphism $\Vb\to T\Vt$ is described as 
$u_{1,i}e_1+ u_{2,i}e_2$ in terms of a basis $\{e_1, e_2\}$ 
of $\Hom_{\T_W}(\Vb,T(\Vt))$. 
\end{itemize}
\end{thm}

\vspace*{0.3cm}

The AR-triangles (\ref{AR1}) and (\ref{AR2}) imply that 
the morphism $V_{i,j}\to V_{i,j-1}$ for each $i$ and $j=2,\dots,\alpha_i$ 
and the morphism $\Vt\to V_0$, respectively, 
are irreducible.  
Furthermore, by $\chi(\vb,\vt)=-2$, one has 
$\hom_{\T_W}(\Vb,T\Vt)=\hom_{\T_W}(T^{-1}(\Vb),\Vt)=2$ and 
can consider the cone of the morphisms 
$T^{-1}(\Vb)\to\Vt$ 
parameterized by $\P^1=\P(\Hom_{\T_W}(\Vb,T\Vt))$. 
Then, the triangle (\ref{coneXX*}) implies that 
this $\P^1$ has $r$ special points $\lambda_i=[u_{1,i}:u_{2,i}]$, 
$i=1,\dots,r$, which correspond to $V_{i,1}$.

As we shall discuss in the proof of the main theorem 
in subsection \ref{ssec:proof-main}, 
Theorem \ref{thm:key} implies the main theorem (Theorem \ref{thm:main}) 
since 
we can obtain strongly exceptional collections 
corresponding to quivers 
$\vec{\Delta}^T_W$ and $\vec{\Delta'}_W$ 
and the AR-triangles for the collections. 
Also, we can replace quiver 
$\vec{\Delta}^T_W$ or $\vec{\Delta'}_W$ 
by the one whose orientation of the arrows between 
$V_{i,j}$ and $V_{i,j+1}$ are taken arbitrary, and then 
the parallel statement to the main theorem holds true. 
The parallel statement further holds true even if 
the orientation of the arrows of these quivers is reversed, 
since the contravariant functor $t:\T_W\to \T_W$ is an automorphism on $\T_W$. 
\begin{rem}\label{rem:key}
For a full exceptional collection $\E=(E_1,\dots,E_l)$ 
in a triangulated category $\T$, let $C$ be the inverse matrix of 
$\{\chi_{ij}\}_{i=1,\dots,l}$, where $\chi_{ij}=\chi(E_i,E_j)$ is the 
Euler number. 
Then, it is known \cite{bo:1, bp:1} that the AR-translation $\tau_{AR}$ 
induces an isomorphism $[\tau_{AR}]$ on the Grothendieck group of $\T$, 
which is expressed in terms of the basis $[E_1],\dots,[E_l]$ as 
\begin{equation}\label{AR-coxeter}
 [\tau_{AR}]([E_1],\dots,[E_l]) = -([E_1],\dots,[E_l])\, 
C\cdot ^t\!C^{-1}, 
\end{equation}
where $-C\cdot ^t\!C^{-1}$ is the Coxeter transformation originally 
defined as the product of the reflections associated to each root base 
of the root lattice $(K_0(\T),\chi+{}^t\chi)$. 
In our case $\T=\T_W$, by Corollary \ref{cor:tau-AR}, 
$[(\tau_{AR})^h]=\Id$ holds, 
which implies that the Coxeter transformation 
on the Grothendieck group of $\T_W$ is of finite order 
and in particular of order $h$.

On the other hand, we can see that 
the AR-triangles in Theorem \ref{thm:key} reduce to 
the identity (\ref{AR-coxeter}) 
at the level of the Grothendieck group. 
First, the reduction of 
the AR-triangles in Theorem \ref{thm:key} 
gives the following identity 
\begin{equation}
[\tau_{AR}(E_i) ]+[E_i]= 
\sum_{j=1}^{\nu_W} 
\left( 
-C_{ij}\cdot [\tau_{AR}(E_j)]
-C_{ji}\cdot [E_j]\right) 
\end{equation}
for any $i=1,\dots,\nu_W$. 
It is easy to see that this is equivalent to the identity (\ref{AR-coxeter}). 
\end{rem}

 \subsection{A categorification of the strange duality}\hfill
\label{ssec:strange}

There exist 
fourteen regular systems of weights with $\varepsilon_W=-1$ and $a_0=0$ 
such that $r=3$ for the signature $A_W=(\alpha_1,\dots,\alpha_r)$. 
The singularity defined by a generic polynomial $f_W\in A_2$ 
has been called 
an exceptional unimodular singularity, where the signature $A_W$ 
coincides with the {\em Dolgachev numbers} \cite{arnold,d:1}. 

For any exceptional unimodular singularity, 
the intersection diagram of a distinguished basis for 
the vanishing cycles of the Milnor fiber is given 
\cite{gv:1} (see also \cite{ew,eb:1,gv:2}) 
as $\Delta_{\beta_1,\beta_2,\beta_3}$ (see Figure \ref{fig:GD1}) 
with some triple of positive integers 
$B_W:=(\beta_1,\beta_2,\beta_3)$. 
This triple $B_W$ is called 
the {\em Gabrielov numbers} of the singularity of the polynomial $f_W$. 

The {\em strange duality}, found by Arnord, is a duality 
between these fourteen exceptional unimodular singularities 
stating the existence of an exceptional unimodular singularity $(f_W)^*$, 
that is, a regular system of weights $W^*$ of exceptional type, 
satisfying 
\begin{equation*}
 A_W=B_{W^*} ,\qquad B_W=A_{W^*}
\end{equation*}
for any $W$ of exceptional unimodular type such that $(W^*)^*=W$. 

This duality is now interpreted in various ways: 
Kawai-Yang \cite{kawai-yang} explained this strange duality 
in terms of the duality of 
orbifoldized Poincar\'e polynomials, i.e., 
the topological mirror symmetry. 
On the other hand, the $*$ duality, 
introduced as a duality of weight systems 
in \cite{sa:duality}, includes the strange duality. 
Though the $*$ duality was originally defined 
in terms of the characteristic polynomials of the Milnor monodromy, 
the relation of it 
with the topological mirror symmetry is also discussed in \cite{t:1}. 

Now, Theorem \ref{thm:main} gives 
another interpretation, say, 
a categorification of the strange duality 
at the level of the Grothendieck group: 
\begin{cor}\label{cor:strange}
The following isomorphism of abelian group holds$:$
\begin{equation*}
 K_0(\T_W)\simeq (H_2(f^{-1}_{W^*}(1),\Z),-I_{W^*}). 
\end{equation*}
\end{cor}
This can be thought of the {\em homological} 
mirror symmetry at the level of the Grothendieck group. 
In order to discuss this kind of duality 
at the level of triangulated categories, 
we need to define a suitable Fukaya category 
of the vanishing cycles of the Milnor fiber, 
which we leave for one of future directions.

 \section{Proof of Theorems \ref{thm:main} and \ref{thm:key}}
\label{sec:proof}

 \subsection{Proof of Theorem \ref{thm:key}}\hfill
\label{ssec:proof-key}

In this subsection, we give a proof of Theorem \ref{thm:key} 
in the following order. 
\begin{itemize}
 \item 
We first give a way to construct a collection 
$\{\E\}:=
\coprod_{i=1}^r\{V_{i,2},\dots,{V}_{i,\alpha_i}\}\coprod\{\Vt,\Vb,V_0\}$ 
of indecomposable objects 
which have the grading matrices listed as in section \ref{sec:list} 
so that $\{\E\}$ has AR-triangles (\ref{AR1}), (\ref{AR2}), (\ref{AR3}), 
(\ref{AR3'}) and the triangle (\ref{coneXX*}). 
Thus, Statements (i) and (iv) are completed there. 

 \item 
Using Lemma \ref{lem:hom-calculation} 
on the existence of morphisms, 
we show Statement (ii), which implies that 
$\{\E\}$ forms an exceptional collection $\E$, 
and also complete Statement (iv'). 

 \item
Statement (iii) is then clear by construction except for 
$\Vb$. For $\Vb$, we show Statement (iii) 
from the AR-triangle (\ref{AR3'}).  

 \item 
Finally, we show that the exceptional collection $\E$ is full 
by using Statements (i), (ii), (iv) and (iv') together with 
Corollary \ref{cor:mainI}.  
\end{itemize}

 \vspace*{0.4cm}

\noindent
{\bf The construction of $\{\E\}$:}\quad
We first find 
a candidate for ${V}_0$ and ${V}_{i,\alpha_i}$, $i=1,\dots,r$, 
by hands, which are listed in section \ref{sec:list}. 
We choose ${V}_0$ so that it is isomorphic to $R/\fm$ 
in $\Dbsing(R_W)$ up to grading shifts. 
Then, $\Vt$ is obtained as $AR(V_0)$. 
Also, given ${V}_{i,\alpha_i}$, 
we obtain 
${V}_{i,\alpha_i-1}=AR(\tau^{-1}{V}_{i,\alpha_i})$, 
${V}_{i,\alpha_i-2}\oplus\tau({V}_{\alpha_i})=AR(\tau({V}_{i,\alpha_i}))$, 
and repeating this procedure 
yields ${V}_{i,2},\dots,{V}_{i,\alpha_i}$ 
and also ${V}_{i,1}$. 
The remaining object we should find is then $\Vb$. 
It is obtained as the cone of $\Vt\to V_{i,1}$ for some $i=1,\dots,r$.
In fact, we have the isomorphic object for any $i$. 

By construction, the collection $\{\E\}$ satisfies the 
AR-triangles (\ref{AR1}) and (\ref{AR2}). 
The existence of the triangle (\ref{coneXX*}) also 
follows from the construction above. 
We shall discuss on the morphism $\Vb\to T\Vt$ 
in the triangle (\ref{coneXX*}) later. 
The AR-triangles (\ref{AR3}) and (\ref{AR3'}) can be checked 
by direct calculations. 

At the level of grading matrices (cf. section \ref{sec:list}), 
all these AR-triangles can be checked easily. 

Now, Statements (i) and (iv) are completed.

\vspace*{0.3cm}

\noindent
{\bf On the structure of $\{\E\}$: }\quad 
We shall need some explicit data of morphisms 
between two graded matrix factorizations. 
For this purpose, we introduce the notion of phase. 
\begin{defn}[Phase of an indecomposable graded matrix factorization]
A graded matrix factorization $\olF\in\tHMF{A}{f_W}$ is called {\em reduced} 
if it has minimal rank in its isomorphism class in $\tMF{A}{f_W}$. 
For an indecomposable graded matrix factorization $\olF\in\tHMF{A}{f_W}$, 
the {\em phase} $\phi(\olF)$ of $\olF$ is defined by 
\begin{equation*}
 \phi(\olF):= \ov{\rank(Q,S)}\Tr(S), 
\end{equation*}
where $(Q,S)$ is a reduced graded matrix factorization 
which is isomorphic to $\olF$ in $\tHMF{A}{f_W}$ and is reduced. 
\end{defn}
Note that, for an indecomposable object $\olF\in\tHMF{A}{f_W}$, 
the phase $\phi(\olF)$ does not depend on the choice of $(Q,S)$. 
\begin{defn}\label{defn:sp}
For two indecomposable objects $\olF,\olFF\in\tHMF{A}{f_W}$, 
we denote $\phi(\olF,\olFF):=\phi(\olFF)-\phi(\olF)$. 
Define the {\em spectrum} $\sp(\olF,\olFF)$ 
by the following multi-set of rational numbers: 
\begin{equation*}
 \sp(\olF,\olFF)
 :=\{\phi(\olF,\tau^n(\olFF))^{\hom_{\tHMF{A}{f_W}}(\olF,\tau^n(\olFF))}\ |\ n\in\Z\}, 
\end{equation*}
where $\phi(\olF,\tau^n(\olFF))^{\hom_{\tHMF{A}{f_W}}(\olF,\tau^n(\olFF))}$ 
indicates that we include $\hom_{\tHMF{A}{f_W}}(\olF,\tau^n(\olFF))$ 
copies of $\phi(\olF,\tau^n(\olFF))$ in $\sp(\olF,\olFF)$. 
In particular, we denote $\sp(\olF,\olF)=\sp(\olF)$. 
\end{defn}
By definition, for any two indecomposable objects 
$\olF,\olFF\in\tHMF{A}{f_W}$, 
one has 
\begin{equation*}
 \begin{split}
 &\sp(\olF,\tau(\olFF))=\sp(\tau(\olF),\olFF)=\sp(\olF,\olFF),\\ 
 &\sp(t(\olFF),t(\olF))=\sp(\olF,\olFF), \\
 &\sp(T(\olF),T(\olFF))=\sp(\olF,\olFF). 
 \end{split}
\end{equation*}
However in general 
$\sp(\olF,T(\olFF))\ne \sp(\olF,\olFF)$; instead of it, 
by the Serre duality, the following holds: 
\begin{equation*}
\sp(\olFF,T(\olF))=\sp(\olFF,\S(\olF))=
\left\{\left(1-\frac{2\varepsilon_W}{h}\right)-p\ \Big|\ 
 p\in\sp(\olF,\olFF)\right\}
\end{equation*}
(with $\varepsilon_W=-1$). 
\begin{lem}\label{lem:hom-calculation}
Given a triangulated category $\T_W$ 
of a regular system of weights 
$W=(a,b,c;h)$ with $\varepsilon_W=-1$, $a_0=0$, 
let the spectrum 
$\sp(V,V')$, $V,V'\in\T_W$, 
be $\{p_0\le p_1\le\cdots\le p_k\}$, $i=0,1,\dots,k$, 
for some $k\in\Z_{\ge 0}$. 
Then, one has 
\begin{equation*}
 0\le p_0\le \cdots\le p_k\le 1-\frac{2\varepsilon_W}{h},
\end{equation*}
for any 
$V,V'\in \coprod_{i=1}^r\{V_{i,\alpha_i}\}\coprod\{V_0\}$. 
In particular 
\begin{itemize}
 \item[(a)] $\sp(V_0)=\{0\le 2(a-\varepsilon_W)/h
 \le (b-\varepsilon_W)/h\le 2(c-\varepsilon_W)/h\}$, 
 \item[(b)] for $\sp(V_{i,\alpha_i})$, $p_0=0$ and $p_1=2\alpha_i/h$, 
\item[(c)] 
$\sp(V_0,V_{i,\alpha_i})=\{p_0\le\cdots < 1/2-\alpha_i/h < 
1/2+(\alpha_i+2)/h<\cdots\le p_k\}$, \\
where $(2\alpha_i+2)/h-1/2<p_0$, 
 \item[(d)] for $\sp(V_{i,\alpha_i}, V_{j,\alpha_j})$, $i\ne j$, 
$p_0=(\alpha_i+\alpha_j-2\varepsilon_W)/h$. 
\end{itemize} 
\qed\end{lem}
A few remarks about (c) are in order. 
The rational number $1/2-\alpha_i/h$ is always greater than zero for any 
regular system of weights $W$ with $\varepsilon_W=-1$ and $a_0=0$. 
By the transpose of $V_0$ and $V_{i,\alpha_i}$, 
$\sp(V_0,V_{i,\alpha_i})=\sp(V_{i,\alpha_i},TV_0)$ holds. 
Thus, by the Serre duality, the condition 
$(2\alpha_i+2)/h-1/2<p_0$ is equivalent to 
that $p_k< 3/2-2\alpha_i/h$. 

Now, we calculate the dimension of the space of morphisms 
to show Statement (ii) and to complete to show Statement (iv'). 
Notice that, by construction, 
the phase of $V_{i,j}$ is $\phi(V_{i,j})=(j-1)/h$ for any $i=1,\dots,r$ and 
$j=1,\dots,\alpha_i$. 
On the other hand, the phases of $\Vt$ and $V_0$ are 
$\phi(\Vt)=-1/2$ and $\phi(V_0)=-(1/2)-(1/h)$. 
(See subsection \ref{ssec:MF-grading}). 
Then, for all 
$V,V'\in\coprod_{i=1}^r\{V_{i,1},\dots,V_{i,\alpha_i}\}\coprod\{V_0,\Vt\}
\ ( =(\{\E\}\backslash \{\Vb\})\coprod (\coprod_{i=1}^r\{V_{i,1}\}))$, 
we can calculate $\hom_{\T_W}({V},T^n{V'})$ with any $n\in\Z$ 
due to Corollary \ref{cor:AR}. 
Consequently, we obtain the followings: 
for $i,i'\in\{1,\dots,r\}$, $j\in\{1,\dots,\alpha_i\}$, 
$j'\in\{1,\dots,\alpha_{i'}\}$, and $k,k'\in\{0,1\}$, 
\begin{equation}\label{hom-V_k}
 \hom_{\T_W}(V_k,T^n(V_{k'}))=
 \begin{cases}
 1 & n=0\ \text{and}\ k\ge k' \\
 0 & \text{otherwise} 
 \end{cases},
\end{equation}
\begin{equation}\label{hom-V_ii'}
 \hom_{\T_W}({V}_{i,j},T^n({V}_{i',j'}))=
 \begin{cases}
 1 & (n=0, \ i=i',\ j\ge j')\ \text{or}\ (n=1,\ i=i',\ j'=1,\ \forall j)  \\
 0 & \text{otherwise}
 \end{cases},
\end{equation}
\begin{equation}\label{hom-V_i1}
\hom_{\T_W}({V}_{i,j},T^n(V_k))=
 \begin{cases}
 1 & n=1 \ \text{and any}\ i,j,k \\
 0 & \text{otherwise} 
 \end{cases}, 
\end{equation}
\begin{equation}\label{hom-V_1i}
 \hom_{\T_W}(V_k,T^n({V}_{i,j}))=
 \begin{cases}
 1 & n=0\ \text{and any}\ i, \ \ j=1,\ \ k=1 \\
 0 & \text{otherwise}
 \end{cases}.
\end{equation}
The equation (\ref{hom-V_k}) follows from Lemma \ref{lem:hom-calculation} (a). 
The equations (\ref{hom-V_ii'}) and (\ref{hom-V_i1}) 
follow from Lemma \ref{lem:hom-calculation} (c). 
The equations (\ref{hom-V_1i}) with $i=i'$ and $i\ne i'$ follow from 
Lemma \ref{lem:hom-calculation} (b) and (d), respectively. 
For $k=1$ and $j=1$, 
eq.(\ref{hom-V_i1}) and eq.(\ref{hom-V_1i}) 
are equivalent under the transpose $t$. 
The equation (\ref{hom-V_ii'}) implies 
that $V_{i,1}$ and $T^n(V_{i',1})$ are not isomorphic to each other 
for any $n$ if $i\ne i'$.

By eq.(\ref{hom-V_k}) with $k=k'$ and 
eq.(\ref{hom-V_ii'}) with $i=i'$ and $j=j'$, 
we can see that 
all elements in $\{\E\}\backslash\{\Vb\}$ 
are exceptional objects. 
Furthermore, for any 
$V\in\coprod_{i=1}^r\{V_{i,1},\dots,V_{i,\alpha_i}\}\coprod\{V_0,\Vt\}$, 
since $\hom_{\T_W}(V,V)=1=\hom_{\T_W}(V,\S(V))$, 
the spectrums satisfy the following rules
\begin{equation*}
 \begin{split}
 & \sp(V,AR(V'))=\{p-1/h\ |\ p\in\sp(V,V'),\ p\ne 0\}\coprod
 \{p+1/h\ |\ p\in\sp(V,V'),\ p\ne 1+2/h\}, \\
 & \sp(AR(V),V')=\{p-1/h\ |\ p\in\sp(V',V),\ p\ne 0\}\coprod
 \{p+1/h\ |\ p\in\sp(V',V),\ p\ne 1+2/h\}  
 \end{split}
\end{equation*}
for any 
$V,V'\in\coprod_{i=1}^r\{V_{i,1},\dots,V_{i,\alpha_i}\}\coprod\{V_0,\Vt\}$, 
which are obtained by rewriting Corollary \ref{cor:AR} directly. 

The remaining thing to show Statement (ii) is to 
calculate $\hom_{\T_W}(V,V')$ for the case 
$V=\Vb$ and/or $V'=\Vb$. 
First, by applying the functor $\Hom_{\T_W}(\, \cdot\, ,\Vt)$ 
to the triangle (\ref{coneXX*}), we get 
\begin{equation*}
 \hom_{\T_W}(\Vb,T^n\Vt)= 
 \begin{cases}
 2 & n=1 \\
 0 & \text{otherwise} 
 \end{cases}, 
\end{equation*}
where we use eq.(\ref{hom-V_i1}) with $j=1$, $k=1$, 
and eq.(\ref{hom-V_k}) with $k=k'=1$. 
This implies that the morphism 
$\Vb\to T\Vt$ in the triangle (\ref{coneXX*}) 
is described as that stated in Theorem \ref{thm:key} 
since $V_{i,1}$ and $V_{i',1}$ are not isomorphic 
to each other for $i\ne i'$. 

By applying the functor $\Hom_{\T_W}(\Vt,\, \cdot\, )$ 
to the triangle (\ref{coneXX*}), 
we get 
\begin{equation}\label{hom-V_1b}
 \hom_{\T_W}(T^n\Vt,\Vb)=0 \qquad \text{for any}\ \ n\in\Z, 
\end{equation}
where we use eq.(\ref{hom-V_1i}) with $j=1$ and eq.(\ref{hom-V_k}) 
with $k=k'=1$. 

In a similar way, by applying $\Hom_{\T_W}(\, \cdot\, ,{V}_0)$ and 
$\Hom_{\T_W}({V}_0,\, \cdot\, )$ to the triangle (\ref{coneXX*}) 
we have 
\begin{equation*}
 \hom_{\T_W}(T^{-n}(\Vb),V_0)= 
 \begin{cases}
 2 & n=1 \\
 0 & \text{otherwise} 
 \end{cases}
\end{equation*}
and 
\begin{equation*}
\hom_{\T_W}(V_0,T^n(\Vb))=0
\qquad \text{for any}\ \ n\in\Z. 
\end{equation*}
Here, the former equation follows from that 
$\hom_{\T_W}(T^{-n}(\Vt),V_0)$ is equal to one for $n=0$ and zero otherwise 
by eq.(\ref{hom-V_k}), and that 
$\hom_{\T_W}(T^{-n}(V_{i',1}),V_0)$ is equal to one for $n=1$ 
and zero otherwise 
by eq.(\ref{hom-V_i1}). 
The latter equation follows from 
that $\hom_{\T_W}(V_0,T^n(\Vt))=
\hom_{\T_W}({V}_0,T^n({V}_{i,1}))=0$ for any $n\in\Z$ 
by eq.(\ref{hom-V_k}) and eq.(\ref{hom-V_1i}), respectively. 

In order to compute $\hom_{\T_W}({V}_{i,2},T^n(\Vb))$, 
we apply the functor $\Hom_{\T_W}({V}_{i,2},\, \cdot\, )$ 
to the triangle (\ref{coneXX*}): 
\begin{equation}\label{coneXX*2}
 \Vt \to {V}_{i',1} \to \Vb \to T\Vt
\end{equation}
for $i'\ne i$. 
The resulting exact sequence is 
\begin{equation*}
 \begin{split}
 \cdots &\to\Hom_{\T_W}({V}_{i,2},\Vt)
 \to\Hom_{\T_W}({V}_{i,2},{V}_{i',1})\to\Hom_{\T_W}({V}_{i,2},\Vb) 
 \to\Hom_{\T_W}({V}_{i,2},T\Vt) \\
 & \to\Hom_{\T_W}({V}_{i,2},T({V}_{i',1}))\to\Hom_{\T_W}({V}_{i,2},T(\Vb))
 \to\Hom_{\T_W}({V}_{i,2},T^2(\Vt))\to\cdots, 
 \end{split}
\end{equation*}
where, $\hom_{\T_W}({V}_{i,2},T^n({V}_{i',1}))=0$ 
by eq.(\ref{hom-V_ii'}) with $j=2$, $j'=1$, and 
$\hom_{\T_W}({V}_{i,2},T^n(\Vt))$ is equal to one for $n=1$ and zero otherwise 
by eq.(\ref{hom-V_i1}) with $j=2$ and $k=1$. 
This implies 
\begin{equation*}
 \hom_{\T_W}({V}_{i,2},T^n(\Vb))=
\begin{cases}
 1 & n=0 \\
 0 & \text{otherwise}
\end{cases}. 
\end{equation*}
Similarly, for $j>2$, 
by considering the functor $\Hom_{\T_W}({V}_{i,j},\, \cdot\,)$ 
we obtain that 
$\hom_{\T_W}({V}_{i,j},T^n(\Vb))$ is equal to one for $n=0$ 
and zero otherwise. 
Conversely, applying the functor 
$\Hom_{\T_W}(\, \cdot\, ,{V}_{i,j})$ with $j\ge 2$ 
to the triangle (\ref{coneXX*2}) with $i'=i$ 
leads to the result $\hom_{\T_W}(T^n(\Vb),{V}_{i,j})=0$ 
for any $n\in\Z$, 
since 
$\hom_{\T_W}(T^n(\Vt),{V}_{i,j})=\hom_{\T_W}(T^n({V}_{i,1}),{V}_{i,j})=0$ 
for any $n\in\Z$ by 
eq.(\ref{hom-V_1i}) and eq.(\ref{hom-V_ii'}), respectively.

Finally, $\hom_{\T_W}(\Vb,T^n(\Vb))$ is computed as follows. 
Applying $\Hom_{\T_W}(V_{i,1},\, \cdot\, )$ to the triangle (\ref{coneXX*2}) 
with $i\ne i'$, we obtain 
\begin{equation}\label{hom-V_i1-b}
 \hom_{\T_W}(V_{i,1},T^n(\Vb))=
 \begin{cases}
 1 & n=0 \\
 0 & \text{otherwise} 
 \end{cases}
\end{equation}
since $\hom_{\T_W}(V_{i,1},T^n(\Vt))$ is equal to one for $n=1$ 
and zero otherwise 
(eq.(\ref{hom-V_i1}) with $j=1$ and $k=1$), 
and $\hom_{\T_W}(V_{i,1},T^n(V_{i',1}))=0$ 
(eq.(\ref{hom-V_ii'}) with $j=1$ and $j'=1$). 
Then, applying $\Hom_{\T_W}(\,\cdot\, ,\Vb)$ to the triangle (\ref{coneXX*}) 
yields 
\begin{equation*}
 \hom_{\T_W}(\Vb,T^n(\Vb))=
 \begin{cases}
 1 & n=0 \\
 0 & \text{otherwise} 
 \end{cases}
\end{equation*}
due to eq.(\ref{hom-V_1b}) and eq.(\ref{hom-V_i1-b}). 

Statement (ii) has been completed.

 \vspace*{0.3cm}

Statement (iii), that is, 
the properties under the transpose can be checked for 
$V_0$ and $V_{i,\alpha_i}$, $i=1,\dots,r$, 
by the explicit form of the graded matrix factorizations and 
for $V_{i,j}$, $j=2,\dots,\alpha_i-1$, by construction above. 

The statement $t(\Vb)=\tau^{-1}(\Vb')$ is related to the AR-triangle 
(\ref{AR3'}) via the octahedral axiom of triangulated categories. 
Recall that, for a triangulated category $\T$, the octahedral axiom 
states the existence of the triangle 
$Z'\to Y'\to X'\to T(Z')$ 
(with appropriate compatibility of morphisms) 
for any given objects $X,Y,Z\in\T$ and 
a composition $X\to Y\to Z$ of morphisms, 
where $X'$, $Y'$ and $Z'$ are defined by the triangles 
$X\to Y\to Z'\to T(X)$, $Y\to Z\to X'\to T(Y)$ and $X\to Z\to Y'\to T(X)$, 
respectively. 
Now, apply the octahedral axiom for the case $\T=\T_W$ with 
$X=\oplus_{i=1}^r T^{-1}(V_{i,2})$, $Y=T^{-1} AR(\Vb)$ and 
$Z=T^{-1}(\Vb)$. Then, one obtains the triangle 
$Z'\to Y'\to X'\to T(Z')$ with $Z'\simeq (\tau(\Vt))^{\oplus 2}$, 
$Y'\simeq\Vb'$ and $X'\simeq \tau(\Vb)$. 
Namely, $\Vb'$ is the mapping cone 
$\Vb'\simeq \tau C(T^{-1}\Vb\to\Vt^{\oplus 2})$. 
Next, apply the octahedral axiom again to the case 
$X=T^{-1}\Vb$, $Y=\Vt^{\oplus 2}$ and $Z=\Vt$, 
where $T^{-1}\Vb\to\Vt$ is the map defining the triangle 
$T^{-1}\Vb\to\Vt\to V_{i,1}\to \Vb$ for some $i\in\{1,\dots,r\}$. 
Then, one obtains the triangle 
\begin{equation}\label{compare1}
\tau^{-1}(\Vb')\to V_{i,1}\to T\Vt\to T\tau^{-1}(\Vb')  
\end{equation} 
as $Z'\to Y'\to X'\to T(Z')$. 
On the other hand, as the transpose of the triangle (\ref{coneXX*}), 
one obtains 
\begin{equation}\label{compare2}
 t(\Vb)\to V_{i,1}\to T\Vt\to Tt(\Vb).  
\end{equation}
By comparing these two triangles (\ref{compare1}) (\ref{compare2}), 
one can see that $t(\Vb)$ is isomorphic to $\tau^{-1}(\Vb')$.

 \vspace*{0.3cm}

\noindent
{\bf The exceptional collection $\E$ is full :}\quad 
Statements (i), (ii), (iv) and (iv') for the collection $\{\E\}$ 
together with Corollary \ref{cor:mainI} leads 
that $\E$ forms a {\em full} exceptional collection in ${\T_W}$ 
as follows. 

For the exceptional collection $\E=(E_1,\dots,E_{\nu_W})$, 
let $\la\E\ra:=\la E_1,\dots,E_{\nu_W}\ra$ 
denote the smallest full triangulated category 
including objects $\{\E\}$. 
First, in the AR triangle (\ref{AR1}), one has 
${V}_0,\Vt\in\la\E\ra$, which leads that 
$T^n(\tau({V}_0))\in\la\E\ra$ for any $n\in\Z$ 
since $T\tau {V}_0$ is isomorphic to 
the cone of $\Vt\to {V}_0$. 
Next, 
in the AR triangle (\ref{AR3}), one has 
${V}_0,\Vt\in \la\E\ra$ and 
$V_2\in\la\E\ra$ since $TV_2\in\la\E\ra$. 
This implies that $T^n(\tau\Vt)\in\la\E\ra$ for any $n\in\Z$ 
since $T(\tau\Vt)$ is isomorphic to the cone of 
$V_2\oplus {V}_0\to\Vt$. 
In a similar way, by the AR triangle (\ref{AR3'}), one has 
$T^n(\tau\Vb)\in\la\E\ra$. 
Then, by the triangle (\ref{coneXX*}), 
one sees that $T^n(\tau {V}_{i,1})\in\la\E\ra$ for each $i=1,\dots,r$. 
Next, in the AR-triangle (\ref{AR1}) 
\begin{equation*}
  \tau {V}_{i,j} \to \tau {V}_{i,j-1}\oplus {V}_{i,j+1} \to {V}_{i,j} \to 
 T\tau {V}_{i,j}
\end{equation*}
with $j=2$, we already know that 
$\tau {V}_{i,1},{V}_{i,3},{V}_{i,2}\in\la\E\ra$. 
Thus, one obtains $\tau {V}_{i,2}\in\la\E\ra$. 
The AR-triangle (\ref{AR1}) with $j=3$ then implies 
that $\tau {V}_{i,3}\in\la\E\ra$, 
and repeating this argument leads that $\tau {V}_{i,j}\in\la\E\ra$ 
for any $j=1,\dots,\alpha_i$ with any $i=1,\dots,r$. 

Thus, what we obtained is $\tau {V}\in\la\E\ra$ for any $V\in\{\E\}$. 
Therefore, one can repeat the same procedure for $\tau {V}$, $V\in\{\E\}$, 
and then obtain that $\tau^n {V}\in\la\E\ra$ 
for any $n\in\Z_{\ge 0}$. 
Recall that, in this triangulated category ${\T_W}$, 
the identity $T^2=\tau^h$ holds, which implies that 
$\tau^n {V}\in\la\E\ra$, $V\in\{\E\}$, for any $n\in\Z$. 

Now, 
since ${V}_0$ corresponds to $R/\fm$ in ${\T_W}\simeq \Dbsing(R)$ 
up to grading shift, 
Theorem \ref{thm:mainI} 
and in particular Corollary \ref{cor:mainI} can be applied 
to the exceptional collection $\la\E\ra$, 
which implies that $\la\E\ra$ is full. 
\qed

 \subsection{Proof of Theorem \ref{thm:main}}\hfill 
\label{ssec:proof-main}

In this subsection, we give a proof of the main theorem 
(Theorem \ref{thm:main}) using the structure theorem 
(Theorem \ref{thm:key}) shown in the previous subsection. 

Before discussing each case $\vec{\Delta}_W^T$ or 
$\vec{\Delta}_W'$ separately, 
let us first prepare the following two lemmas 
(Lemma \ref{lem:comp1} and Lemma \ref{lem:comp2}) 
which follow from the AR-triangles in the structure theorem 
(Theorem \ref{thm:key}). 
\begin{lem}\label{lem:comp1}
For a fixed $i\in\{1,\dots,r\}$ and 
$j\in\{3,4,\dots,\alpha_i\}$, 
the composition of a nonzero element in $\Hom_{\T_W}(V_{i,j},V_{i,j-1})$ 
and a nonzero element in $\Hom_{\T_W}(V_{i,j-1}, V)$ 
gives a nonzero element in $\Hom_{\T_W}(V_{i,j},V)$ 
for any $V\in\{V_{i,j-2},\dots,V_{i,1}\}\coprod\{\Vb,TV_1,TV_0\}$. 
\end{lem}
\begin{pf}
Applying $\Hom_{\T_W}(\,\cdot\, ,\tau_{AR}V)$ to the AR-triangle (\ref{AR1}) 
gives arise to the following short exact sequence 
\begin{equation}\label{comp-AR1}
 \begin{split}
 0\to\Hom_{\T_W}(V_{i,j},\tau_{AR}V)\to
 &\Hom_{\T_W}(V_{i,j+1},\tau_{AR}V) \\
 \oplus&\Hom_{\T_W}(\tau_{AR}V_{i,j-1},\tau_{AR}V)\to 
 \Hom_{\T_W}(\tau_{AR}V_{i,j},\tau_{AR}V)\to 0. 
 \end{split}
\end{equation}
Consider the case $j=\alpha_j$. Then, since the term 
$\Hom_{\T_W}(V_{i,j+1},\tau_{AR}V)$ is absent, 
the map $\Hom_{\T_W}(\tau_{AR}V_{i,\alpha_i-1},\tau_{AR}V)\to 
\Hom_{\T_W}(\tau_{AR}V_{i,\alpha_i},\tau_{AR}V)$ is surjective and 
in particular bijective since 
$\hom_{\T_W}(\tau_{AR}V_{i,\alpha_i-1},\tau_{AR}V)=
\hom_{\T_W}(\tau_{AR}V_{i,\alpha_i},\tau_{AR}V)=1$. 
This gives the statement of this lemma 
for the case $j=\alpha_i$ together with 
that $\Hom_{\T_W}(V_{i,\alpha_i},\tau_{AR}V)=0$. 
Next, consider the short exact sequence (\ref{comp-AR1}) for the case 
$j=\alpha_i-1$. 
As we saw just now, since $\Hom_{\T_W}(V_{i,\alpha_i},\tau_{AR}V)=0$, 
the map $\Hom_{\T_W}(\tau_{AR}V_{i,\alpha_i-2},\tau_{AR}V)\to 
\Hom_{\T_W}(\tau_{AR}V_{i,\alpha_i-1},\tau_{AR}V)$ is bijective, 
which gives the statement of this lemma 
for the case $j=\alpha_i-1$, and 
also that $\Hom_{\T_W}(V_{i,\alpha_i-1},\tau_{AR}V)=0$. 
Repeating this procedure gives the statement of this lemma 
for all $j\in\{3,\dots,\alpha_i\}$. 
\qed\end{pf}
\begin{lem}\label{lem:comp2}
For any $V\in\coprod_{i=1}^r\{V_{i,2},\dots,V_{i,\alpha_i}\}\coprod
\{\Vb,\Vb'\}$, 
the composition of a nonzero element in $\Hom_{\T_W}(V,TV_1)$ 
and a nonzero element in $\Hom_{\T_W}(TV_1,TV_0)$ 
gives a nonzero element in $\Hom_{\T_W}(V,TV_0)$. 
In particular, this induces an isomorphism 
$\Hom_{\T_W}(V,TV_1)\simeq\Hom_{\T_W}(V,TV_0)$. 
\end{lem}
\begin{pf}
We may apply $\Hom_{\T_W}(V,\,\cdot\, )$ to the AR-triangle 
(\ref{AR2}) and then obtain the short exact sequence 
\begin{equation*}
 0\to\Hom_{\T_W}(V,\tau_{AR}TV_0)\to 
     \Hom_{\T_W}(V,TV_1)\to
     \Hom_{\T_W}(V,TV_0)\to 0. 
\end{equation*}
Since $\hom_{\T_W}(V,TV_1)=\hom_{\T_W}(V,TV_0)$, 
the map $\Hom_{\T_W}(V,TV_1)\to\Hom_{\T_W}(V,TV_0)$ 
is bijective, which implies Lemma \ref{lem:comp2}. 
\qed\end{pf}

Now, we show Theorem \ref{thm:main} 
for the quiver $\vec{\Delta}_W^T$. 
Consider the collection 
$$
\{\E^T\}:=
\coprod_{i=1}^r\{V_{i,1},\dots,V_{i,\alpha_i}\}\coprod \{TV_0,T\Vt,\Vb\}.
$$ 
By Theorem \ref{thm:key}, 
this forms a strongly exceptional collection 
by giving an appropriate ordering. 

By Lemma \ref{lem:comp1} and Lemma \ref{lem:comp2}, 
in order to obtain the composition law of 
morphisms between $\{\E^T\}$, 
the remaining thing is only to check 
the relations corresponding to $\rho(v_{i,2},v_1)$. 
These relations are obtained by applying $\Hom(V_{i,2},\,\cdot\, )$ to 
the triangle (\ref{coneXX*}); 
the resulting exact sequence is 
\begin{equation*}
 0\to\C\to
 \Hom_{\T_W}(V_{i,2}, \Vb)\overset{u_{1,i}e_1+u_{2,i}e_2}{\longrightarrow}
 \Hom_{\T_W}(V_{i,2}, T(\Vt))\to\C\to 0. 
\end{equation*}
Here, 
$(u_{1,i}e_1+u_{2,i}e_2):\Hom_{\T_W}(V_{i,2}, \Vb)\to\Hom_{\T_W}(V_{i,2}, T(\Vt))$ 
is a zero map: 
\begin{equation}\label{relation1-proof}
 (u_{1,i}e_1+u_{2,i}e_2)\circ f_i=0,\qquad i=1,\dots,r, 
\end{equation}
since $\hom_{\T_W}(V_{i,2}, \Vb)=1$ and $\hom_{\T_W}(V_{i,2}, T(\Vt))=1$. 
This implies the relation (\ref{relation1}): 
\begin{equation*}
 (u_{1,i}\rho_1(\vb,\vt) + u_{2,i}\rho_2(\vb,\vt))\circ\rho(v_{i,2},\vb)), 
\end{equation*}
where we identify $\rho_1(\vb,\vt)$ and $\rho_2(\vb,\vt)$ 
with $e_1$ and $e_2$, respectively. 
Thus, Theorem \ref{thm:main} has been completed for 
the quiver $\vec{\Delta}_W^T$.  

\vspace*{0.3cm}

Next, we show Theorem \ref{thm:main} for the quiver 
$\vec{\Delta}_W'$. 
Consider the collection 
\begin{equation*}
 \{\E'\}:=
 \coprod_{i=1}^r
\{V_{i,1},\dots,V_{i,\alpha_i}\}\coprod \{TV_0,T\Vt,\Vb'\}. 
\end{equation*}
Recall that $\Vb'\in{\T_W}$ is defined by the triangle 
\begin{equation}\label{Vb-t}
 \Vb'\to\oplus_{i=1}^r {V}_{i,2}\overset{\oplus_i f_i}{\longrightarrow}
 \Vb\to T(\Vb') 
\end{equation}
and we can set $\Vb'=\tau(t(\Vb))$ due to Theorem \ref{thm:key} (iii). 
\begin{lem}\label{lem:Vb-sym}
The transpose of the triangle (\ref{Vb-t}) 
is isomorphic to the triangle (\ref{Vb-t}) itself. 
Equivalently, 
the morphism $\Vb'\to\oplus_{i=1}^r {V}_{i,2}$ 
in the triangle (\ref{Vb-t}) is given by $\oplus_i \tau^{-1} t(f_i)$. 
\end{lem}
\begin{pf}
By applying $\tau^{-1} t$ to the triangle (\ref{Vb-t}), 
one obtains the triangle 
\begin{equation*}
 \Vb'\overset{\oplus_i \tau^{-1} t(f_i)}{\longrightarrow}
 \oplus_{i=1}^r {V}_{i,2}\to\Vb\to T(\Vb'),  
\end{equation*}
which shows the statement of this lemma. 
\qed\end{pf}

Now, 
by applying $\Hom_{\T_W}(\,\cdot\, ,V_{i,j})$ to the triangle (\ref{Vb-t}) 
for any $i=1,\dots, r$ and $j=1,\dots,\alpha_i$, 
one obtains that 
$\hom_{\T_W}(T^{-n}(\Vb'),V_{i,j})$ is equal to one only for 
$n=0$ and $j=2$, and is zero otherwise. 
Also, applying $\Hom_{\T_W}(V_{i,j},\,\cdot\, )$ and 
$\Hom_{\T_W}(V_k,\, \cdot\, )$, $k=0,1$, 
to the triangle (\ref{Vb-t}) leads that 
$\hom_{\T_W}(V_{i,j},\Vb')=\hom_{\T_W}(V_k,\Vb')=0$ 
for any $i$, $j$ and $k$.

Let us apply $\Hom_{\T_W}(\,\cdot\, , T\Vt)$ 
to the triangle (\ref{Vb-t}). 
Then, we obtain the following exact sequence 
\begin{equation}\label{relation2-proof}
 0\to \Hom_{\T_W}(\Vb,T\Vt)
 \overset{\oplus_i f_i}{\longrightarrow}
 \Hom_{\T_W}(\oplus_{i=1}^r  V_{i,2}, T\Vt)
 \overset{\oplus_i \tau t(f_i)}{\longrightarrow}
 \Hom_{\T_W}(\Vb', T\Vt)\to 0
\end{equation}
which gives relations of the path algebra $\PA'$. 
Here, in order to obtain the exact sequence above, 
we used the fact $\Hom_{\T_W}(T^{-n}\Vb',T\Vt)=0$ for any $n\ne 0$. 
The $n=-1$ case, 
$\Hom_{\T_W}(T\Vb',T\Vt)=0$, is nontrivial, 
which is equivalent to 
\begin{equation}\label{fact}
 \Hom_{\T_W}(T\tau(\Vt),\Vb)=0 
\end{equation}
by the transpose $t$. 
This equality (\ref{fact}) can be shown in the following two steps. 
First, by applying the functor $\Hom_{\T_W}(T\tau(\Vt),\,\cdot\,)$ 
to the triangle 
\begin{equation*}
 \oplus_{i=1}^rV_{i,2}\to AR(\Vb)\to 
 T\tau\Vt^{\oplus 2}\to T(\oplus_{i=1}^rV_{i,2}), 
\end{equation*}
one obtains $\Hom_{\T_W}(T\tau(\Vt),AR(\Vb))=0$ because 
the map $\Hom_{\T_W}(T\tau(\Vt),(T\tau\Vt)^{\oplus 2})\to
\Hom_{\T_W}(T\tau(\Vt), T(\oplus_{i=1}^rV_{i,2}))$ is injective 
and 
\begin{equation*}
\Hom_{\T_W}(T\tau(\Vt), (\oplus_{i=1}^rV_{i,2}))\simeq
\Hom_{\T_W}(T(\Vt), \tau^{-1}(\oplus_{i=1}^rV_{i,2}))\simeq
\Hom_{\T_W}((\oplus_{i=1}^rV_{i,2}),\Vt)=0, 
\end{equation*}
where the second isomorphism follows from the transpose.  
Next, by applying the functor $\Hom_{\T_W}(T\tau(\Vt),\,\cdot\,)$ 
to the AR-triangle of $\Vb$ (eq.(\ref{AR3'})), one obtains eq.(\ref{fact}). 

Thus, one obtains the exact sequence (\ref{relation2-proof}). 
By Lemma \ref{lem:Vb-sym}, one has the morphism 
$\oplus_i\tau^{-1} t(f_i): \Hom_{\T_W}(\oplus_{i=1}^r  V_{i,2}, T\Vt)
\to\Hom_{\T_W}(\Vb', T\Vt)$. We identify $\tau^{-1} t(f_i)$ with 
$\rho(\vb,v_{i,2})$ for $i=1,\dots,r$. 
On the other hand, 
the exact sequence (\ref{relation2-proof}) implies 
the following two relations 
\begin{equation*}
  \sum_{i=1}^r( e_1\circ f_i)\circ \tau^{-1}t(f_i) =0,\qquad 
  \sum_{i=1}^r( e_2\circ f_i)\circ \tau^{-1}t(f_i) =0. 
\end{equation*}
We can assume that $|u_{1,i}|^2+|u_{2,i}|^2=1$ for each $i=1,\dots,i$ 
since all the triangles in $\T_W$ depend on the ratio $u_{1,i}:u_{2,i}$ only. 
Then, using the relation (\ref{relation1-proof}), one obtains 
\begin{equation*}
 e_1\circ f_i =((1-u_{1,i}^*u_{1,i}) e_1 -u_{1,i}^*u_{2,i}e_2)\circ f_i
 = u_{2,i} (u_{2,i}^* e_1-u_{1,i}^* e_2)\circ f_i, 
\end{equation*}
where $u_{1,i}^*$ and $u_{2,i}^*$ are the complex conjugates of 
$u_{1,i}$ and $u_{2,i}$, respectively. 
We can set $g_i=(u_{2,i}^* e_1-u_{1,i}^* e_2)\circ f_i$ and 
identify $g_i$ with $\rho(v_{i,2},v_1)$. 
Thus, one of the relations (\ref{relation2}), 
$\sum_i  u_{2,i}\cdot(\rho(v_{i,2},v_1)\circ\rho(\vb,v_{i,2}))$, 
is obtained. 
In a similar way, 
the other relation of eq.(\ref{relation2}), 
$- \sum_i u_{1,i}\cdot(\rho(v_{i,2},v_1)\circ\rho(\vb,v_{i,2}))$, 
is obtained from $\sum_{i=1}^r( e_2\circ f_i)\circ \tau^{-1}t(f_i) =0$.

Now, together with Lemma \ref{lem:comp2}, 
we can conclude that $\{\E'\}$ forms a strongly exceptional collection.  
Furthermore, together with Lemma \ref{lem:comp1}, 
we checked all the composition law. 
Theorem \ref{thm:main} has been completed for 
the quiver $\vec{\Delta}_W'$.  
\qed
\begin{rem}\label{rem:AR'}
For the strongly exceptional collection 
$\E'$ associated to the quiver $\vec{\Delta}_W'$, 
the AR-triangles can be described in terms of 
the indecomposable objects $\{\E'\}$ and their grading shifts. 
In particular, the transpose 
of the AR-triangles (\ref{AR3}) and (\ref{AR3'}) yield 
\begin{equation*}
 \begin{split}
  & AR(\Vt)\simeq \tau\Vt\oplus C((\Vb')^{\oplus 2}
 \to\oplus_{i=1}^r V_{i,2}),\\
  & AR(\Vb')\simeq T^{-1}\tau C(\oplus_{i=1}^r V_{i,2}
 \to (T\Vt)^{\oplus 2}), 
 \end{split}
\end{equation*}
where the morphisms in the mapping cones are also determined by the transpose. 
\end{rem}

 \section{Graded matrix factorizations associated to the vertices of 
the quiver $\vec{\Delta}_W$}
\label{sec:list}

In this section, we list up the graded matrix factorizations 
which form the exceptional collection $\E$ and are 
in the reduced form, i.e., 
representatives having minimum rank in the isomorphism classes.

 \subsection{The grading matrices}\hfill
\label{ssec:MF-grading}

We list up the grading matrices for the reduced indecomposable 
graded matrix factorizations of the corresponding quivers $\vec{\Delta}_W$.

For the grading matrix 
$S:=\diag\{s_1,\dots,s_r;\sb_1,\dots,\sb_r\}$, 
recall that $s_i\in 2\Z/h$ and $\sb_i\in 2\Z/h-1$ for $i=1,\dots,r$. 
For some $\phi\in 2\Z/h-1/2$, rewrite 
\begin{equation*}
 \diag\{s_1,\dots,s_r;\sb_1,\dots,\sb_r\} 
 =\diag\{s'_1,\dots,s'_r;\sb'_1,\dots,\sb'_r\}+\phi\cdot\1
\end{equation*}
where $s'_i=s_i-\phi$ and $\sb'_i=\sb_i-\phi$ for $i=1,\dots,r$. 
Except for $\Vb$, we take this $\phi$ so that 
$S-\phi\cdot\1$ is traceless. 
We describe the data of the grading matrix as 
\begin{equation*}
 [h s'_1,\dots,h s'_r;h \sb'_1,\dots,h \sb'_r]_{h\phi} .
\end{equation*}
In particular, we use the following simplified notation 
\begin{equation*}
 (q_1,\dots, q_\nu;\qb_1,\dots,\qb_\nu)_{h\phi} 
 :=[q_1,-q_1,\dots, q_\nu,-q_\nu;
 \qb_1,-\qb_1,\dots,\qb_r,-\qb_r]_{h\phi}\ ,
\end{equation*}
for a grading matrix of this kind. 
Furthermore, we denote 
\begin{equation*}
 (q_1,\dots,q_\nu)_{h\phi}
 :=(q_1,\dots, q_\nu;\qb_1,\dots, \qb_\nu)_{h\phi} 
\end{equation*}
if $q_i=\qb_i$ for $i=1,\dots, \nu$. 
Similarly, if $s'_i=\sb'_i$ for any $i=1,\dots,r$, 
we denote 
\begin{equation*}
 [h s'_1,\dots,h s'_r]_{h\phi}
 :=[h s'_1,\dots,h s'_r;h \sb'_1,\dots,h \sb'_r]_{h\phi}. 
\end{equation*}
If we attach the grading matrices below to the vertices 
of the corresponding quiver $\vec{\Delta}_W$, we see 
some interesting phenomenological rules 
due to the Auslander-Reiten triangles in Theorem \ref{thm:key} 
at the level of the Grothendieck group (see eq.(\ref{rem:key})) 
as in the case that $W$ is of type ADE 
(see \cite[section 5, Table 2]{KSTade}).

\vspace*{0.2cm}

\noindent
$\bullet$\ $W=(6,14,21;42)$, $f_W=x^7+y^3+z^2$, $A_W=(2,3,7)$, 
\begin{equation*}
 \begin{split}
 & V_{1,2}: (3,7,11)_1 ,\\
 & V_{2,2}: (3,5,9,11)_1\ ,\quad V_{2,3}: (4,10)_2 ,\\
 & V_{3,2}: (3,5,7,9,11,13,27)_1 ,\quad 
   V_{3,3}: (4,6,8,10,12,26)_2 ,\quad 
   V_{3,4}: (5,7,9,11,25)_3 ,\\ 
 & \qquad V_{3,5}: (6,8,10,24)_4 ,\quad 
   V_{3,6}: (7,9,23)_5 ,\quad 
   V_{3,7}: (8,22)_6 ,\\
 & V_0: (8,22)_{-22} ,\quad \Vt: (7,9,23)_{-21} ,\quad 
\Vb: [28,14,12^2,10,8,6,4,-2^2,-4,-6,-8,-10].
 \end{split}
\end{equation*}

\noindent
$\bullet$\ $W=(4,10,15;30)$, $f_W=y^3+yx^5+z^2$, $A_W=(2,4,5)$, 
\begin{equation*}
 \begin{split}
 & V_{1,2}: (3,7)_1 ,\qquad 
 V_{2,2}: (3,5,7)_1 ,\quad V_{2,3}: (4,6)_2 ,\quad V_{2,4}: (5)_3 ,\\
 & V_{3,2}: (3,5,7,9,19)_1 ,\quad 
   V_{3,3}: (4,6,8,18)_2 ,\quad 
   V_{3,4}: (5,7,17)_3 ,\quad V_{3,5}: (6,16)_4 ,\\
 & V_0: (6,16)_{-16} ,\quad \Vt: (5,7,17)_{-15} ,\quad 
 \Vb:[20,10,8^2,6,4,-2^2,-4,-6] .
 \end{split}
\end{equation*}

\noindent
$\bullet$\ $W=(3,8,12;24)$, $f_W=x^4z+y^3+z^2$, $A_W=(3,3,4)$, 
\begin{equation*}
 \begin{split}
 & V_{1,2}: (3,5)_1 ,\quad V_{1,3}: (4)_2 ,\qquad 
 V_{2,2}: (3,5)_1 ,\quad V_{2,3}: (4)_2 ,\\
 & V_{3,2}: (3,5,7,15)_1 ,\quad 
   V_{3,3}: (4,6,14)_2 ,\quad 
   V_{3,4}: (5,13)_3 , \\
 & V_0: (5,13)_{-13} ,\quad \Vt: (4,6,14)_{-12} ,\quad 
 \Vb: [16,8,6^2,4,-2^2,-4] .
 \end{split}
\end{equation*}

\noindent
$\bullet$\ $W=(6,8,15;30)$, $f_W=x^5+xy^3+z^2$, $A_W=(2,3,8)$, 
\begin{equation*}
 \begin{split}
 & V_{1,2}: (3,7,11)_1 ,\qquad 
 V_{2,2}: (3,5,9,11)_1 ,\quad V_{2,3}: (4,10)_2 , \\
 & V_{3,2}: (3,5,7,9,11,13,15)_1 ,\quad 
   V_{3,3}: (4,6,8,10,12,14)_2 ,\\
 & V_{3,4}: (5,7,9,11,13)_3 ,\quad
   V_{3,5}: (6,8,10,12)_4 ,\quad 
   V_{3,6}: (7,9,11)_5 ,\quad
   V_{3,7}: (8,10)_6 ,\quad
   V_{3,8}: (9)_7 ,\\
 & V_0: (2,16)_{-16} ,\quad \Vt: (1,3,17)_{-15} ,\quad 
 \Vb: [16,14,12^2,10,8,6,4,-2^2,-4,-6,-8,-10] .
 \end{split}
\end{equation*}

\noindent
$\bullet$\ $W=(4,6,11;22)$, $f_W=yx^4+xy^3+z^2$, $A_W=(2,4,6)$, 
\begin{equation*}
 \begin{split}
 & V_{1,2}: (3,7)_1 ,\qquad 
 V_{2,2}: (3,5,7)_1 ,\quad V_{2,3}: (4,6)_2 ,\quad V_{2,4}: (5)_3 ,\\
 & V_{3,2}: (3,5,7,9,11)_1 ,\quad 
   V_{3,3}: (4,6,8,10)_2 ,\quad 
 V_{3,4}: (5,7,9)_3 ,\quad
   V_{3,5}: (6,8)_4 ,\quad 
   V_{3,6}: (7)_5 ,\\
 & V_0: (2,12)_{-12} ,\quad \Vt: (1,3,13)_{-11} ,\quad 
 [12,10,8^2,6,4,-2^2,-4,-6] .
 \end{split}
\end{equation*}

\noindent
$\bullet$\ $W=(3,5,9;18)$, $f_W=x^3z+xy^3+z^2$, $A_W=(3,3,5)$, 
\begin{equation*}
 \begin{split}
 & V_{1,2}: (3,5)_1 ,\quad V_{1,3}: (4)_2 ,\qquad 
 V_{2,2}: (3,5)_1 ,\quad V_{2,3}: (4)_2 ,\\
 & V_{3,2}: (3,5,7,9)_1 ,\quad 
   V_{3,3}: (4,6,8)_2 ,\quad V_{3,4}: (5,7)_3 ,\quad
   V_{3,5}: (6)_4 ,\\
 & V_0: (2,10)_{-10} ,\quad \Vt: (1,3,11)_{-9} ,\quad 
 \Vb: [10,8,6^2,4,-2^2,-4] .
 \end{split}
\end{equation*}

\noindent
$\bullet$\ $W=(4,5,10;20)$, $f_W=x^5+y^2z+z^2$, $A_W=(2,5,5)$, 
\begin{equation*}
 \begin{split}
 & V_{1,2}: (3,7)_1 ,\\
 & V_{2,2}: (3,5,7,9)_1 ,\quad V_{2,3}: (4,6,8)_2 ,\quad 
   V_{2,4}: (5,7)_3 ,\quad V_{2,5}: (6)_4 , \\
 & V_{3,2}: (3,5,7,9)_1 ,\quad  
   V_{3,3}: (4,6,8)_2 ,\quad
 V_{3,4}: (5,7)_3 ,\quad
   V_{3,5}: (6)_4 ,\\
 & V_0: (1,11)_{-11} ,\quad \Vt: (0,2,12)_{-10} ,\quad 
 \Vb: [10^2,8^2,6,4,-2^2,-4,-6] .
 \end{split}
\end{equation*}

\noindent
$\bullet$\ $W=(3,4,8;16)$, $yx^4+y^2z+z^2$, $A_W=(3,4,4)$, 
\begin{equation*}
 \begin{split}
 & V_{1,2}: (3,5)_1 ,\quad V_{1,3}: (4)_2 ,\\ 
 & V_{2,2}: (3,5,7)_1 ,\quad V_{2,3}: (4,6)_2 ,\quad V_{2,4}: (5)_3 ,
 \qquad
 V_{3,2}: (3,5,7)_1 ,\quad 
   V_{3,3}: (4,6)_2 ,\quad 
   V_{3,4}: (5)_3 ,\\
 &V_0: (1,9)_{-9} ,\quad \Vt: (0,2,10)_{-8} ,\quad 
 \Vb: [8^2,6^2,4,-2^2,-4] .
 \end{split}
\end{equation*}

\noindent
$\bullet$\ $W=(6,8,9;24)$, $f_W=x^4+y^3+xz^2$, $A_W=(2,3,9)$, 
\begin{equation*}
 \begin{split}
 & V_{1,2}: (3,7,11;1,3,5)_1 ,\qquad 
 V_{2,2}: (3,5,9,11;1,3^2,5)_1 ,\quad V_{2,3}: (4,10;2,4)_2 , \\
 & V_{3,2}: (3,5,7,9,11,13,15,17;1^2,3^2,5^2,7,9)_1 ,\quad 
   V_{3,3}: (4,6,8,10,12,14,16;0,2^2,4^2,6,8)_2 ,\\
 & V_{3,4}: (5,7,9,11,13,15;1^2,3^2,5,7)_3 ,\quad 
   V_{3,5}: (6,8,10,12,14;0,2^2,4,6)_4 ,\\
 & V_{3,6}: (7,9,11,13;1^2,3,5)_5 ,\quad 
   V_{3,7}: (8,10,12;0,2,4)_6 ,\quad
   V_{3,8}: (9,11;1,3)_7 ,\quad 
   V_{3,9}: (10;2)_8 ,\\
 & V_0: [13,-1,-5,-7;7,5,1,-13]_{-13} ,
 \quad \Vt: [14,0,-2,-4,-6^2,-8;8,6^2,4,2,0,-14]_{-12} ,\\ 
 & \Vb: 
[18,16,14,12^2,10,8,6,4,-2^2,-4,-6,-8,-10;
10,8,6^3,4^3,2^2,0^2,-2^2,-4] .
 \end{split}
\end{equation*}

\noindent
$\bullet$\ $W=(4,6,7;18)$, $f_W=x^3y+y^3+xz^2$, $A_W=(2,4,7)$, 
\begin{equation*}
 \begin{split}
 & V_{1,2}: (3,7;1,3)_1 ,\qquad
 V_{2,2}: (3,5,7;1^2,3)_1 ,\quad V_{2,3}: (4,6;0,2)_2 ,\quad 
   V_{2,4}: (5;1)_3 , \\
 & V_{3,2}: (3,5,7,9,11,13;1^2,3^2,5,7)_1 ,\quad 
   V_{3,3}: (4,6,8,10,12;0,2^2,4,6)_2 ,\\
 & V_{3,4}: (5,7,9,11;1^2,3,5)_3 ,\quad
   V_{3,5}: (6,8,10;0,2,4)_4 ,\quad 
   V_{3,6}: (7,9;1,3)_5 ,\quad
   V_{3,7}: (8;2)_6 ,\\
 & V_0: [10,0,-4,-6;6,4,0,-10]_{-10} ,\quad 
   \Vt: [11,1,-1,-3,-5^2,-7;7,5^2,3,1,-1,-11]_{-9} ,\\
 &\Vb:[14,12,10,8^2,6,4,-2^2,-4,-6;8,6,4^3,2^3,0^2,-2] .
 \end{split}
\end{equation*}

\noindent
$\bullet$\ $W=(3,5,6;15)$, $f_W=x^3z+y^3+xz^2$, $A_W=(3,3,6)$, 
\begin{equation*}
 \begin{split}
 & V_{1,2}: (3,5;0,2)_1 ,\quad V_{1,3}: (4;1)_2 , \quad
 V_{2,2}: (3,5;0,2)_1 ,\quad V_{2,3}: (4;1)_2 , \\
 & V_{3,2}: (3,5,7,9,11;0,2^2,4,6)_1 ,\quad 
   V_{3,3}: (4,6,8,10;1^2,3,5)_2 ,\\
 & V_{3,4}: (5,7,9;0,2,4)_3 ,\quad 
   V_{3,5}: (6,8;1,3)_4 ,\quad 
   V_{3,6}: (7;2)_5 , \\
 & V_0: [17/2,1/2,-7/2,-11/2;11/2,7/2,-1/2,-17/2]_{-17/2} ,\\
 &  \Vt: [19/2,3/2,-1/2,-5/2,(-9/2)^2,-13/2;
13/2,(9/2)^2,5/2,1/2,-3/2,-19/2]_{-15/2} ,\\
 & \Vb: [12,10,8,6^2,4,-2^2,-4;7,5,3^3,1^3,-1]  .
 \end{split}
\end{equation*}

\noindent
$\bullet$\ $W=(4,5,6;16)$, $f_W=x^4+y^2z+z^2x$, $A_W=(2,5,6)$, 
\begin{equation*}
 \begin{split}
 & V_{1,2}: (3,7;1,3)_1 ,\\
 & V_{2,2}: (3,5,7,9;1^2,3,5)_1 ,\quad V_{2,3}: (4,6,8;0,2,4)_2 ,
   \quad V_{2,4}: (5,7;1,3)_3 ,\quad V_{2,5}: (6;2)_4 , \\
 & V_{3,2}: (3,5,7,9,11;1^2,3^2,5)_1 ,\quad 
   V_{3,3}: (4,6,8,10;0,2^2,4)_2 ,\quad 
   V_{3,4}: (5,7,9;1^2,3)_3 ,\\
 & V_{3,5}: (6,8;0,2)_4 ,\quad 
   V_{3,6}: (7;1)_5 , \\
 & V_0: [9,-1,-3,-5;5,3,1,-9]_{-9} ,\quad 
   \Vt: [10,0,-2^2,-4^2,-6;6,4^2,2^2,0,-10]_{-8} ,\\
 & \Vb: [12,10^2,8^2,6,4,-2^2,-4,-6;6^2,4^3,2^3,0^2,-2] .
 \end{split}
\end{equation*}

\noindent
$\bullet$\ $W=(3,4,5;13)$, $f_W=x^3y+y^2z+z^2x$, $A_W=(3,4,5)$, 
\begin{equation*}
 \begin{split}
 & V_{1,2}: (3,5;0,2)_1 ,\quad V_{1,3}: (4;1)_2 ,\quad 
 V_{2,2}: (3,5,7;0,2^2,4)_1 ,\quad V_{2,3}: (4,6;1,3)_2 ,
   \quad V_{2,4}: (5;2)_3 ,\\
 & V_{3,2}: (3,5,7,9;0,2^2,4)_1 ,\quad 
   V_{3,3}: (4,6,8;1^2,3)_2 ,\quad 
   V_{3,4}: (5,7;0,2)_3 ,\quad
   V_{3,5}: (6;1)_4 , \\
 & V_0: [15/2,-1/2,-5/2,-9/2;9/2,5/2,1/2,-15/2]_{-15/2} ,\\
 &  \Vt: [17/2,1/2,(-3/2)^2,(-7/2)^2,-11/2;
 11/2,(7/2)^2,(3/2)^2,-1/2,-17/2]_{-13/2} ,\\
 & \Vb: [10,8^2,6^2,4,-2^2,-4;5^2,3^3,1^3,-1] .
 \end{split}
\end{equation*}

\noindent
$\bullet$\ $W=(3,4,4;12)$, $f_W=x^4+yz(y-z)$, $A_W=(4,4,4)$, 
\begin{equation*}
 \begin{split}
 & V_{1,2}: (3,5,7;1^2,3)_1 ,\quad V_{1,3}: (4,6;0,2)_2 ,
 \quad V_{1,4}: (5;1)_3 ,\\
 & V_{2,2}: (3,5,7;1^2,3)_1 ,\quad V_{2,3}: (4,6;0,2)_2 ,
   \quad V_{2,4}: (5;1)_3 ,\\
 & V_{3,2}: (3,5,7;1^2,3)_1 ,\quad 
   V_{3,3}: (4,6;0,2)_2 ,\quad 
   V_{3,4}: (5;1)_3 ,\\
 & V_0: [7,-1,-3,-3;3,3,1,-7]_{-7} ,\quad 
   \Vt: [8,0,-2^3,-4^2; 4^2,2^3,0,-8]_{-6} ,\\
 & \Vb: [8^3,6^2,4,-2^2,-4;4^3,2^4,0^2] .
 \end{split}
\end{equation*}

\noindent
$\bullet$\ $W=(2,6,9;18)$, $f_W=y(y-x^3)(y-\lambda x^3)+z^2$, 
$\lambda\ne 0,1$, $A_W=(2,2,2,3)$, 
\begin{equation*}
 \begin{split}
 & V_{1,2}: (3)_1 ,\quad V_{2,2}: (3)_1 ,
 \quad V_{3,2}: (3)_1 ,\quad 
 V_{4,2}: (3,5,11)_1 ,\quad V_{4,3}: (4,10)_2 , \\
 & V_0: (4,10)_{-10} ,\quad \Vt: (3,5,11)_{-9} ,\quad 
 \Vb: [12,6,4^2,-2^2] .
 \end{split}
\end{equation*}

\noindent
$\bullet$\ $W=(2,4,7;14)$, $f_W=xy(y-x^2)(y-\lambda x^2)+z^2$, 
$\lambda\ne 0,1$, $A_W=(2,2,2,4)$, 
\begin{equation*}
 \begin{split}
 & V_{1,2}: (3)_1 ,\quad V_{2,2}: (3)_1 ,
 \quad V_{3,2}: (3)_1 ,\quad 
 V_{4,2}: (3,5,7)_1 ,\quad V_{4,3}: (4,6)_2 , \quad V_{4,4}: (5)_3 , \\
 & V_0: (2,8)_{-8} ,\quad 
   \Vt: (1,3,9)_{-7} ,\quad \Vb: [8,6,4^2,-2^2] .
 \end{split}
\end{equation*}

\noindent
$\bullet$\ $W=(2,4,5;12)$, $f_W=y(y-x^2)(y-\lambda x^2)+xz^2$, 
$\lambda\ne 0,1$, $A_W=(2,2,2,5)$, 
\begin{equation*}
 \begin{split}
 & V_{1,2}: (3;1)_1 ,\quad V_{2,2}: (3;1)_1 ,
 \quad V_{3,2}: (3;1)_1 ,\\
 & V_{4,2}: (3,5,7,9;1^2,3,5)_1 ,\quad V_{4,3}: (4,6,8;0,2,4)_2 , \quad 
 V_{4,4}: (5,7;1,3)_3 ,\quad V_{4,5}: (6;2)_4 ,\\
 & V_0: [7,1,-3,-5;5,3,-1,-7]_{-7} ,\quad 
   \Vt: [8,2,0,-2,-4^2,-6;6,4^2,2,0,-2,-8]_{-6} ,\\
 & \Vb: [10,8,6,4^2,-2^2;6,4,2^3,0^2] .
 \end{split}
\end{equation*}

\noindent
$\bullet$\ $W=(2,3,6;12)$, $f_W=(y^2-x^3)(y^2-\lambda x^3)+z^2$, 
$\lambda\ne 0,1$, $A_W=(2,2,3,3)$, 
\begin{equation*}
 \begin{split}
 & V_{1,2}: (3)_1 ,\quad V_{2,2}: (3)_1 ,\quad V_{3,2}: (3,5)_1 ,\quad 
 V_{3,3}: (4)_2 ,\quad 
V_{4,2}: (3,5)_1 ,\quad V_{4,3}: (4)_2 ,\\
 & V_0: (1,7)_{-7} ,\quad 
 \Vt: (0,2,8)_{-6} ,\quad \Vb: [6^2,4^2,-2^2] .
 \end{split}
\end{equation*}

\noindent
$\bullet$\ $W=(2,3,4;10)$, $f_W=x(z-x^2)(z-\lambda x^2)+y^2z$, 
$\lambda\ne 0,1$, $A_W=(2,2,3,4)$, 
\begin{equation*}
 \begin{split}
 & V_{1,2}: (3;1)_1 ,\quad V_{2,2}: (3;1)_1 , 
 \quad V_{3,2}: (3,5;1,3)_1 ,\quad V_{3,3}: (4;2)_2 ,\\
 & V_{4,2}: (3,5,7;1^2,3)_1 ,\quad V_{4,3}:(4,6;0,2)_2 ,\quad 
 V_{4,4}: (5;1)_3 ,\\
 & V_0: [6,0,-2,-4;4,2,0,-6]_{-6} ,\quad 
   \Vt: [7,1,-1^2,-3^2,-5;5,3^2,1^2,-1,-7]_{-5} ,\\
 & \Vb: [8,6^2,4^2,-2^2;4^2,2^3,0^2] .
 \end{split}
\end{equation*}

\noindent
$\bullet$\ $W=(2,3,3;9)$, $f_W=x^3y+z(z-y)(z-\lambda y)$, 
$\lambda\ne 0,1$, $A_W=(2,3,3,3)$, 
\begin{equation*}
 \begin{split}
 & V_{1,2}: (3;0)_1 ,\quad V_{2,2}: (3,5;0,2)_1 ,
 \quad V_{2,3}: (4;1)_2 ,\\
 & V_{3,2}: (3,5;0,2)_1 ,\quad V_{3,3}: (4;1)_2 ,
 \quad V_{4,2}:(3,5;0,2)_1 ,\quad V_{4,3}: (4;1)_2 ,\\
 & V_0: [11/2,-1/2,-5/2,-5/2;5/2,5/2,1/2,-11/2]_{-11/2} ,\\
 & \Vt: [13/2,1/2,(-3/2)^3,(-7/2)^2;(7/2)^2,(3/2)^3,-1/2,-13/2]_{-9/2} ,\\
 & \Vb: [6^3,4^2,-2^2;3^3,1^4] .
 \end{split}
\end{equation*}

\noindent
$\bullet$\ $W=(2,2,5;10)$, 
$f_W=xy(x-y)(y-\lambda_1 x)(y-\lambda_2 x)+z^2$, 
$\lambda_1\ne 0,1$, $\lambda_1\ne\lambda_2$, $A_W=(2,2,2,2,2)$, 
\begin{equation*}
 \begin{split}
 & V_{1,2}: (3)_1 ,\quad V_{2,2}: (3)_1 ,
 \quad V_{3,2}: (3)_1 , \quad
 V_{4,2}: (3)_1 ,\quad V_{5,2}: (3)_1 ,\\
 & V_0: (0,6)_{-6} ,\quad 
   \Vt: (1^2,7)_{-5} ,\quad \Vb: [6,4^3,-2^2] .
 \end{split}
\end{equation*}

\noindent
$\bullet$\ $W=(2,2,3;8)$, 
$f_W=y(y-x)(y-\lambda_1 x)(y-\lambda_2 x)+xz^2$, 
$\lambda_1\ne 0,1$, $\lambda_1\ne\lambda_2$, $A_W=(2,2,2,2,3)$, 
\begin{equation*}
 \begin{split}
 & V_{1,2}: (3;1)_1 ,\quad V_{2,2}: (3;1)_1 , 
 \quad V_{3,2}: (3;1)_1 ,\quad  
V_{4,2}: (3;1)_1 ,\quad V_{5,2}: (3,5;1^2)_1 ,\quad V_{5,3}: (4;0)_2 , \\
 & V_0: [5,-1^2,-3;3,1^2,-5]_{-5} ,\quad 
   \Vt: [6,0^2,-2^3,-4;4,2^3,0^2,-6]_{-4} ,\\
 & \Vb: [6^2,4^3,-2^2;4,2^4,0^2] .
 \end{split}
\end{equation*}

 \subsection{Matrix factorizations $V_0, V_{\alpha_i}, \Vb$} \hfill 
\label{ssec:MF-list}\hfill

We present the explicit form of a reduced graded matrix factorization 
for $V_0, V_{\alpha_i}$, and $\Vb$. 
The remaining reduced graded matrix factorizations 
are obtained by considering the AR-triangles 
(\ref{AR1}) and (\ref{AR2}). 

Since the grading matrices are already listed up in the 
previous subsection, 
we present only the part $Q=\left(\bps \0 & q_0\\ q_1 & \0 \eps\right)$ 
of each reduced graded matrix factorization $(Q,S)$.

\noindent
$\bullet$\ $W=(6,14,21;42)$, $A_W=(2,3,7)$. 
For the polynomial $f_W=x^7+y^3+z^2$, 
\begin{equation*}
 V_0 :
 {q_0}
 = {q_1} =\bp z & y^2 & x^6 & 0 \\
      y & -z & 0 & x^6\\ 
      x & 0 & -z & -y^2 \\
      0 & x & -y & z \ep\ ,\quad 
 V_{1,2} :
 {q_0} = {q_1} =
{\small 
\bp 
 z  & 0  & 0  & y^2   & x^5  & x^3y \\
 0  & z  & 0  & -x^2y & y^2  & -x^5 \\
 0  & 0  & z  & x^4   & -x^2y& -y^2 \\ 
 y  & 0  & x^3& -z    & 0    & 0    \\
 x^2& y  & 0  & 0     & -z   & 0    \\
 0  &-x^2& -y & 0     & 0    & -z    
\ep\ ,
}
\end{equation*}
\begin{equation*}
 V_{2,3} :
 {q_0} = {q_1} =
\bp 
  z  & x^4& -y^2& 0 \\ 
  x^3& -z & 0   & y^2 \\
  -y & 0  & -z  & x^4 \\
  0  & y  & x^3 & z 
\ep\ ,\qquad 
 V_{3,7} : q_0=q_1=q_0(V_0)=q_1(V_0)\ , 
\end{equation*}
\begin{equation*}
  \Vb: {q_0}={q_1}= 
{\footnotesize 
 \left(
\begin{array}{cccccccccccccc}
 z & -y^2 & 0 & 0 & 0 & 0 & -x^3y & 0 & -x^6 & 0 & 0 & 0 & 0 & 0 \\
 -y & -z & 0 & 0 & 0 & 0 & 0 & 0 & 0 & 0 & -x^5 & 0 & x^3y & 0 \\
 0 & xy & z & 0 & 0 & 0 & x^4 & 0 & -y^2 & 0 & 0 & 0 & 0 & 0 \\
 0 & -xy & 0 & z & 0 & 0 & 0 & 0 & 0 & y^2 & 0 & x^5 & 0 & x^2y \\
 x^2 & 0 & 0 & 0 & -z & 0 & 0 & 0 & 0 & 0 & -y^2 & 0 & -x^5 & 0 \\
 0 & x^2 & 0 & 0 & 0 & z & 0 & 0 & 0 & -x^2y & 0 & y^2 & 0 & -x^5 \\
 0 & 0 & x^2 & 0 & 0 & 0 & -z & 0 & 0 & 0 & x^2y & 0 & -y^2 & 0  \\
 0 & 0 & 0 & 0 & x^3 & 0 & 0 & z & 0 & x^4 & 0 & -x^2y & 0 & -y^2 \\
 -x & 0 & -y & 0 & 0 & 0 & 0 & 0 & -z & 0 & 0 & 0 & 0 & 0 \\
 -x & 0 & 0 & y & 0 & 0 & 0 & x^3 & 0 & -z & 0 & 0 & x^4 & 0 \\
 0 & -x^2 & 0 & 0 & -y & 0 & 0 & 0 & 0 & 0 & z & 0 & 0 & 0 \\
 0 & 0 & 0 & x^2 & 0 & y & 0 & 0 & 0 & 0 & 0 & -z & 0 & 0 \\
 0 & 0 & 0 & 0 & -x^2 & 0 & -y & 0 & -x^3 & 0 & 0 & 0 & z & 0 \\
 0 & 0 & 0 & 0 & 0 & -x^2 & 0 & -y & 0 & 0 & -x^3 & 0 & 0 & -z
\end{array}
\right)} \ .
\end{equation*}

\noindent
$\bullet$\ $W=(4,10,15;30)$, $A_W=(2,4,5)$. 
For the polynomial $f_W=y^3+yx^5+z^2$, 
\begin{equation*}
 V_0: {q_0} 
 = {q_1} 
 =
\bp z & y^2 & x^4y & 0 \\ 
    y & -z & 0 & x^4y \\
    x & 0 & -z & -y^2 \\ 
    0 & x & -y & z \ep\ ,\quad 
 V_{1,2} :
{q_0} = {q_1} =
\bp z & 0 & y^2 & x^3y \\
    0 & -z & x^2y & -y^2 \\ 
    y & x^3 & z & 0 \\
    x^2 & -y & 0 & z 
\ep\ ,
\end{equation*}
\begin{equation*}
 V_{2,4} :
 {q_0} = {q_1} =
\bp z & x^5+y^2 \\
    y & -z \ep\ ,\qquad 
 V_{3,5} : q_0=q_1= q_0(V_0)=q_1(V_0)\ ,  
\end{equation*}
\begin{equation*}
 \Vb: 
 {q_0}
 = {q_1} =
{\footnotesize 
\bp -z & y^2 & 0 & 0 & 0 & 0 & x^4y & 0 & 0 & -x^2y^2 \\
    y & z & 0 & 0 & 0 & x^2y & 0 & 0 & 0 & 0 \\ 
    0 & -xy & -z & 0 & 0 & 0 & y^2 & 0 & 0 & x^3y \\
    0 & -xy & 0 & -z & x^4 & 0 & 0 & -y^2 & 0 & 0 \\
    0 & 0 & 0 & xy & z & 0 & 0 & 0 & -y^2 & 0 \\
    0 & x^3 & 0 & 0 & 0 & -z & -x^2y & 0 & 0 & 0 \\
    x & 0 & y & 0 & 0 & 0 & z & 0 & 0 & 0 \\
    -x & 0 & 0 & -y & 0 & -x^3 & 0 & z & -x^4 & 0 \\ 
    0 & -x^2 & 0 & 0 & -y & 0 & 0 & -xy & -z & 0 \\
    0 & 0 & x^2 & 0 & 0 & y & 0 & 0 & 0 & z
\ep
}\ .
\end{equation*}

\noindent
$\bullet$\ $W=(3,8,12;24)$, $A_W=(3,3,4)$. 
For the polynomial $f_W=x^4z+y^3+z^2$, 
\begin{equation*}
 V_0 : 
{q_0}={q_1}=
\bp z & y^2 & x^3z & 0 \\ 
    y & -z & 0 & x^3z \\ 
    x & 0 & -z & -y^2 \\
    0 & x & -y & z \ep\ , 
\end{equation*}
\begin{equation*}
 V_{1,3} :
 {q_0}=
 \bp z+x^4 & y^2 \\ y & -z\ep\ ,\qquad 
 {q_1}=
 \bp z & y^2 \\ y & -(z+x^4)\ep\ ,  
\end{equation*}
\begin{equation*}
 V_{2,3} :
 {q_0}=
 \bp z & y^2 \\ y & -(z+x^4)\ep\ ,\qquad 
 {q_1}=
 \bp z+x^4 & y^2 \\ y & -z \ep\ , 
\end{equation*}
\begin{equation*}
 V_{3,4} : q_0=q_1=q_0(V_0)=q_1(V_0)\ , 
\end{equation*}
\begin{equation*}
 \Vb : 
{q_0}=
{\footnotesize 
 \bp 
 -z & y^2 & 0 & 0 & 0 & 0 & 0 & 0 \\
 y & z+x^4 & 0 & 0 & 0 & 0 & 0 & 0 \\
 0 & xy & -z & 0 & 0 & -y^2 & 0 & 0 \\
 0 & -xy & 0 & z+x^4 & 0 & 0 & y^2 & x^3y \\ 
 -x^2 & 0 & 0 & 0 & z+x^4 & 0 & 0 & -y^2 \\
 x & 0 & -y & 0 & 0 & z+x^4 & 0 & 0 \\
 -x & 0 & 0 & y & x^3 & 0 & -z & 0 \\
 0 & -x^2 & 0 & 0 & -y & 0 & 0 & -(z+x^4)
 \ep
}\ , 
\end{equation*}
\begin{equation*}
{q_1}=
{\footnotesize 
 \bp 
 -(z+x^4) & y^2 & 0 & 0 & 0 & 0 & 0 & 0 \\
 y & z & 0 & 0 & 0 & 0 & 0 & 0 \\
 0 & xy & -(z+x^4) & 0 & 0 & -y^2 & 0 & 0 \\
 0 & xy & 0 & z & 0 & 0 & y^2 & x^3y \\ 
 -x^2 & 0 & 0 & 0 & z & 0 & 0 & -y^2 \\
 x & 0 & -y & 0 & 0 & z & 0 & 0 \\
 x & 0 & 0 & y & x^3 & 0 & -(z+x^4) & 0 \\
 0 & -x^2 & 0 & 0 & -y & 0 & 0 & -z 
 \ep
}\ . 
\end{equation*}

\noindent
$\bullet$\ $W=(6,8,15;30)$, $A_W=(2,3,8)$. 
For the polynomial $f_W=x^5+xy^3+z^2$, 
\begin{equation*}
 V_0 : 
{q_0}={q_1}=
\bp z & xy^2 & x^4 & 0 \\ 
    y & -z & 0 & x^4 \\ 
    x & 0 & -z & -xy^2 \\
    0 & x & -y & z \ep\ , \ \ 
 V_{1,2} : 
{q_0}={q_1}=
{\small 
\bp z & 0 & 0 & xy^2 & x^4 & -x^2y \\ 
    0 & z & 0 & -x^2y & xy^2 & x^4 \\
    0 & 0 & z & x^3 & -x^2y & xy^2 \\
    y & 0 & x^2 & -z & 0 & 0 \\ 
    x & y & 0 & 0 & -z & 0 \\
    0 & x & y & 0 & 0 & -z \ep
}\ , 
\end{equation*}
\begin{equation*}
 V_{2,3} :
{q_0}={q_1}=
\bp z & x^3 & xy^2 & 0 \\ 
    x^2 & -z & 0 & -xy^2 \\ 
    y & 0 & -z & x^3 \\
    0 & -y & x & z \ep\ ,\quad 
 V_{3,8} :
{q_0}={q_1}=
\bp z & x^4+y^3 \\ 
    x & -z \ep\ ,
\end{equation*}
\begin{equation*}
 \Vb\ : 
 {q_0}={q_1}=
{\footnotesize 
 \left(
 \begin{array}{cccccccccccccc}
  z & -y^2 & 0 & 0 & 0 & 0 & -x^2y & 0 & -x^4 & 0 & 0 & 0 & 0 & 0 \\
  -xy & -z & 0 & 0 & 0 & 0 & 0 & 0 & 0 & 0 & -x^4 & 0 & x^2y & 0 \\
  0 & -xy & z & 0 & 0& 0 & -x^3 & 0 & xy^2 & 0 & 0 & 0 & 0 & 0 \\
  0 & 0 & 0 & z & y^2 & 0 & 0 & 0 & 0 & xy^2 & 0 & x^4 & 0 & -x^3y \\
  0 & 0 & 0 & 0 & -z & 0 & 0 & 0 & 0 & 0 & -xy^2 & 0 & -x^4 & 0 \\
  0 & 0 & 0 & 0 & 0 & z & y^2 & 0 & 0 & -x^2y & 0 & xy^2 & 0 & x^4 \\
  0 & 0 & -x^2 & 0 & 0 & 0 & -z & 0 & 0 & 0 & x^2y & 0 & -xy^2 & 0 \\
  0 & 0 & 0 & 0 & 0 & 0 & -xy & z & 0 & x^3 & 0 & -x^2y & 0 & xy^2 \\
  -x & 0 & y & 0 & 0 & 0 & 0 & 0 & -z & 0 & 0 & 0 & 0 & 0 \\
  0 & 0 & -y & y & 0 & 0 & 0 & x^2 & 0 & -z & y^2 & 0 & 0 & 0 \\
  0 & -x & 0 & 0 & -y & 0 & 0 & 0 & xy & 0 & z & 0 & 0 & 0 \\
  0 & 0 & 0 & x & 0 & y & 0 & 0 & 0 & 0 & 0 & -z & y^2 & 0 \\
  0 & 0 & 0 & 0 & -x & 0 & -y & 0 & 0 & 0 & 0 & 0 & z & 0 \\
  0 & 0 & 0 & 0 & 0 & x & 0 & y & 0 & 0 & 0 & 0 & 0 & -z
 \end{array}
 \right)}\ . 
\end{equation*}

\noindent
$\bullet$\ $W=(4,6,11;22)$, $A_W=(2,4,6)$. 
For the polynomial $f_W=yx^4+xy^3+z^2$, 
\begin{equation*}
 V_0 : 
{q_0}={q_1}=
\bp z & -x^4 & y^3 & 0 \\ 
    y & -z & 0 & y^3 \\ 
    x & 0 & -z & -x^4 \\
    0 & x & -y & z \ep\ ,\qquad  
 V_{1,2} :
{q_0}={q_1}=
 \bp
 z & 0 & xy^2 & x^3y \\
 0 & z & x^2y & -xy^2 \\
 y & x^2 & -z & 0 \\
 x & -y & 0 & -z
 \ep\ ,
\end{equation*}
\begin{equation*}
 V_{2,4} :
{q_0}={q_1}=
\bp z & x(x^3+y^2) \\ 
    y & -z \ep\ ,\qquad 
 V_{3,6} :
{q_0}={q_1}=
\bp z & y(x^3+y^2) \\ 
    x & -z \ep\ , 
\end{equation*}
\begin{equation*}
 \Vb: {q_0}={q_1}=
{\footnotesize 
 \bp
 z & -y^2 & 0 & 0 & -xy & 0 & -x^3y & 0 & 0 & 0\\
 -xy & -z & 0 & 0 & 0 & 0 & 0 & 0 & -x^3y & 0 \\
 0 & xy & z & 0 & -x^3 & 0 & -xy^2 & 0 & 0 & 0 \\
 0 & 0 & 0 & z & y^2 & 0 & 0 & xy^2 & 0 & x^3y \\
 0 & 0 & -xy & 0 & -z & 0 & 0 & 0 & -xy^2 & 0 \\
 0 & 0 & 0 & 0 & xy & z & 0 & x^2y & 0 & -xy^2 \\
 -x & 0 & -y & 0 & 0 & 0 & -z & 0 & 0 & 0 \\
 0 &  0 & y & y & 0 & x^2 & 0 & -z & y^2 & 0 \\
 0 & -x &  0 & 0 & -y & 0 & xy & 0 & z & 0 \\
 0 & 0 & 0 & x & 0 & -y & 0 & 0 & 0 & -z
\ep
}\ . 
\end{equation*}

\noindent
$\bullet$\ $W=(3,5,9;18)$, $A_W=(3,3,5)$. 
For the polynomial $f_W=x^3z+xy^3+z^2$, 
\begin{equation*}
 V_0 :
{q_0} ={q_1}
 =\bp
   z & xy^2 & x^2z & 0\\
   y & -z & 0 & x^2z \\
   x & 0 & -z & -xy^2 \\
   0 & x & -y & z
  \ep\ ,
\end{equation*}
\begin{equation*}
 V_{1,3} :
{q_0}=
\bp x^3+z & xy^2 \\ 
    y & -z \ep\ ,\qquad 
{q_1}= 
\bp z & xy^2 \\ 
    y & -(x^3+z) \ep\ ,
\end{equation*}
\begin{equation*}
 V_{2,3} :
{q_0}=
\bp z & xy^2 \\ 
    y & -(x^3+z) \ep\ ,\qquad 
{q_1}= 
\bp x^3+z & xy^2 \\ 
    y & -z \ep\ ,
\end{equation*}
\begin{equation*}
 V_{3,5} :
{q_0}={q_1}=
\bp z & y^3+x^2z \\ 
    x & -z \ep\ , 
\end{equation*}
\begin{equation*}
 \Vb: {q_0}
 ={q_1}
 =
{\footnotesize 
\bp
   z & y^2 & 0 & 0 & -xz & 0 & x^2z & 0 \\
   0 & -z & -y^2 & 0 & -x^2y & yz & 0 & -x^2z \\
   0 & 0 & -z & -x^3 & 0 & -xy^2 & 0 & x^3y \\
   0 & -xy & -z & z & 0 & 0 & xy^2 & 0 \\
   0 & 0 & 0 & xy & -z & -xz & 0 & -xy^2 \\
   0 & 0 & -y & 0 & -x^2 & z & 0 & 0 \\
   x & 0 & -y & y & -x^2 & 0 & -z & 0 \\
   0 & -x & 0 & 0 & -y & 0 & xy & -z
  \ep
}\ . 
\end{equation*}

\noindent
$\bullet$\ $W=(4,5,10;20)$, $A_W=(2,5,5)$. 
For the polynomial $f_W=x^5+y^2z+z^2$, 
\begin{equation*}
 V_0 :
{q_0} ={q_1}=
 \bp
 z & yz & x^4 & 0 \\
 y & -z & 0 & x^4 \\
 x & 0 & -z & -yz \\
 0 & x & -y & z
 \ep\ , \quad 
 V_{1,2} :
 q_0=q_1=
 \bp z & x^3 & yz & 0 \\
     x & -z & 0 & -yz \\
     y & 0 & -z & x^3 \\
     0 & -y & x & z \ep\ ,
\end{equation*}
\begin{equation*}
 V_{2,5} :
{q_0}=
\bp y^2+z & x^4 \\ 
    x & -z \ep\ ,\qquad 
{q_1}= 
\bp z & x^4 \\ 
    x & -(y^2+z) \ep\ ,
\end{equation*}
\begin{equation*}
 V_{3,5} :
{q_0}=
\bp z & x^4 \\ 
    x & -(y^2+z) \ep\ ,\qquad 
{q_1}= 
\bp y^2+z & x^4 \\ 
    x & -z \ep\ ,
\end{equation*}
\begin{equation*}
 \Vb: {q_0}={q_1}=
{\footnotesize 
 \bp
 -z & -y^2 & 0 & 0 & 0 & 0 & x^4 & 0 & x^3y & 0\\
 -y^2 & z & 0 & 0 & -x^3 & 0 & 0 & 0 & 0 & 0 \\
 0 & xy & -z & 0 & 0 & 0 & y^3 & 0 & -x^4 & 0 \\
 xy & 0 & 0 & z & 0 & x^3 & 0 & y^3 & 0 & 0 \\
 0 & -x^2 & 0 & 0 & -z & 0 & -xy^2 & 0 & -y^3 & 0 \\
 0 & 0 & 0 & x^2 & 0 & -z & -x^2y & 0 & -xy^2 & -y^3 \\
 x & 0 & y & 0 & 0 & 0 & z & 0 & 0 & 0 \\
 0 &  x & 0 & y & 0 & 0 & 0 & -z & 0 & x^3 \\
 0 & 0 &  -x & 0 & -y & 0 & 0 & 0 & z & 0 \\
 0 & 0 & 0 & 0 & x & -y & 0 & x^2 & 0 & z
\ep
}\ . 
\end{equation*}

\noindent
$\bullet$\ $W=(3,4,8;16)$, $A_W=(3,4,4)$. 
For the polynomial $f_W=yx^4+y^2z+z^2$, 
\begin{equation*}
 V_0 :
{q_0} ={q_1}
 =\bp
   z & yz & x^3y & 0 \\
   y & -z & 0 & x^3y \\
   x & 0 & -z & -yz \\
   0 & x & -y & z
  \ep\ ,\quad  
 V_{1,3} :
{q_0}={q_1}=
 \bp z & x^4+yz \\
     y & -z \ep\ ,
\end{equation*}
\begin{equation*}
 V_{2,4} :
{q_0}=
\bp y^2+z & x^3y \\ 
    x & -z \ep\ ,\qquad 
{q_1}= 
\bp z & x^3y \\ 
    x & -(y^2+z) \ep\ ,
\end{equation*}
\begin{equation*}
 V_{3,4} :
{q_0}=
\bp z & x^3y \\ 
    x & -(y^2+z) \ep\ ,\qquad 
{q_1}= 
\bp y^2+z & x^3y \\ 
    x & -z \ep\ ,
\end{equation*}
\begin{equation*}
 \Vb: {q_0} ={q_1}
 =
{\footnotesize 
\bp
   z & z & 0 & 0 & x^2y & 0 & -x^3y & 0  \\
   y^2 & -z & 0 & 0 & 0 & x^3y & 0 & x^2y^2 \\
   xy & 0 & -z & 0 & 0 & -yz & 0 & x^3y \\
   0 & 0 & z & z & x^3 & 0 & yz & 0 \\
   0 & 0 & 0 &  xy & -z & 0 & 0 & -yz \\
   0 & x & -y & 0 & 0 & z & 0 & 0 \\
   -x & 0 & 0 & y & 0 & z & -z & 0 \\
   0 & 0 & x & 0 & -y & 0 & xy & z 
  \ep
}\ . 
\end{equation*}

\noindent
$\bullet$\ $W=(6,8,9;24)$, $A_W=(2,3,9)$. 
For the polynomial $f_W=x^4+y^3+xz^2$, 
\begin{equation*}
 V_0 : 
{q_0}=
\bp xz & y^2 & x^3 & 0 \\ 
    y & -z & 0 & x^3 \\ 
    x & 0 & -z & -y^2 \\
    0 & x & -y & xz \ep\ ,\qquad 
{q_1}=
 \bp z & y^2 & x^3 & 0 \\ 
     y & -xz & 0 & x^3 \\
     x & 0 & -xz & -y^2 \\
     0 & x & -y & z 
 \ep\ , 
\end{equation*}
\begin{equation*}
 V_{1,2} :
{q_0}=
{\small
\bp xz & y^2 & 0 & x^3 & 0 & -x^2y \\ 
    0 & -xy & xz & y^2 & 0 & x^3 \\
    0 & x^2 & 0 & -xy & xz & y^2 \\
    y & -z & 0 & 0 & x^2 & 0 \\ 
    x & 0 & y & -z & 0 & 0 \\
    0 & 0 & x & 0 & y & -z \ep }\ ,\qquad 
{q_1}=
{\small
 \bp z & 0 & 0 & y^2 & x^3 & -x^2y \\ 
     y & 0 & x^2 & -xz & 0 & 0 \\ 
     0 & z & 0 & -xy & y^2 & x^3 \\
     x & y & 0 & 0 & -xz & 0  \\
     0 & 0 & z & x^2 & -xy & y^2 \\
     0 & x & y & 0 & 0 & -xz \ep
}\ , 
\end{equation*}
\begin{equation*}
 V_{2,3} :
{q_0}=
\bp xz & y^2 & -x^3 & 0 \\ 
    x^2 & 0 & xz & -y^2 \\
    y & -z & 0 & x^2 \\ 
    0 & x & y & z \ep\ ,\qquad 
{q_1}=
 \bp z & x^2 & y^2 & 0 \\ 
     y & 0 & -xz & x^3 \\ 
     -x & z & 0 & y^2 \\
     0 & -y & x^2 & xz \ep\ , 
\end{equation*}
\begin{equation*}
 V_{3,9}  :
{q_0}=
\bp y^2 & x^3+z^2 \\ x & -y \ep\ ,\qquad 
{q_1}=
 \bp y & x^3+z^2 \\ x & -y^2 \ep\ , 
\end{equation*}
\begin{equation*}
 \Vb : 
{q_0}=
{\footnotesize 
 \left(
 \begin{array}{ccccccccccccccc}
 -y^2 & -yz & z^2 & -x^3 & 0 & 0 & 0 & 0 & 0 & 0 & xz^2 & 0 & 0 & 0 & 0 \\
 xz & -y^2 & yz & 0 & 0 & -x^3 & 0 & 0 & 0 & 0 & x^2y & 0 & 0 & 0 & 0 \\
 xy & xz & y^2 & 0 & 0 & 0 & 0 & 0 & -x^3 & 0 & 0 & 0 & x^2y & 0 & 0 \\
 0 & -xy & xz & xz & 0 & y^2 & 0 & 0 & 0 & 0 & x^3 & 0 & 0 & 0 & 0 \\
 0 & xy & 0 & 0 & xz & 0 & 0 & y^2 & 0 & 0 & x^3 & x^3 & 0 & 0 & -x^2y \\
 x^2 & 0 & xy & 0 & 0 & 0 & xz & 0 & y^2 & 0 & 0 & 0 & x^3 & 0 & 0 \\
 0 & -x^2 & 0 & 0 & 0 & 0 & 0 & -xy & 0 & xz & 0 & y^2 & 0 & 0 & x^3 \\
 0 & 0 & -x^2 & -x^2 & 0 & 0 & 0 & 0 & -xy & 0 & xz & 0 & y^2 & 0 & 0 \\
 0 & 0 & 0 & 0 & 0 & 0 & x^2 & x^2 & 0 & 0 & 0 & -xy & 0 & xz & y^2 \\
 -x & 0 & 0 & y & 0 & -z & 0 & 0 & 0 & 0 & 0 & 0 & 0 & 0 & 0 \\
 x & 0 & y & 0 & y & 0 & 0  &-z & 0 & 0 & 0 & 0 & x^2 & x^2 & 0 \\
 0 & -x & 0 & 0 & 0 & 0 & y & 0 & -z & 0 & 0 & 0 & 0 & 0 & 0 \\
 0 & 0 & 0 & 0 & x & 0 & 0 & 0 & 0 & y & 0 & -z & 0 & 0 & 0 \\
 0 & 0 & 0 & 0 & 0 & -x & x & 0 & 0 & 0 & y & 0 & -z & 0 & 0 \\
 0 & 0 & 0 & 0 & 0 & 0 & 0 & 0 & x & x & 0 & 0 & 0 & y & -z 
 \end{array}
 \right)\ ,
}
\end{equation*}
\begin{equation*}
{q_1}=
 {\footnotesize 
 \left(
 \begin{array}{ccccccccccccccc}
 -y & z & 0 & 0 & 0 & 0 & 0 & 0 & 0 & -x^3 & 0 & 0 & 0 & 0 & 0 \\
 0 &  -y & z & 0 & 0 & 0 & 0 & 0 & 0 & 0 & 0 & -x^3 & 0 & x^2y & 0  \\
 x & 0 & y & 0 & 0 & 0 & 0 & -x^2 & 0 & 0 & 0 & 0 & 0 & 0 & 0 \\
 -x & 0 & 0 & z & 0 & 0 & 0 & 0 & 0 & y^2 & 0 & 0 & 0 & 0 & 0 \\
 0 & 0 & -y & 0 & z & 0 & 0 & 0 & 0 & 0 & y^2 & 0 & x^3 & 0 & -x^2y \\
 0 & -x & 0 & y & 0 & 0 & 0 & 0 & 0 & -xz & 0 & 0 & 0 & 0 & 0 \\
 0 & -x & 0 & 0 & 0 & z & 0 & 0 & 0 & 0 & 0 & y^2 & 0 & x^3 & 0 \\
 0 & 0 & 0 & 0 & y & 0 & 0 & 0 & x^2 & 0 & -xz & 0 & 0 & -x^3 & 0 \\
 0 & 0 & -x & 0 & 0 & y & 0 & 0 & 0 & 0 & 0 & -xz & 0 & 0 & 0 \\
 0 & 0 & x & 0 & 0 & 0 & z & 0 & 0 & 0 & -xy & 0 & y^2 & 0 & x^3 \\
 0 & 0 & 0 & x & 0 & 0 & 0 & z & 0 & 0 & 0 & -xy & 0 & y^2 & 0 \\
 0 & 0 & 0 & 0 & x & 0 & y & 0 & 0 & 0 & 0 & 0 & -xz & 0 & 0 \\
 0 & 0 & 0 & 0 & 0 & x & 0 & y & 0 & x^2 & 0 & 0 & 0 & -xz & 0 \\
 0 & 0 & 0 & 0 & 0 & -x & 0 & 0 & z & 0 & x^2 & 0 & -xy & 0 & y^2 \\
 0 & 0 & 0 & 0 & 0 & 0 & x & 0 & x & 0 & 0 & -x^2 & 0 & 0 & -xz
 \end{array}
 \right)}\ . 
\end{equation*}

\noindent
$\bullet$\ $W=(4,6,7;18)$, $A_W=(2,4,7)$. 
For the polynomial $f_W=x^3y+y^3+xz^2$, 
\begin{equation*}
 V_0 :
{q_0}=
\bp
 xz & y^2 & x^2y & 0 \\
 y & -z & 0 & x^2y \\
 x & 0 & -z & -y^2 \\
 0 & x & -y & xz
\ep\ ,\qquad 
 {q_1}=
\bp
 z & y^2 & x^2y & 0 \\
 y & -xz & 0 & x^2y \\
 x & 0 & -xz & -y^2 \\
 0 & x & -y & z
\ep\ ,
\end{equation*}
\begin{equation*}
 V_{1,2} :
{q_0}=
\bp xz & y^2 & 0 & x^2y \\
    0 & -xy & xz & y^2 \\
    y & -z & -x^2 & 0 \\
    x & 0 & y & -z  \ep\ ,\qquad 
{q_1}=
\bp z & 0 & y^2 & -x^2y \\
    y & -x^2 & -xz & 0 \\
    0 & z & -xy & y^2 \\
    x & y & 0 & -xz \ep\ , 
\end{equation*}
\begin{equation*}
 V_{2,4} :
{q_0}=
\bp xz & x^3+y^2 \\ y & -z \\ \ep\ ,\qquad 
{q_1}=
 \bp z & x^3+y^2 \\ y & -xz \\ \ep\ , 
\end{equation*}
\begin{equation*}
 V_{3,7} :
{q_0}=
\bp x^2y+z^2 & y^2 \\ y & -x \\ \ep\ ,\qquad 
{q_1}=
 \bp x & y^2 \\ y & -(x^2y+z^2) \\ \ep\ , 
\end{equation*}
\begin{equation*}
 \Vb\ : 
{q_0}=
 {\footnotesize 
\left(
 \begin{array}{ccccccccccc}
 -y^2 & -yz & z^2 & -x^2y & 0 & 0 & x^2z & 0 & 0 & 0 & 0 \\
  xz & -y^2 & yz & 0 & 0 & -x^2y & x^2y & 0 & 0 & 0 & 0 \\
  xy & xz & y^2 & 0 & 0 & 0 & 0 & 0 & x^2y & 0 & 0 \\
  0 & -xy & xz & xz & 0 & y^2 & x^3 & 0 & 0 & 0 & 0 \\
  0 & -xy & 0 & 0 & xz & xy & 0 & y^2 & 0 & 0 & x^2y \\
  0 & 0 & -xy & -xy & 0 & 0 & xz & 0 & y^2 & 0 & 0 \\
  0 & 0 & 0 & 0 & 0 & 0 & -xy & -xy & 0 & xz & y^2 \\
  -x & 0 & 0 & y & 0 & -z & 0 & 0 & 0 & 0 & 0 \\
  0 & 0 & -y & -y & y & 0 & 0 & -z & 0 & -x^2 & 0  \\
  0 & x & 0 & 0 & 0 & -y & y & 0 & -z & 0 & 0 \\
  0 & 0 & 0 & 0 & x & 0 & 0 & 0 & -y & y & -z 
 \end{array}
 \right)
}\ , 
\end{equation*}
\begin{equation*}
{q_1}=
 {\footnotesize 
\left(
 \begin{array}{ccccccccccc}
  -y & z & 0 & 0 & 0 & 0 & 0 & -x^2y & 0 & 0 & 0 \\
  0 & -y & z & 0 & 0 & 0 & 0 & 0 & 0 & x^2y & 0 \\
  x & 0 & y & 0 & 0 & -x^2 & 0 & 0 & 0 & 0 & 0 \\
  -x & 0 & 0 & z & 0 & 0 & 0 & y^2 & 0 & 0 & 0 \\
  0 & 0 & y & 0 & z & 0 & 0 & 0 & y^2 & 0 & x^2y \\
  0 & -x & 0 & y & 0 & 0 & 0 & -xz & 0 & 0 & 0 \\
  0 & 0 & 0 & y & 0 & z & 0 & 0 & 0 & y^2 & 0 \\
  0 & -x & 0 & 0 & y & 0 & -x^2 & 0 & -xz & 0 & 0 \\
  0 & 0 & x & 0 & 0 & y & 0 & xy & 0 & -xz & 0  \\
  0 & 0 & 0 & 0 & 0 & y & z & xy & -xy & 0 & y^2  \\
  0 & 0 & 0 & 0 & x & 0 & y & 0 & 0 & xy & -xz
 \end{array}
 \right)
}\ . 
\end{equation*}

\noindent
$\bullet$\ $W=(3,5,6;15)$, $A_W=(3,3,6)$. 
For the polynomial $f_W=x^3z+y^3+xz^2$, 
\begin{equation*}
 V_0 :
{q_0}=
\bp
 xz & y^2 & x^2z & 0 \\
 y & -z & 0 & x^2z \\
 x & 0 & -z & -y^2 \\
 0 & x & -y & xz
\ep\ ,\qquad 
 {q_1}=
\bp
 z & y^2 & x^2z & 0 \\
 y & -xz & 0 & x^2z \\
 x & 0 & -xz & -y^2 \\
 0 & x & -y & z
\ep\ ,
\end{equation*}
\begin{equation*}
 V_{1,3} :
{q_0}=
\bp xz & y^2 \\ y & -(x^2+z) \ep\ ,\qquad 
{q_1}=
 \bp x^2+z & y^2 \\ y & -xz \ep\ , 
\end{equation*}
\begin{equation*}
 V_{2,3} :
{q_0}=
\bp x(x^2+z) & y^2 \\ y & -z \ep\ ,\qquad 
{q_1}=
 \bp z & y^2 \\ y & -x(x^2+z) \ep\ , 
\end{equation*}
\begin{equation*}
 V_{3,6} :
{q_0}=
\bp y^2 & z(x^2+z) \\ x & -y \ep\ ,\qquad 
{q_1}=
 \bp y & z(x^2+z) \\ x & -y^2 \ep\ ,
\end{equation*}
\begin{equation*}
 \Vb: 
{q_0}=
{\footnotesize 
\bp
 -y^2 & -yz & x^z+z^2 & 0 & 0 & 0 & 0 & 0 & xyz \\
 xz & -y^2 & x^y+yz & x^2y & 0 & -x^2 & 0 & 0 & xy^2 \\
 xy & xz & y^2 & 0 & 0 & 0 & 0 & 0 & -x^2z \\
 0 & -xy & x^3+xz & x^3+xz & 0 & y^2 & 0 & 0 & x^2y \\
 0 & 0 & -x^3 & -x^3 & xz & 0 & 0 & y^2 & -x^2y \\
 x^2 & 0 & xy & 0 & 0 & 0 & xz & 0 & y^2 \\
 -x & 0 & 0 & y & 0 & -z & 0 & 0 & 0 \\
 0 & 0 & 0 & 0 & y & -x^2 & x^2 & -x^2-z & 0 \\
 0 & -x & 0 & 0 & 0 & 0 & y & 0 & -z 
\ep
}\ ,
\end{equation*}
\begin{equation*}
 {q_1}=
{\footnotesize 
\bp
 -y & z & 0 & 0 & 0 & 0 & -x^2z & 0 & 0 \\
 0 & -y & z & 0 & 0 & xy & x^2y & 0 & -x^2z \\
 x &0 & y & 0 & 0 & 0 & 0 & 0 & 0 \\
 -x & 0 & 0 & z & 0 & 0 & y^2 & 0 & 0 \\
 0 & 0 & 0 & x^2 & x^2+z & 0 & 0 & y^2 & -x^2y \\
 0 & -x & 0 & y & 0 & 0 & -xz & 0 & 0 \\
 0 & -x & 0 & 0 & 0 & x^2+z & x^3 & 0 & y^2 \\
 0 & 0 & 0 & 0 & y & x^2 & x^3 & -xz & 0 \\
 0 & 0 & -x & 0 & 0 & y & 0 & 0 & -xz 
\ep
}\ . 
\end{equation*}

\noindent
$\bullet$\ $W=(4,5,6;16)$, $A_W=(2,5,6)$. 
For the polynomial $f_W=x^4+y^2z+xz^2$, 
\begin{equation*}
 V_0 :
{q_0}=
\bp
 xz & yz & x^3 & 0 \\
 y & -z & 0 & x^3 \\
 x & 0 & -z & -yz \\
 0 & x & -y & xz
\ep\ ,\qquad 
 {q_1}=
\bp
 z & yz & x^3 & 0 \\
 y & -xz & 0 & x^3 \\
 x & 0 & -xz & -yz \\
 0 & x & -y & z
\ep\ ,
\end{equation*}
\begin{equation*}
 V_{1,2} :
{q_0}=
\bp xz & yz & x^3 & 0 \\ 
    -x^2 & 0 & xz & yz \\
    -y & z & 0 & -x^2 \\ 
    0 & x & -y & z \ep\ ,\qquad 
{q_1}=
 \bp z & -x^2 & -yz & 0 \\
     y & 0 & xz & x^3 \\
     x & z & 0 & -yz \\ 
     0 & y & -x^2 & xz  \ep\ , 
\end{equation*}
\begin{equation*}
 V_{2,5} :
{q_0}=
\bp y^2+xz & x^3 \\ x & -z  \\ \ep\ ,\qquad 
{q_1}=
 \bp z & x^3 \\ x & -(y^2+xz) \\ \ep\ , 
\end{equation*}
\begin{equation*}
 V_{3,6}  :
{q_0}=
\bp yz & x^3+z^2  \\ x & -y \ep\ ,\qquad 
{q_1}=
 \bp y & x^3+z^2 \\ x & -yz \ep\ , 
\end{equation*}
\begin{equation*}
 \Vb : 
{q_0}=
{\footnotesize 
\left(
\begin{array}{ccccccccccc}
 -yz & yz & x^3+z^2 & -x^3 & 0 & 0 & 0 & 0 & -x^2z & 0 & 0 \\
    -xz & xz & -yz & 0 & 0 & 0 & -x^3 & 0 & 0 & 0 & 0 \\
    xz & y^2 & yz & 0 & 0 & -x^3 & 0 & 0 & -x^2y & 0 & 0 \\
    0 & -xy & 0 & -xz & 0 & -yz & 0 & 0 & x^3 & 0 & 0 \\
    0 & xy & 0 & 0 & xz & 0 & 0 & yz & -x^3 & x^3 & 0 \\
    0 & x^2 & 0 & 0 & 0 & xz & xz & 0 & yz & 0 & 0 \\
    0 & 0 & 0 & 0 & -x^2 & 0 & 0 & 0 & -xz & xz & yz \\
    -x & 0 & -y & y & 0 & -z & 0 & 0 & 0 & 0 & 0 \\
    0 & x & 0 & 0 & -y & 0 & 0 & z & 0 & 0 & -x^2 \\
    0 & 0 & 0 & x & 0 & 0 & -y & 0 & z & 0 & 0  \\
    0 & 0 & 0 & 0 & 0 & -x & -x & x & 0 & -y & z 
\end{array}\right)
}\ ,
\end{equation*}
\begin{equation*}
\qquad {q_1}=
{\footnotesize 
\left(
 \begin{array}{ccccccccccc}
     -y & 0 & z & 0 & 0 & 0 & 0 & -x^3 & 0 & 0 & 0 \\
     0 & z & z & 0 & 0 & x^2 & 0 & 0 & 0 & 0 & 0 \\
     x & -y & 0 & 0 & 0 & 0 & 0 & 0 & 0 & x^3 & 0 \\
     0 & -y & 0 & -z & 0 & 0 & 0 & yz & 0 & x^3 & 0 \\
     0 & 0 & 0 & -z & z & 0 & -x^2 & 0 & -yz & 0 & 0 \\
     0 & 0 & -x & -y & 0 & 0 & 0 & -xz & 0 & 0 & 0 \\
     0 & -x & 0 & 0 & 0 & z & 0 & xz & 0 & -yz & 0 \\
     0 & -x & 0 & 0 & y & 0 & 0 & xz & xz & 0 & x^3 \\
     0 & 0 & 0 & x & 0 & y & 0 & 0 & 0 & xz & 0 \\
     0 & 0 & 0 & 0 & x & 0 & z & 0 & 0 & xz & -yz \\
     0 & 0 & 0 & 0 & 0 & x & y & 0 & -x^2 & 0 & xz
 \end{array}\right)
}\ . 
\end{equation*}

\noindent
$\bullet$\ $W=(3,4,5;13)$, $A_W=(3,4,5)$. 
For the polynomial $f_W=x^3y+y^2z+xz^2$, 
\begin{equation*}
 V_0 :
{q_0}=
\bp
 xz & yz & x^2y & 0 \\
 y & -z & 0 & x^2y \\
 x & 0 & -z & -yz \\
 0 & x & -y & xz
\ep\ , \qquad 
 {q_1}=
\bp
 z & yz & x^2y & 0 \\
 y & -xz & 0 & x^2y \\
 x & 0 & -xz & -yz \\
 0 & x & -y & z
\ep\ , 
\end{equation*}
\begin{equation*}
 V_{1,3} :
{q_0}=
\bp xz+y^2 & x^3 \\ y & -z \ep\ ,\qquad 
{q_1}=
 \bp z & x^3 \\ y & -(xz+y^2) \ep\ , 
\end{equation*}
\begin{equation*}
 V_{2,4} :
{q_0}=
\bp y^2 & x^2y+z^2 \\ x & -z \ep\ ,\qquad 
{q_1}=
 \bp z & x^2y+z^2 \\ x & -y^2 \ep\ , 
\end{equation*}
\begin{equation*}
 V_{3,5} :
 {q_0}=
\bp yz & x^2y+z^2 \\ x & -y \ep\ ,\qquad 
{q_1}=
 \bp y & x^2y+z^2 \\ x & -yz \ep\ , 
\end{equation*}
\begin{equation*}
 \Vb : 
 {q_0}=
{\footnotesize 
\bp 
 yz & yz & 0 & -x^2y-z^2 & 0 & 0 & 0 & xy^2 & 0  \\
 xz & -y^2 & 0 & yz & 0 & -x^2y & 0 & 0 & xy^2 \\
 xz & xz & 0 & yz & 0 & 0 & 0 & x^2y & 0 \\
 0 & -xy & xz & xz & 0 & yz & 0 & 0 & x^2y \\
 0 & 0 & 0 & 0 & xz+y^2 & 0 & x^3 & 0 & -x^2y \\
 0 & -x^2 & -xy & -xy & 0 & xz & 0 & y^2+xz & 0 \\
 -x & 0 & y & 0 & 0 & -z & 0 & 0 & 0 \\
 0 & -x & 0 & 0 & y & z & -z & z & 0 \\
 0 & 0 & -x & -x & 0 & 0 & 0 & y & z
\ep 
}\ , 
\end{equation*}
\begin{equation*}
{q_1}=
{\footnotesize 
 \bp 
 y & z & 0 & 0 & 0 & 0 & -x^2y & 0 & -xy^2 \\
 0 & -z & z & 0 & 0 & -xy & 0 & 0 & xy^2 \\
 x & 0 & 0 & z & 0 & 0 & yz & 0 & -x^2y \\
 -x & 0 & y & 0 & 0 & 0 & 0 & 0 & 0 \\
 0 & 0 & 0 & 0 & z & -x^2 & 0 & x^3 & x^2y \\
 0 & -x & 0 & y & 0 & 0 & -xz & 0 & 0 \\
 0 & 0 & 0 & 0 & y & z & 0 & -y^2-xz & 0  \\
 0 & 0 & x & 0 & 0 & z & xz & 0 & 0 \\
 0 & 0 & 0 & x & 0 & -y & 0 & 0 & y^2+xz
\ep
}\ . 
\end{equation*}

\noindent
$\bullet$\ $W=(3,4,4;12)$, $A_W=(4,4,4)$. 
For the polynomial $f_W=x^4+yz(y-z)$, 
\begin{equation*}
 V_0 : 
{q_0}=
\bp -yz & yz & x^3 & 0 \\ 
    y & -z & 0 & x^3 \\ 
    x & 0 & -z & -yz \\
    0 & x & -y & -yz \ep\ ,\qquad 
{q_1}=
 \bp z & yz & x^3 & 0 \\ 
     y & yz & 0 & x^3 \\
     x & 0 & yz & -yz \\
     0 & x & -y & z 
 \ep\ ,
\end{equation*}
\begin{equation*}
 V_{1,4} :
{q_0}=
\bp yz & x^3 \\ x & -(y-z) \ep\ ,\qquad 
{q_1}=
 \bp (y-z) & x^3 \\ x & -yz \ep\ , 
\end{equation*}
\begin{equation*}
 V_{2,4} :
{q_0}=
\bp (y-z)z & x^3 \\ x & -y \ep\ ,\qquad 
{q_1}=
 \bp y & x^3 \\ x & -(y-z)z \ep\ , 
\end{equation*}
\begin{equation*}
 V_{3,4} :
{q_0}=
\bp y(y-z) & x^3 \\ x & -y \ep\ ,\qquad 
{q_1}=
 \bp y & x^3 \\ x & -y(y-z) \ep\ , 
\end{equation*}
\begin{equation*}
 \Vb : 
{q_0}=
{\footnotesize 
\bp 
 -(y-z)z & 0 & 0 & -x^3 & 0 & 0 & 0 & -x^2z & 0 \\
 0 & y(y-z) & 0 & 0 & -x^3 & 0 & 0 & -x^2y & 0 \\
 0 & -y^2 & yz & 0 & 0 & 0 & x^3 & 0 & x^2y \\
 xz & xy & 0 & -yz & yz & 0 & 0 & -x^3 & 0 \\
 0 & xy & 0 & 0 & 0 & yz & 0 & 0 & -x^3 \\
 -x^2 & 0 & 0 & xy & 0 & 0 & 0 & yz & 0 \\
 -x & -x & 0 & y & -z & 0 & 0 & 0 & 0 \\
 0 & 0 & x & 0 & 0 & y & -(y-z) & 0 & 0 \\
 0 & 0 & 0 & 0 & x & -x & 0 & y & -(y-z)
\ep
}\ ,
\end{equation*}
\begin{equation*}
{q_1}=
{\footnotesize  
\bp 
 -y & 0 & 0 & 0 & 0 & -x^2 & 0 & 0 & 0 \\
 0 & z & 0 & 0 & 0 & x^2 & -x^3 & 0 & 0 \\
 0 &  y & y-z & 0 & 0 & 0 & 0 & x^3 & x^2y \\
 -x & 0 & 0 & z & 0 & 0 & yz & 0 & 0 \\
 0 & -x & 0 & y & 0 & 0 & yz & 0 & 0 \\
 0 & -x & 0 & 0 & y-z & 0 & 0 & 0 & -x^3 \\
 0 & 0 & x & 0 & y & 0 & 0 & -yz & 0  \\
 0 & 0 & 0 & -x & 0 & y-z & -xy & 0 & 0 \\
 0 & 0 & 0 & 0 & -x & y & -xy & 0 & -yz
 \ep 
}\ . 
\end{equation*}

\noindent
$\bullet$\ $W=(2,6,9;18)$, $A_W=(2,2,2,3)$. 
For the polynomial 
$f_W=yy'y''+z^2$, $y':=y-x^3$, $y'':=y-\lambda x^3$, $\lambda\ne 0,1$, 
\begin{equation*}
 V_0 : 
{q_0}={q_1}=
\bp z & yy'' & -x^2yy'' & 0 \\ 
    y & -z & 0 & -x^2yy'' \\ 
    x & 0 & -z & -yy'' \\
    0 & x & -y & z \ep\ , 
\end{equation*}
\begin{equation*}
 V_{1,2} :
{q_0}={q_1}=
\bp z & y'y'' \\ y & -z \ep\ ,\quad 
 V_{2,2} :
{q_0}={q_1}=
\bp z & yy'' \\ y' & -z \ep\ ,\quad 
 V_{3,2} :
{q_0}={q_1}=
\bp z & yy' \\ y'' & -z \ep\ ,
\end{equation*}
\begin{equation*}
 V_{4,3} :
 q_0=q_1= q_0(V_0)=q_1(V_0)\ ,  
\end{equation*}
\begin{equation*}
 \Vb : 
{q_0}={q_1}=
{\small
\bp z & -y^2 & 0 & 0 & x^2yy'' & 0 \\ 
    -y'' & -z & x^2y'' & 0 & 0 & 0 \\ 
    0 & -xy & z & 0 & yy'' & 0 \\
    0 & x^4 & 0 & z & 0 & y'y'' \\
    -x & 0 & y & 0 & z & 0 \\
    0 & 0 & x^3 & y & 0 & -z \ep
}\ . 
\end{equation*}

\noindent
$\bullet$\ $W=(2,4,7;14)$, $A_W=(2,2,2,4)$. 
For the polynomial 
$f_W=xyy'y''+z^2$, $y':=y-x^2$, $y'':=y-\lambda x^2$, $\lambda\ne 0,1$, 
\begin{equation*}
 V_0 : 
{q_0}={q_1}=
\bp 
 z & xyy'' & -x^2yy'' & 0 \\ 
 y & -z & 0 & -x^2yy'' \\ 
 x & 0 & -z & -xyy'' \\
 0 & x & -y & z 
\ep\ , 
\end{equation*}
\begin{equation*}
 V_{1,2} :
{q_0}={q_1}=
\bp z & xy'y'' \\ y & -z \ep\ ,\quad 
 V_{2,2} :
{q_0}={q_1}=
\bp z & xyy'' \\ y' & -z \ep\ , 
\end{equation*}
\begin{equation*}
 V_{3,2} :
{q_0}={q_1}=
\bp z & xyy' \\ y'' & -z \ep\ ,\quad 
 V_{4,4} :
{q_0}={q_1}=
\bp z & yy'y'' \\ x & -z \ep\ , 
\end{equation*}
\begin{equation*}
 \Vb : {q_0}={q_1}=
{\small
 \bp 
 x & -y^2 & 0 & 0 & x^2yy'' & 0 \\ 
 -xy'' & -z & x^2y'' & 0 & 0 & 0 \\
 0 & -xy & z & 0 & xyy'' & 0 \\
 0 & x^3 & 0 & z & -x^3y'' & xy'y'' \\
 -x & 0 & y & 0 & -z & 0 \\
 0 & 0 & x^2 & y & 0 & -z
 \ep
}\ . 
\end{equation*}

\noindent
$\bullet$\ $W=(2,4,5;12)$, $A_W=(2,2,2,5)$. 
For the polynomial $f_W=yy'y''+xz^2$, $y':=y-x^2$, $y'':=y-\lambda x^2$, 
$\lambda\ne 0,1$, 
\begin{equation*}
 V_0 : 
{q_0}=
\bp 
 xz & yy'' & -xyy'' & 0 \\ 
 y & -z & 0 & -xyy'' \\ 
 x & 0 & -z & -yy'' \\
 0 & x & -y & xz 
 \ep\ ,\quad 
{q_1}=
 \bp 
 z & yy'' & -xyy'' & 0 \\ 
 y & -xz & 0 & -xyy'' \\
 x & 0 & -xz & -yy'' \\
 0 & x & -y & z 
 \ep\ , 
\end{equation*}
\begin{equation*}
 V_{1,2} :
{q_0}=
\bp xz & y'y'' \\ y & -z \\ \ep\ ,\qquad 
{q_1}= 
\bp z & y'y'' \\ y & -xz \\ \ep\ , 
\end{equation*}
\begin{equation*}
 V_{2,2} :
{q_0}=
\bp xz & yy'' \\ y' & -z \\ \ep\ ,\qquad 
{q_1}= 
\bp z & yy'' \\ y' & -xz \\ \ep\ , 
\end{equation*}
\begin{equation*}
 V_{3,2} :
{q_0}=
\bp xz & yy' \\ y'' & -z \ep\ ,\qquad 
{q_1}= 
\bp z & yy' \\ y'' & -xz \ep\ , 
\end{equation*}
\begin{equation*}
 V_{4,5} :
{q_0}=
\bp y'y'' & z^2 \\ x & -y \ep\ ,\qquad 
{q_1}= 
\bp y & z^2 \\ x & -y'y'' \ep\ , 
\end{equation*}
\begin{equation*}
 \Vb : 
{q_0}=
{\footnotesize 
\bp 
 yy'' & yz & -xyy'' & -z^2 & 0 & 0 & 0 \\
 xz & -y^2 & 0 & yz & 0 & xyy'' & 0 \\
 xy'' & xz & -x^2y'' & y'y'' & 0 & 0 & 0 \\
 0 & -xy & xz & xz & 0 & yy'' & 0 \\
 0 & -x^3 & 0 & 0 & xz & x^2y'' & y'y'' \\
 -x & 0 & y & 0 & 0 & -z & 0 \\ 
 0 & 0 & -x^2 & -x^2 & y & 0 & -z 
\ep
}\ , 
\end{equation*}
\begin{equation*}
{q_1}=
{\footnotesize 
 \bp 
 y & z & 0 & 0 & 0 & xyy'' & 0 \\
 0 & -y'' & z & xy'' & 0 & 0 & 0 \\
 x & 0 & 0 & z & 0 & yy'' & 0  \\
 -x & 0 & y & 0 & 0 & 0 & 0 \\
 0 & 0 & x^2 & 0 & z & x^2y' & y'y'' \\
 0 & -x & 0 & y & 0 & -xz & 0 \\
 0 & 0 & 0 & -x^2 & y & 0 & -xz 
 \ep 
}\ . 
\end{equation*}

\def \ll{\sqrt{-\lambda}}

\noindent
$\bullet$\ $W=(2,3,6;12)$, $A_W=(2,2,3,3)$. 
For the polynomial $f_W=YY'+z^2$, $Y:=y^2-x^3$, $Y':=y^2-\lambda x^3$, 
$\lambda\ne 0,1$, 
\begin{equation*}
 V_0 :
 {q_0}={q_1}= 
  \bp
 z & yY' & -x^2Y' & 0 \\
 y & -z & 0 & -x^2Y' \\
 x & 0 & -z & -yY' \\
 0 & x & -y & z
 \ep\ , 
\end{equation*}
\begin{equation*}
 V_{1,2} :
 {q_0}=
 \bp (z+\sqrt{-\lambda}x^3) & y(Y'-x^3) \\ 
     y & -(z-\sqrt{-\lambda}x^3) \ep\ ,\qquad 
 {q_1}=
 \bp (z-\sqrt{-\lambda}x^3) & y(Y'-x^3) \\ 
     y & -(z+\sqrt{-\lambda}x^3) \ep\ , 
\end{equation*}
\begin{equation*}
 V_{2,2} :
 {q_0}=
 \bp (z-\sqrt{-\lambda}x^3) & y(Y'-x^3) \\ 
     y & -(z+\sqrt{-\lambda}x^3) \ep\ ,\qquad 
 {q_1}=
 \bp (z+\sqrt{-\lambda}x^3) & y(Y'-x^3) \\ 
     y & -(z-\sqrt{-\lambda}x^3) \ep\ ,  
\end{equation*}
\begin{equation*}
 V_{3,3} :
 {q_0}=
 \bp z+\sqrt{-1}y^2 & -x^2(\lambda y^2+Y') \\ 
     x & -(z-\sqrt{-1}y^2)\ep\ ,\qquad 
 {q_1}=
 \bp z-\sqrt{-1}y^2 & -x^2(\lambda y^2+Y') \\ 
     x & -(z+\sqrt{-1}y^2)\ep\ ,
\end{equation*}
\begin{equation*}
 V_{4,3} :
 {q_0}=
 \bp z-\sqrt{-1}y^2 & -x^2(\lambda y^2+Y') \\ 
     x & -(z+\sqrt{-1}y^2)\ep\ ,\qquad 
 {q_1}=
 \bp z+\sqrt{-1}y^2 & -x^2(\lambda y^2+Y') \\ 
     x & -(z-\sqrt{-1}y^2)\ep\ ,
\end{equation*}
\begin{equation*}
 \Vb: {q_0}=
{\small  
\bp
 z & Y' & -x^2y & \ll x^2y & 0 & x^2Y' \\
 y^2 & -z & \ll x^2y & 0 & -x^2Y' & 0 \\
 0 & 0 & -(z-\ll x^3) & y^2 & 0 & 0 \\
 0 & 0 & Y'-x^3 & z+\ll x^3 & 0 & 0 \\
 0 & x & -y & 0 & z & y^2 \\
 -x & 0 & 0 & y & Y' & -z
 \ep
}\ , 
\end{equation*}
\begin{equation*}
 {q_1}=
{\small
 \bp
 z & Y' & -x^2y & -\ll x^2y & 0 & x^2Y' \\
 y^2 & -z & -\ll x^2y & 0 & -x^2Y' & 0 \\
 0 & 0 & -(z+\ll x^3) & y^2 & 0 & 0 \\
 0 & 0 & Y'-x^3 & z-\ll x^3 & 0 & 0 \\
 0 & x & -y & 0 & z & y^2 \\
 -x & 0 & 0 & y & Y' & -z
 \ep 
}\ . 
\end{equation*}

\noindent
$\bullet$\ $W=(2,3,4;10)$, $A_W=(2,2,3,4)$. 
For the polynomial 
$f_W=xz'z''+y^2z$, $z':=z-x^2$, $z'':=z-\lambda x^2$, 
$\lambda\ne 0,1$, 
\begin{equation*}
 V_0 : 
{q_0}=
\bp 
 xz'' & yz & -x^2z'' & 0 \\ 
 y & -z & 0 & -x^2z'' \\ 
 x & 0 & -z & -yz \\
 0 & x & -y & xz'' \ep\ ,\quad 
{q_1}=
 \bp 
 z & yz & -x^2z'' & 0 \\ 
 y & -xz'' & 0 & -x^2z'' \\
 x & 0 & -xz'' & -yz \\
 0 & x & -y & z 
 \ep\ , 
\end{equation*}
\begin{equation*}
 V_{1,2} :
{q_0}=
\bp xz'' & yz\\ y & -z' \ep\ ,\qquad 
{q_1}= 
\bp z' & yz\\ y & -xz'' \ep\ , 
\end{equation*}
\begin{equation*}
 V_{2,2} :
{q_0}=
\bp xz' & yz\\ y & -z'' \ep\ ,\qquad 
{q_1}= 
\bp z'' & yz\\ y & -xz' \ep\ , 
\end{equation*}
\begin{equation*}
 V_{3,3} :
{q_0}=
 \bp y^2 & z'z''\\ x & -z \ep\ ,\qquad 
{q_1}= 
 \bp z & z'z'' \\ x & -y^2 \ep\ ,
\end{equation*}
\begin{equation*}
 V_{4,4} :
{q_0}=
 \bp yz & z'z''\\ x & -y \ep\ ,\qquad 
{q_1}= 
 \bp y & z'z''\\ x & -yz \ep\ , 
\end{equation*}
\begin{equation*}
 \Vb : 
{q_0}=
{\footnotesize 
\bp 
 yz & yz & -\lambda x^2z' & -zz' & 0 & 0 & 0 \\
 xz'' & -y^2 & (\lambda+1)x^2y & yz' & 0 & x^2z'' & 0 \\
 xz & xz'' & 0 & yz & 0 & \lambda x^2z & 0 \\
 0 & -xy & xz' & xz' &  0 & yz & 0 \\
 0 & 0 & xz & xz & xz'' & 0 & yz \\
 -x & 0 & y & 0 & 0 & -z & 0 \\
 0 & -x & 0 & 0 & y & z & -z' 
\ep
}\ , 
\end{equation*}
\begin{equation*}
{q_1}
{\footnotesize 
 \bp 
 y & z & 0 & 0 & 0 & x^2z'' & 0 \\
 0 & -z & z' & 0 & 0 & (\lambda-1)x^2z & 0 \\
 x & 0 & 0 & z & 0 & yz & 0 \\
 -x & 0 & y & -\lambda x^2 & 0 & 0 & 0 \\
 0 & 0 & 0 & -z & z' & 0 & yz \\
 0 & -x & 0 & y & 0 & -xz'' & 0 \\
 0 & 0 & -x & 0 & y & -xz & -xz''
 \ep 
}\ . 
\end{equation*}

\noindent
$\bullet$\ $W=(2,3,3;9)$, $A_W=(2,3,3,3)$. 
For the polynomial $f_W=x^3y+zz'z''$, $z':=z-y$, $z'':=z-\lambda y$, 
$\lambda\ne 0,1$, 
\begin{equation*}
 V_0 : 
{q_0}=
\bp zz'' & -zz'' & x^2y & 0 \\ 
    y & -z & 0 & x^2y \\ 
    x & 0 & -z & zz'' \\
    0 & x & -y & zz'' \ep\ ,\quad 
{q_1}=
 \bp z & -zz'' & x^2y & 0 \\ 
     y & -zz'' & 0 & x^2y \\
     x & 0 & -zz'' & zz'' \\
     0 & x & -y & z 
 \ep\ , 
\end{equation*}
\begin{equation*}
 V_{1,2} :
{q_0}=
\bp x^3 & z'z'' \\ z & -y \ep\ ,\qquad 
{q_1}= 
\bp y & z'z'' \\ z & -x^3 \ep\ , 
\end{equation*}
\begin{equation*}
 V_{2,3} :
{q_0}=
\bp z'z'' & x^2y\\ x & -z \ep\ ,\qquad 
{q_1}= 
\bp z & x^2y \\ x & -z'z'' \ep\ , 
\end{equation*}
\begin{equation*}
 V_{3,3} :
{q_0}=
\bp zz'' & x^2y\\ x & -z' \ep\ ,\qquad 
{q_1}= 
\bp z' & x^2y \\ x & -zz'' \ep\ , 
\end{equation*}
\begin{equation*}
 V_{4,3} :
{q_0}=
\bp zz' & x^2y\\ x & -z'' \ep\ ,\qquad 
{q_1}= 
\bp z'' & x^2y \\ x & -zz' \ep\ , 
\end{equation*}
\begin{equation*}
 \Vb : 
{q_0}=
{\footnotesize 
 \bp 
 zz'' & x^3-zz'' & \lambda zz'' & x^2 z'' & 0 & 0 & 0 \\
 -yz & z^2 & -\lambda zz'' & -x^2y & 0 & 0 & 0 \\
 -y^2 & yz & -yz+zz'' & 0 & -x^2y & 0 & 0 \\
 -xy & xz & -xz & zz'' & -zz'' & 0 & 0 \\
 -xz'' & xz'' & \lambda xz'' & 0 & 0 & x^3 & z'z'' \\
 0 & 0 & -x & y & -z & 0 & 0  \\
 0 & -x &0 & -z'' & 0 & z & -y 
 \ep
}\ ,
\end{equation*}
\begin{equation*}
{q_1}=
{\footnotesize 
 \bp 
 z & z & -\lambda z & -x^2 & 0 & x^2z & 0 \\
 y & z'' & 0 & 0 & 0 & 0 & 0 \\
 0 & -y & z & 0 & 0 & -x^2y & 0 \\
 0 & -x & 0 & z & 0 & -zz'' & 0  \\
 0 & 0 & -x & y & 0 & -zz'' & 0 \\
 0 & 0 & 0 & z'' & y & -(z'')^2 & z'z'' \\
 -x & 0 & 0 & 0 & z & 0 & -x^3 
 \ep 
}\ . 
\end{equation*}

\noindent
$\bullet$\ $W=(2,2,5;10)$, $A_W=(2,2,2,2,2)$. 
For the polynomial 
$f_W=-xyy_1y_2y_3+z^2$, $y_1:=y-x$, $y_2:=y-\lambda_1 x$, 
$y_3:=y-\lambda_2 x$, 
$\lambda_1,\lambda_2\ne 0,1$, $\lambda_1\ne\lambda_2$, 
\begin{equation*}
 V_0 : 
{q_0}={q_1}=
\bp z & -xyy_2y_3 
      & xyy_2y_3 & 0 \\ 
    y & -z & 0 & xyy_2y_3 \\ 
    x & 0 & -z & xyy_2y_3 \\
    0 & x & -y & z 
\ep\ , 
\end{equation*}
\begin{equation*}
 V_{1,2} :
{q_0}={q_1}=
\bp z & -yy_1y_2y_3 \\ x & -z \ep\ ,\quad 
 V_{2,2} :
{q_0}={q_1}=
\bp z & -xy_1y_2y_3 \\ y & -z \ep\ , 
\end{equation*}
\begin{equation*}
 V_{3,2} :
{q_0}={q_1}=
\bp z & xyy_2y_3 \\ y_1 & -z \ep\ ,  \qquad 
 V_{4,2} :
{q_0}={q_1}=
\bp z & -xyy_1y_3 \\ y_2 & -z \ep\ ,  
\end{equation*}
\begin{equation*}
 V_{5,2} :
{q_0}={q_1}=
\bp z & -xyy_1y_2 \\ y_3 & -z \ep\ ,\qquad 
\end{equation*}
\begin{equation*}
 \Vb\ : {q_0}={q_1}= 
{\small 
 \bp z & -xy_2y_3 & 0 & 0 & 0 & 0 \\
     yy_1 & -z & 0 & 0 & 0 & 0 \\
     xy_1 & 0 & -z & 0 & 0 & xy_1y_2y_3 \\
     yy_2 & 0 & 0 & z & xyy_2y_3 & 0 \\
     0 & -y_2 & 0 & -y_1 & -z & 0 \\
     0 & x & -y & 0 & 0 & z \ep
}\ . 
\end{equation*}

\noindent
$\bullet$\ $W=(2,2,3;8)$, $A_W=(2,2,2,2,3)$. 
For the polynomial $f_W=yy_1y_2y_3+xz^2$, $y_1:=y-x$, $y_2:=y-\lambda_1 x$, 
$y_3:=y-\lambda_2 x$ and $\lambda:=\lambda_1$, 
$\lambda_1\ne 0,1$, $\lambda_1\ne\lambda_2$, 
\begin{equation*}
 V_0 : 
{q_0}=
\bp 
 xz & yy_2y_3 & 
      -yy_2y_3 & 0 \\ 
    y & -z & 0 & -yy_2y_3 \\ 
    x & 0 & -z & -yy_2y_3 \\ 
    0 & x & -y & xz \ep\ , \quad 
{q_1}=
 \bp 
 z & yy_2y_3 & 
     -yy_2y_3 & 0 \\ 
     y & -xz & 0 & -yy_2y_3 \\ 
     x & 0 & -xz & -yy_2y_3 \\ 
     0 & x & -y & z 
 \ep\ , 
\end{equation*}
\begin{equation*}
 V_{1,2} :
{q_0}=
\bp xz & y_1y_2y_3 \\ y & -z  \ep\ ,\qquad 
{q_1}=
 \bp z & y_1y_2y_3 \\ y & -xz \ep\ , 
\end{equation*}
\begin{equation*}
 V_{2,2} :
{q_0}=
\bp xz & yy_2y_3 \\ y_1 & -z \ep\ ,\qquad 
{q_1}=
 \bp z & yy_2y_3 \\ y_1 & -xz \ep\ , 
\end{equation*}
\begin{equation*}
 V_{3,2} :
{q_0}=
\bp xz & yy_1y_3 \\ y_2 & -z \ep\ ,\qquad 
{q_1}=
 \bp z & yy_1y_3 \\ y_2 & -xz \ep\ , 
\end{equation*}
\begin{equation*}
 V_{4,2} :
{q_0}=
\bp xz & yy_1y_2 \\ y_3 & -z \ep\ ,\qquad 
{q_1}=
 \bp z & yy_1y_2 \\ y_3 & -xz \ep\ , 
\end{equation*}
\begin{equation*}
 V_{5,3} :
{q_0}=
\bp z^2 & y_1y_2y_3 \\ y & -x  \ep\ ,\qquad 
{q_1}=
 \bp x & y_1y_2y_3 \\ y & -z^2 \ep\ , 
\end{equation*}
\begin{equation*}
 \Vb : 
{q_0}=
{\footnotesize 
 \bp 
 xz & y_1y_2y_3 & xy_2y_3 & -xy_2y_3 & 0 & x^2y_3z & 0 \\
 -yz & z^2 & -yy_2y_3 & yy_2y_3 & 0 & -yy_3z & 0 \\
 -y^2 & yz & xz & 0 & 0 & -\lambda xyy_3 & 0 \\
 -xy & xz & 0 & xz & 0 & (y_1y_2-\lambda x^2)y_3 & 0 \\
 -x^2 & 0 & 0 & 0 & xz & (y_2-x)xy_3 & y_1y_2y_3 \\
 0 & 0 & -x & y & 0 & -z & 0 \\
 0 & -x & 0 & -x & y & 0 & -z
 \ep}\ , 
\end{equation*}
\begin{equation*}
{q_1}=
{\footnotesize 
 \bp 
 z & 0 & (x-y_2)y_3 & -\lambda xy_3 & 0 & xy_3 z & 0 \\
 y & x & 0 & 0 & 0 & 0 & 0 \\
 0 & -y & z & 0 & 0 & yy_2y_3 & 0 \\
 0 & -x & 0 & z & 0 & yy_2y_3 & 0 \\
 x & 0 & 0 & 0 & z & xy_2y_3 & y_1y_2y_3 \\
 0 & 0 & -x & y & 0 & -xz & 0 \\
 0 & 0 & 0 & -x & y & 0 & -xz
 \ep }\ . 
\end{equation*}


%

\end{document}